\documentclass{article}
\usepackage{booktabs}
\usepackage{etoolbox}
\patchcmd{\thebibliography}
  {\settowidth}
  {\setlength{\itemsep}{-3pt }\settowidth}
  {}{}
\apptocmd{\thebibliography}
  {\small}
  {}{}

\usepackage[utf8]{inputenc}
\usepackage[english]{babel}
\usepackage[T1]{fontenc}
\usepackage{lmodern}
\usepackage{mathtools}

\usepackage{xcolor}
\usepackage{caption}
\usepackage[caption=false]{subfig}

\usepackage[a4paper,top=3cm,bottom=2cm,left=3cm,right=3cm,marginparwidth=1.75cm]{geometry}
\usepackage{mathrsfs, amsthm}
\allowdisplaybreaks[4]
\usepackage{comment}
\usepackage{amsmath,amssymb}
\usepackage{graphicx}
\usepackage{bm}
\usepackage{bbm}
\usepackage{float}
\usepackage[normalem]{ulem}
\usepackage[colorlinks=true, allcolors=blue]{hyperref}
\usepackage[capitalize,nameinlink]{cleveref}

\newtheorem{theorem}{Theorem}[section]

\newtheorem{prop}[theorem]{Proposition}
\newtheorem{assp}[theorem]{Assumption}

\theoremstyle{definition}

\newtheorem{rem}[theorem]{Remark}

\numberwithin{equation}{section}

\Crefname{prop}{Proposition}{Propositions}
\Crefname{cor}{Corollary}{Corollaries}
\Crefname{assp}{Assumption}{Assumptions}
\Crefname{rem}{Remark}{Remarks}
\Crefname{defn}{Definition}{Definitions}
\Crefname{example}{Example}{Examples}

\newcommand{\E}{\mathbb{E}}
\newcommand{\N}{\mathbb{N}}

\newcommand{\RR}{\mathbb{R}}

\newcommand{\rd}{\mathrm d}
\newcommand{\op}{\text{op}}
\def\tr{\text{tr}}

\def\<{\langle}
\def\>{\rangle}

\def\nn{\nonumber}

\def\bc{\begin{center}}       \def\ec{\end{center}}
\def\ba{\begin{array}}        \def\ea{\end{array}}
\def\be{\begin{equation}}     \def\ee{\end{equation}}
\def\bea{\begin{eqnarray}}    \def\eea{\end{eqnarray}}
\def\beaa{\begin{eqnarray*}}  \def\eeaa{\end{eqnarray*}}

\usepackage{amsopn}

\newcommand{\email}[1]{\href{mailto:#1}{\nolinkurl{#1}}}
\newenvironment{keywords}{\par\medskip\noindent\textbf{Keywords.}\ }{\par}
\newenvironment{MSCcodes}{\par\smallskip\noindent\textbf{MSC codes.}\ }{\par}

\hypersetup{
pdftitle={Explicit geometric integrators for Langevin equations with superquadratic potentials},
pdfauthor={C. Chen, X. Chen, J. Hong, and M.V. Tretyakov}
}

\title{Explicit geometric integrators for Langevin equations with superquadratic potentials}
\author{Chuchu Chen\thanks{State Key Laboratory of Mathematical Sciences, Academy of Mathematics and Systems Science, Chinese Academy of Sciences, Beijing 100190, and School of Mathematical Sciences, University of Chinese Academy of Sciences, Beijing 100049, China
  (\email{chenchuchu@lsec.cc.ac.cn}, \email{chenxinyu@amss.ac.cn}, \email{hjl@lsec.cc.ac.cn}).} \and Xinyu Chen\footnotemark[1] \and Jialin Hong\footnotemark[1]
\and Michael V. Tretyakov\thanks{School of Mathematical Sciences,
University of Nottingham,
Nottingham, UK
  (\email{Michael.Tretyakov@nottingham.ac.uk}).}}
\date{\today}

\begin{document}

\maketitle

\begin{abstract}
We propose a new class of explicit methods for Langevin equations with a superquadratic potential that simultaneously preserve the conformal symplecticity and ergodicity as well as achieving high-order convergence. The construction is based on a mollified potential approach, in which the original potential is truncated to quadratic growth at infinity and globally mollified to control the superquadratic behavior while maintaining the underlying properties of the original dynamics. By appropriately choosing the cut-off radius, the far-field quadratic coefficient, and the step size, we obtain the coercivity of the mollified potential and derive its derivative estimates. These estimates allow us to prove uniform-in-time weak convergence for the proposed (including higher-order) explicit conformal symplectic numerical methods. 
Since the mollified potential is identical to the original one in a ball whose radius grows as the step size decreases, the computational cost of the new methods is essentially the same as that of the usual schemes (e.g., BAOAB) applied to the original potential. 
Numerical experiments are provided to illustrate the theoretical results.
\end{abstract}



	\section{Introduction}

    We construct convergent explicit numerical methods for (kinetic) Langevin equations that preserve \textit{both} the intrinsic geometric structure and ergodicity in the superquadratic setting.
	
    Langevin equations have been a workhorse in molecular dynamics for a few decades (see, e.g., \cite{leimkuhler_mathews_15} and references therein) and other applications (see, e.g., \cite{Coffey2012book}). They have also been successfully used for sampling \cite{Neal2011} and optimization \cite{ChenFoxGuestrin2014} tasks. 
    The challenge of those tasks is threefold. First, the related applicable problems of interest are high dimensional with computationally expensive forces, which rules out the use of implicit methods. Implicit schemes possess good stability and structural properties \cite{CDHZ2025SINUM,CHS2022, Sun2017SINUM, Talay2002}. 
    However, implementation of implicit schemes, as a rule, requires solving nonlinear equations at each time step, leading to substantial computational cost in large-scale applications. 
	Second, the SDEs need to be simulated over long periods of time, and hence geometric integrators are required. Third, potentials are usually heavily non-linear, which renders the SDEs' drift non-globally Lipschitz. In this case classical explicit numerical methods for SDEs fail to converge \cite{Hutzenthaler2011EulerDiverge} or to be ergodic \cite{Mattingly2002SPA}.
    The first two challenges were addressed by the construction of explicit quasi-symplectic/conformal symplecticity preserving numerical methods for Langevin equations (see \cite{Owhadi2010SINUM,HongWangbook2019,BenM13,MilsteinTretyakov2003Quasisym,Milstein2004book} and references therein), which are known to demonstrate excellent accuracy over long integration time inheriting those properties from deterministic and stochastic Hamiltonian systems \cite{HLW02,Hongsun2022,MRT02,MRT03}. 
    There have been a substantial number of papers addressing the third challenge within a general SDE setting via explicit tamed or balanced schemes or adaptive time stepping (see e.g.  \cite{FG19,HJ15,KL18,Mattingly2002SPA,TZ13,Handy24,Wang2025}). Although such explicit schemes converge in the case of SDEs with non-globally Lipschitz coefficients, they fail to preserve intrinsic geometric structures such as conformal symplecticity of Langevin equations.    
    The authors of \cite{Dai2025arXiv} proposed an explicit splitting scalar auxiliary variable scheme and proved strong convergence of order one and weak error bounds of order one with constants growing polynomially in time, but conformal symplecticity and ergodicity of their scheme are not established.
    A different approach to this challenge was proposed in \cite{MilsteinTretyakov2005SINUM,MilsteinTretyakov2007PhysD} (see also \cite{Milstein2004book}), where approximate trajectories leaving a sufficiently large ball are rejected,  allowing the safe use of any method of weak approximation for SDEs with non-globally Lipschitz coefficients including  explicit geometric integrators for Langevin equations. Note that under this approach schemes may not converge to the exact solution when the time step goes to zero. This approach can work well for finite time weak approximation of SDEs and for computing ergodic limits via ensemble averaging but not via time averaging. 
    The reason for this approach to work is that trajectories of regular (and moreover ergodic) SDEs live in that ball with high probability, which is also exploited here for constructing new numerical methods.  

    The above discussion leads to the interesting question of practical importance: \textit{how to construct explicit conformal symplectic schemes for Langevin equations whose numerical integration errors are uniform in time in the case of potentials growing faster than quadratically at infinity}.  
    To address this, we propose a new approach to constructing geometric integrators for Langevin equations with superquadratic potentials. Its main idea consists in appropriately mollifying the potential, which suppresses superquadratic growth while preserving the underlying conformal symplecticity and smoothness. Consider a ball whose radius $R(h)$ grows as the time step $h$ decreases. The mollified potential $U^h$ (see Section~\ref{sec_tamingAndGeneralmethod}) is identical to the given one in this ball, and then it is smoothly modified so that outside the ball of radius $2R(h)$ it equals $\lambda(h)|q|^2$ with $\lambda(h)$ proportional to $h^{-1/2}$, i.e., at infinity it does not grow faster than quadratically.
    Then (Section~\ref{sec_mainresultsofExplicitmethods}), we construct first and second order weakly convergent explicit conformal symplectic integrators for the mollified Langevin equations which form a novel class of geometric integrators for kinetic Langevin dynamics. Since $U^h=U$ in the ball of radius $R(h)$, in which trajectories live with high probability, these integrators compute only $\nabla U$ except at the rare steps when $|q|>R(h)$ (see the truncation counts in Section~\ref{sec_numericalexperiment}). Hence, their \textit{computational cost is essentially the same as that of the usual schemes} (e.g., BAOAB) applied to the original potential.

    In this paper we show that (i) the conformal symplecticity and ergodicity of the original Langevin equations are inherited by the mollified Langevin equations; (ii) their invariant measures are within $e^{-ch^{-1/2}}$ for some $c>0$ of each other; (iii) the solution of the mollified Langevin equations converges weakly to the solution of the original ones uniformly in time with error of order $h^k$ for any $k>0$ (Theorem~\ref{thm_weakconvergencetamed}); (iv) a general class of explicit methods \eqref{eq:generalP}--\eqref{eq:generalQ} applied to the mollified Langevin equations is ergodic and uniform-in-time weakly convergent of order one and two (Theorem~\ref{thm_weakconvergeconformalsym}); and also (v) to illustrate these results, we construct explicit methods of weak order one (the conformal symplectic Euler method and a Lie--Trotter splitting) and of weak order two (a conformal partitioned Runge--Kutta method, a Strang splitting and the BAOAB scheme), which simultaneously preserve  conformal symplecticity and  ergodicity (see Section~\ref{sec_mainresultsofExplicitmethods}). To the best of our knowledge, this is the first construction of convergent explicit methods that simultaneously preserve the conformal symplecticity and  ergodicity in the superquadratic setting. Note that one can also tame the potential instead of the drift, which preserves the conformal symplecticity (see Appendices~\ref{suppsec_tamed} and~\ref{suppsec_ergodicTame}). However, the tamed potential introduces a bias of order $h^\alpha$, and our analysis of numerical methods for it requires $\alpha \le 1$, so that taming the potential cannot provide provably second-order methods.  Further, schemes with the tamed potential are more expensive, since its gradient contains additional terms, including $U$ itself, which have to be computed at every time step.

	The proofs of the results of Sections~\ref{sec_tamingAndGeneralmethod} and~\ref{sec_mainresultsofExplicitmethods} are given in Section~\ref{sec_proofs} and in Appendices~\ref{suppsec_proofmollifiedequ} and~\ref{suppsec_proofofnumerical}. The technical challenge in establishing the uniform-in-time weak convergence for explicit methods in the superquadratic setting is due to the singular dependence on inverse powers of $h$ within the upper bound estimates of $U^h$ and its derivatives which is introduced by the  modification of the potential.
	Through a suitable scaling of the cut-off radius $R(h)$ and the far-field quadratic coefficient $\lambda(h)$ with respect to the step size $h$, we obtain appropriate bounds of $U^h$ and its derivatives of order $h^{-1/2}$, while maintaining a high-order approximation of $\nabla U$ by $\nabla U^h$ (see Proposition~\ref{prop:fh-fS2}).  These estimates allow us to show the coercivity of $U^h$, which compensates for the absence of a global one-sided Lipschitz condition and also mitigates the divergent $h$-dependence introduced by the modification.
	By combining the coercivity with separate estimates for small and large $|q|$, we establish the numerical Lyapunov structure and prove the uniform-in-time weak convergence results.
	Finally, numerical experiments are presented in Section~\ref{sec_numericalexperiment} to illustrate our theoretical results. 
    
\section{Mollified potential approach} \label{sec_tamingAndGeneralmethod}

\subsection{Langevin dynamics}

In this paper, we endow the space $\RR^d$ with the Euclidean norm denoted by  $|\cdot|$. For $M\in (\RR^d)^{\otimes k}$, we denote by  $\|M\|_{\text{op}}$ its operator norm induced by the Euclidean vector norm. 
We use $C$ to denote an arbitrary positive constant which may vary from one place to another.

Consider the standard Langevin equations
\begin{equation}\label{eq:origin}
		\left\{\begin{aligned}
			dP(t) &= -\upsilon P(t)dt-\nabla U(Q(t)) dt+\sigma dW(t), \, P(0)=p_0,
			\\
			dQ(t) &= P(t)dt, \, Q(0)=q_0,
		\end{aligned} \right. 
	\end{equation}
where  $\upsilon>0$ is the friction coefficient,  $\sigma >0$ is a constant, $U>0$ is the potential function	which may exhibit  superquadratic growth, and $W(t)$ is a $d$-dimensional standard Wiener process on a sufficiently rich filtered complete probability space $(\Omega, \mathcal{F},  \{\mathcal{F}_t\}_{t\geq 0}, \mathbb{P})$. 

Recall that, similar to classical dissipative Hamiltonian systems, the stochastic Langevin equation \eqref{eq:origin} possesses intrinsic geometric structure known as the conformal symplecticity (see, e.g., \cite{HongWangbook2019}), which describes the exponential decay of the symplectic form over time:
$$\rd P(t)  \wedge \rd Q(t) = e^{-\upsilon t} \rd P(0) \wedge \rd Q(0).$$
 
Throughout this work, we impose the following assumptions on the potential function $U$. For convenience, we interchangeably use the notation for the force $f:=\nabla U$.
	\begin{assp} \label{assm_CandPoly}
		Assume  $U\in \mathcal C^{\infty}(\RR^d; \RR_+)$ and all its derivatives have at most polynomial growth. Furthermore, it satisfies
		\begin{align*}
			&L_1|q|^{r_1}- C\leq U(q) \leq L_2|q|^{r_1} +C, \\ %
            &|\nabla U(q)| \leq L_3|q|^{r_2}+C,   \\
			&\|D^2U(q)\|_{\text{op}} \leq L_4 |q|^{r_3}+C,
		\end{align*}
		with constants $C>0$ and $L_i>0$, and  $r_i \ge 0$ satisfying $r_1 > 2$, $r_2=r_1-1$, and $r_3\leq r_2-1$. 
	\end{assp}

Let us discuss this assumption. It describes a superquadratic potential belonging to the class of practical interest
for which explicit numerical methods require stabilisation. Superquadratic growth 
$r_{1}>2$ is the regime we are interested in. Indeed,
when $r_{1} = 2,$ the gradient $\nabla U$ grows not faster than linearly, the
standard explicit schemes are stable, and no mollification is needed. The
motivation for this work is the superquadratic
case, where $\nabla U$ grows superlinearly and, as is well known, classical
explicit schemes for the associated Langevin equations fail to converge (see
e.g. \cite{Mattingly2002SPA}). We emphasise that Assumption~\ref{assm_CandPoly}
relaxes the global Lipschitz condition on $\nabla U$  to polynomial
growth. The exponent relations $r_{2}=r_{1}-1$ and $r_{3} \leq r_{2}-1$ reflect
differentiation as each derivative of a pure power $|q|^{r_{1}}$ lowers the
growth exponent by one. Typical examples satisfying
Assumption~\ref{assm_CandPoly} include the Ginzburg--Landau (or $\phi^{4}$)
potential $(|q|^{2}-1)^{2}+1$, its anisotropic generalisations, and polynomial potentials 
with positive-definite leading part.

As usual, we also need a dissipativity condition on the drift.
    
	\begin{assp}\label{assm_ErgodicU}
		There exist constants  $c_0, c_1>0$ such that
			$	\nabla U(q)^\top q \geq c_0 U(q)-c_1.$
	\end{assp}

Under Assumptions~\ref{assm_CandPoly}--\ref{assm_ErgodicU}, the solution ${(P(t),Q(t))}$, $t\ge 0$, 
 of \eqref{eq:origin} is ergodic: for smooth test functions $g\in\mathcal{C}^\infty(\RR^{2d};\RR)$ with derivatives of all orders growing not faster than polynomials
	$$
	\lim_{T\to\infty} \frac 1T \int^T_0 \E [ g(P(t),Q(t))] dt = \int_{\RR^{2d} }g(p,q) \pi (dp,dq),
	$$
    where $\pi$ has density proportional to
$\exp(-\frac{2\upsilon}{\sigma^2}(\frac{1}{2}|p|^2 + U(q)))$ and is the unique Gibbs invariant measure (see e.g. \cite{Mattingly2002SPA,Talay2002}).
 Assumption~\ref{assm_ErgodicU} is the standard condition in the analysis of ergodic Langevin dynamics beyond the convex setting and is implied by convexity at infinity (see e.g. \cite{Mattingly2002SPA}).
 We note that Assumption~\ref{assm_ErgodicU} is consistent with and complementary to Assumption~\ref{assm_CandPoly}. For a power potential $U(q)\sim |q|^{r_{1}}$, one has $\nabla U(q)^{\top }q\sim r_{1}U(q),$ and
Assumption~\ref{assm_ErgodicU} holds with $c_{0}$ up to order $r_{1}$. 
The two assumptions are therefore satisfied simultaneously by the examples above.

\subsection{Mollified Langevin dynamics}

    Recall our discussion from the Introduction about  our mollified potential approach inspired by \cite{MilsteinTretyakov2005SINUM}. The solution of \eqref{eq:origin} lives in a ball of radius $R$ with high probability \cite{HAS80,MilsteinTretyakov2005SINUM}, then we keep the mollified potential $U^h$ equal to $U$ in that ball. For explicit schemes to work, $U^h$ should grow not faster than quadratically at infinity, say for $|q| \ge 2R$, and we need a smooth transition between the two regimes. For convergence of the mollified Langevin to the original one and for convergence of explicit numerical methods, we need $R$ to appropriately grow with decrease of the time step $h$ of numerical schemes and be compatible with the behaviour of $U^h$ for large $|q|$.  
	
    Keeping the above in mind, the \textit{mollified potential} $U^h$ under this approach is defined by 
	\begin{align}\label{eq:molpot} 
		U^h(q):=U(q) \theta\Big(\frac{|q|}{R(h)}\Big) + \Big[1-\theta\Big(\frac{|q|}{R(h)}\Big) \Big] \lambda(h) |q|^2,
	\end{align}
	where $h\in(0,1)$, 
    $$R(h)=R_0h^{-\frac{1}{2(r_2-1)}}$$ 
    with  a fixed reference length $R_0 \ge 1$, 
    $$\lambda(h) =(2R_0)^{r_2-1} L_2 h^{-1/2},$$ and  
    $\theta: [0,\infty) \to [0,1]$ is a smooth function such that $\theta(x) =1$ for $x\in [0,1]$, $\theta(x) =0 $ for $x \geq 2$, and $\theta'$ is bounded  satisfying $\theta'(x) \leq 0$ for all $x \ge 0$.  Here,  we denote by $\theta'$ and $\theta''$ the first-order and second-order derivatives of $\theta$, respectively. 
    
 The fixed reference length $R_0$ is required in the definition of $R(h)$ for dimensional consistency so that $R(h)$ is a length measured in the same unit as the norm $|q|$ of the position. It also gives sufficient flexibility to ensure that the process $Q(t)$ spends all but a negligible fraction of its time  in the ball of radius $R(h)$ even for step sizes $h$ that are not especially small.
 Since the mollified force $f^{h}=\nabla U$ on that ball, this is exactly the aim we set out for the construction of the mollified potential: to disturb the Langevin dynamics as little as possible, which here means that $f^{h}(Q(t))=\nabla U(Q(t))$ with high probability at every $t$.
  The theoretical results of this and the next section hold for every $R_0 \ge 1$. In practice $R_0$ is bounded above by the stability of the numerical scheme, while it is advantageous to take it as large as possible so that ${|q| \le R(h)}$ carries all but a negligible fraction of the Gibbsian invariant mass. Since that fraction decays super-exponentially in $R_0$ (see the proof of Theorem~\ref{thm_weakconvergencetamed}(\romannumeral 2)), fine tuning is not needed, it suffices to verify a posteriori that the proportion of steps at which $Q^h_k$ leaves the ball of radius $R(h)$ is small, say below $10^{-3}$. See also our numerical experiments in Section~\ref{sec_numericalexperiment}. 
 
	Let us provide a simple example of $\theta$ which satisfies the above conditions:
    \begin{align}\label{eq:theta}
        \theta(x)= \left\{
        \begin{aligned} 
        &1, \quad &x \in [0,1],\\
        &\frac{g(2-x)}{g(x-1)+g(2-x)},\quad &x\in(1,2),\\
        & 0, \quad & x \ge 2,
        \end{aligned}
        \right. 
    \end{align} 
    where   $g(t)=e^{-1/t}$ for $t>0$ and $g(t)=0$ otherwise.

The mollified force $f^h=\nabla U^h$ takes the form
	\begin{align} \label{eq:fhS2}
		f^h(q)&= \theta\Big(\frac{|q|}{R(h)}\Big) \nabla U(q)   + \Big[1-\theta\Big(\frac{|q|}{R(h)}\Big) \Big]2\lambda(h) q + \big(U(q)- \lambda(h)  |q|^2 \big) \nabla \Big( \theta\Big(\frac{|q|}{R(h)}\Big)\Big) \nn \\
		&= \theta\Big(\frac{|q|}{R(h)}\Big) \nabla U(q)   + \Big[1-\theta\Big(\frac{|q|}{R(h)}\Big) \Big]2\lambda(h)  q +  \frac{(U(q)- \lambda(h)  |q|^2 ) q}{R(h)|q|} \theta'\Big(\frac{|q|}{R(h)}\Big).
	\end{align}
	
    The proof of the following proposition on estimates for the mollified force is given in Section~\ref{sec_proofpart1}. 
	
	\begin{prop}\label{prop:fh-fS2}
		Under Assumption~\ref{assm_CandPoly}, there is a constant $C>0$ independent of $h$ and $q$ (and, in (i) and (iv), of $k$) such that \\
		(\romannumeral 1)
		$
		| f(q)-f^h(q)| \leq  Ch^k(1+|q|^{(2k+1)r_2 -2k}) 
		$ for any $k>0$; \\ %
		(\romannumeral 2)
		$
		|f^h(q)| \leq C\min\{ 1+|q|^{r_2},  h^{-1/2}(1+|q|)\}
		$; \\
		(\romannumeral  3)
		$\|D^2 U^h(q)\|_{\op} \leq C\min\{1+|q|^{r_2-1},  h^{-1/2}\}$;  \\
		(\romannumeral  4) 
        $\| D^2 U(q)-D^2U^h(q)\|_{\op} \leq  Ch^k(1+|q|^{(2k+1)(r_2 -1)})  $ for any $k>0$; \\
        (\romannumeral  5)
  $\|D^3U^h(q)\|_{\op}\leq C(1+|q|^{\bar m})$ for some $\bar m\geq 0$ independent of $h$.
	\end{prop}
	Since $\theta$ is smooth, $U^h$ inherits the smoothness of $U$. Together with Proposition~\ref{prop:fh-fS2} (\romannumeral 3), this ensures that $f^h$ is globally Lipschitz with the constant $Ch^{-1/2}$. 
	
	The next proposition shows that the coercivity condition can be inherited by the proposed mollified potential (see its proof in Section~\ref{sec_proofpart1}). This property is important in the subsequent analysis of numerical ergodicity.
	
	\begin{prop} \label{prop_ergodicfhS2} 
		Under Assumptions \ref{assm_CandPoly} and \ref{assm_ErgodicU}, the inequality holds
		\begin{align*}
			f^h(q)^\top q \geq c_2 U^h(q)-c_3,
		\end{align*}
        where $c_2 =\min\{c_0 , 2\}$ and $c_3>0$ is independent of $h$.
	\end{prop}

	Having introduced the mollified potential $U^h$, we write the \textit{mollified Langevin equations}
	\begin{equation}\label{eq:tamed}
		\left\{\begin{aligned}
			dP^h(t) &= -\upsilon P^h(t)dt-\nabla U^h(Q^h(t)) dt+\sigma dW(t), \,\, P^h(0)=p_0,
			\\
			dQ^h(t) &= P^h(t)dt, \,\, Q^h(0)=q_0.
		\end{aligned} \right.
	\end{equation}
    
	Using the properties derived in Propositions~\ref{prop:fh-fS2}--\ref{prop_ergodicfhS2}, 
	the conformal symplecticity and ergodicity of the original Langevin equation (\ref{eq:origin}) are inherited 
    by the mollified Langevin dynamics \eqref{eq:tamed}. The invariant measures of (\ref{eq:origin}) and \eqref{eq:tamed} 
    are within $e^{-cR^{r_1}(h)}$ of each other. Further, the solution of \eqref{eq:tamed}  converges weakly to the solution of \eqref{eq:origin} uniformly in time. These results are stated in the below theorem (see its proof in Section~\ref{sec:proffThem21}).
	\begin{theorem} \label{thm_weakconvergencetamed}
		(\romannumeral 1) 
		The phase flow of \eqref{eq:tamed} preserves the conformal symplecticity, i.e., 
		$\rd P^h(t)  \wedge \rd Q^h(t) = e^{-\upsilon t} \rd P(0) \wedge \rd Q(0).$
		\\
        (\romannumeral 2) Let Assumptions \ref{assm_CandPoly} and \ref{assm_ErgodicU} hold. Then, for all sufficiently small $h>0$, the process $\{(P^h(t), Q^h(t))\}_{t\geq 0}$ is ergodic with the unique invariant measure $\pi_h$ being Gibbsian whose density is proportional to $\exp\big(-\frac{2\upsilon}{\sigma^2}\big(\frac12|p|^2+U^h(q)\big)\big)$.
		Furthermore, for any measurable $g:\RR^{2d}\to\RR$ with $|g(p,q)|\leq K(1+|p|^m+|q|^m)$ for some constants $K,m \geq 0$ and any $c\in(0,2\upsilon L_1/\sigma^2)$, there is a constant $C>0$ independent of $h$ such that
		\begin{align}
			|\pi(g)-\pi_h(g)| \leq C\exp\big(-cR^{r_1}(h)\big)=C\exp\Big(-cR_0^{r_1}h^{-\frac{r_1}{2(r_2-1)}}\Big).
		\end{align}
		In particular, for any $k>0$ there is a constant $C(k)>0$ independent of $h$ such that $|\pi(g)-\pi_h(g)|\leq C(k)h^k$.
		\\
		(\romannumeral 3) Let Assumptions \ref{assm_CandPoly} and \ref{assm_ErgodicU} hold. For any given $k>0$ and any  $g\in \mathcal C^{\infty}(\RR^{2d}; \RR)$ with derivatives growing no faster than polynomials at infinity, there exists a constant $C(k)>0$ independent of $h$ (but dependent on the initial value $(p_0,q_0)$) such that for sufficiently small $h>0$:
		\begin{align}
			\sup_{t\geq0}\big|\E[g(P^h(t),Q^h(t))]-\E[g(P(t),Q(t))]\big| \leq C(k)h^k.
		\end{align}
	\end{theorem}

\section{Explicit methods preserving both the conformal symplecticity and ergodicity} \label{sec_mainresultsofExplicitmethods}

We start this section with a generic construction of an approximation for the mollified Langevin equations (\ref{eq:tamed}) (Section~\ref{sec:genmeth}) and its convergence proof. It is an auxiliary construction (not to be used in practice) and an auxiliary result convenient for straightforward proving convergence order in terms of ergodic limits (we limit it to first and second orders, it can potentially be extended to order 3) for any numerical scheme satisfying certain conditions.
Then, we provide an illustration of conformal (hence quasi-symplectic) schemes of order 1 and 2 for (\ref{eq:tamed}) (and hence for (\ref{eq:origin})), we recall that a scheme for Langevin equations (\ref{eq:tamed}) (and hence for (\ref{eq:origin})) is called conformal if at every step it satisfies the following property  \cite{HongWangbook2019}:
\[
  \sum\nolimits_i \mathrm dP^i_{n+1}\wedge \mathrm dQ^i_{n+1}
  \;=\; e^{-\upsilon h}\sum\nolimits_i \mathrm dP^i_{n}\wedge \mathrm dQ^i_{n}.
\]

In Section~\ref{sec:meth1}, we consider two examples of conformal first-order schemes.
In Sections~\ref{sec:meth2a} and~\ref{sec:meth2b}, three types of 2nd order conformal methods are considered. 
In Section~\ref{sec:meth2a}, we get a conformal partitioned Runge–Kutta method obtained from the friction-conjugated stochastic Hamiltonian system. 
In Section~\ref{sec:meth2b},  we consider a scheme based on the  Strang splitting with the OU +
free-motion flow solved exactly and a scheme based on splitting into three subflows, the Ornstein--Uhlenbeck process and the two subflows of the deterministic Hamiltonian system which is the family usually referred to as BAOAB and its permutations \cite{leimkuhler_mathews_15}. For brevity, we discuss only BAOAB, the others (ABOBA, OBABO, etc.) can be treated similarly.

\subsection{A generic approximation and convergence theorem}\label{sec:genmeth}

We define $t_n:=nh$ and $\Delta W_n:=W(t_{n+1})-W(t_n)$.
In what follows we consider one-step maps $(P^h_{n+1},Q^h_{n+1})=\Psi_h\bigl(P^h_n,Q^h_n;\,W(t_n+\cdot)-W(t_n)\bigr)$ whose increments admit the representation
	\begin{align} \label{eq:generalP}
		P^h_{n+1}
		&=P^h_n+ (-\upsilon P^h_n  - f^h(Q^h_n) ) h +\sigma \Delta W_n + (a_1 P^h_n +a_2 f^h(Q^h_n) +a_3 )h^2 +\mathcal{Z}_{p, n,1}  \nn\\
		&\quad +\mathcal{R}_{p,n,1}+ h \bar H_{p,n} (a_{0,p,1}  h P^h_n +a_{0,p,2} h^2 P^h_n +a_{0,q}  h^2 f^h(Q^h_n)+ \mathcal{Z}_{p, n,2}+\mathcal{R}_{p,n,2}) \nn \\
		&= \big(I_d +[-\upsilon I_d+a_{0,p,1} h\bar H_{p,n}] h +[a_1 I_d+ a_{0,p,2} h\bar H_{p,n}]h^2 \big)P^h_n \nn \\
		&\quad + (- h I_d +[a_2 I_d +a_{0,q} h\bar H_{p,n}]h^2 ) f^h(Q^h_n) + \sigma \Delta W_n + a_3 h^2 %
		+ \mathcal{Z}_{p, n,1} \nn \\
		&\quad +h\bar H_{p,n} \mathcal{Z}_{p, n,2}+\mathcal{R}_{p,n} 
	\end{align}
	and
	\begin{align}\label{eq:generalQ}
		Q^h_{n+1} &=Q^h_n+ hP^h_{n} + (b_1 P^h_n +b_2 f^h(Q^h_n) +b_3)h^2  + \mathcal{Z}_{q, n,1} +\mathcal{R}_{q,n,1}\nn\\
		&\quad + h^2 \bar H_{q,n} ( b_{0,p} h P^h_n +\mathcal{Z}_{q, n,2}
		+ \mathcal{R}_{q,n,2})  \nn \\
		&= (h I_d+[b_1 I_d +b_{0,p} h\bar H_{q,n}] h^2) P^h_n  + Q^h_n + b_2 h^2 f^h(Q^h_n) +  b_3 h^2
		+ \mathcal{Z}_{q, n,1}\nn \\
		&\quad + h^2 \bar H_{q,n} \mathcal{Z}_{q, n,2} +\mathcal{R}_{q,n}. 
	\end{align}
	Here, $a_i, b_i, a_{0,p,i}, a_{0,q}$ and $b_{0,p}$ with $i=1,2$ are real constants; $a_3$ and $b_3$ are constant vectors in $\RR^d$;
    $\bar H_{j,n}:=\int_0^1 D^2U^h\big(Q^h_n+s(\mathcal{Q}_{j,n}-Q^h_n)\big)\,d s,\quad j=p,q$ is the averaged Hessian along the segment from $Q^h_n$ to an $\mathcal{F}_{t_{n+1}}$-measurable point $\mathcal{Q}_{j,n}, j=p,q$ with $|\mathcal{Q}_{j,n}-Q^h_n|\leq|Q^h_{n+1}-Q^h_n|+h|P^h_n|$ so that $f^h(\mathcal{Q}_{j,n})-f^h(Q^h_n)=\bar H_{j,n}(\mathcal{Q}_{j,n}-Q^h_n)$ (in the schemes of this section, $\mathcal{Q}_{j,n}$ is the point at which the scheme evaluates the force after $Q^h_n$);
	$\mathcal{Z}_{p, n,i}$ and $ \mathcal{Z}_{q, n,i}$, $i=1,2$, are stochastic convolutions (It\^o integrals with deterministic integrands, hence Gaussian) that possess finite moments of any positive order  for sufficiently small $h>0$, satisfying
	\begin{align*}
		&\E[ \mathcal{Z}_{p, n,i} |\mathcal{F}_{t_n}]=0, \quad \E[ \mathcal{Z}_{q, n,i} |\mathcal{F}_{t_n}]=0, \\ 
		&\E[ |\mathcal{Z}_{p, n,i}|^{\alpha} |\mathcal{F}_{t_n}]\leq Ch^{3\alpha/2}, \quad \E[ |\mathcal{Z}_{q, n,i}|^{\alpha} |\mathcal{F}_{t_n}]\leq Ch^{3\alpha/2}
	\end{align*}
	with $\alpha\geq 2$;
	$\mathcal{R}_{p,n,i}$ and $\mathcal{R}_{q,n,i}$ with $i\in\{1,2\}$  are higher-order remainders satisfying
	$
		|\mathcal{R}_{p,n,i}| \leq Ch^3(1+|P^h_n|+|f^h(Q^h_n)| ), \, |\mathcal{R}_{q,n,1}| \leq Ch^3(1+|P^h_n|+|f^h(Q^h_n)| ), 
	$
	and
	$
		|\mathcal{R}_{q,n,2}| \leq Ch^2(1+|P^h_n|+|f^h(Q^h_n)| ).
	$
	
	Using  Proposition~\ref{prop:fh-fS2} (\romannumeral 3), the remainders
	$\mathcal{R}_{p,n}:= \mathcal{R}_{p,n,1}+h \bar H_{p,n} \mathcal{R}_{p,n,2} $ and $\mathcal{R}_{q,n}:=\mathcal{R}_{q,n,1}+h^2 \bar H_{q,n} \mathcal{R}_{q,n,2}$ satisfy 
	$
		|\mathcal{R}_{p,n}| \leq Ch^3(1+|P^h_n|+|f^h(Q^h_n)| )$ and $|\mathcal{R}_{q,n}|\leq Ch^3(1+|P^h_n|+|f^h(Q^h_n)| ).
	$

    Finally, we require that the noise enters schemes non-degenerately. 

\begin{assp}\label{assp:nondeg}
 There is $m\in\{1,2\}$ such that for every starting point $x=((P^h_0)^\top,(Q^h_0)^\top)^\top\in\RR^{2d}$, we have $((P^h_m)^\top,(Q^h_m)^\top)^\top=\Gamma(x,\xi)$, where $\xi$ is a non-degenerate Gaussian random vector in $\RR^{2d}$ formed from the noise of the first $m$ steps, $\Gamma$ is continuously differentiable, and, for each $x$, $\Gamma_x:=\Gamma(x,\cdot)$ is a bijection of $\RR^{2d}$ with nonsingular Jacobian.   
\end{assp}

    This assumption is satisfied by all the schemes of Sections~\ref{sec:meth1}--\ref{sec:meth2b} and is verified in their proofs.
For example, for the conformal symplectic Euler method \eqref{eq:ConformalSym}, we have $m=2$ and $\xi=(\Delta W_0,\Delta W_1)$.  

We are now positioned to state the main result concerning the ergodicity and weak convergence order of the general numerical methods \eqref{eq:generalP}--\eqref{eq:generalQ}.
	\begin{theorem} \label{thm_weakconvergeconformalsym}
		Let Assumptions \ref{assm_CandPoly}, \ref{assm_ErgodicU} and \ref{assp:nondeg} hold. Then for sufficiently small $h>0$, the chain $\{(P^h_n, Q^h_n)\}_{n\in \N}$ is ergodic with a unique invariant measure $\pi_{h, \text{num}}$.  Furthermore, let the scheme start from $(P^h_0,Q^h_0)=(P(0),Q(0))$ and let $g\in\mathcal C^{\infty}(\RR^{2d};\RR)$ have all its derivatives growing not faster than polynomials at infinity. Then the following weak convergence results hold uniformly in time: \\
		(\romannumeral 1)
		There exists a constant $C>0$ independent of $h$ such that
		\begin{align*}
        \sup_{n\geq0}\big|\E[g(P^h_n,Q^h_n)]-\E[g(P(t_n),Q(t_n))]\big|&\leq Ch,\\
			|\pi(g)-\pi_{h,\text{num}} (g)| & \leq Ch.
		\end{align*}
		(\romannumeral 2)  Let the coefficients in  \eqref{eq:generalP}--\eqref{eq:generalQ} satisfy $a_{0,p,1}=-\frac 12, a_1=\frac{\upsilon^2}{2}, a_2=\frac{\upsilon}{2},  a_3=0$,  $b_1= -\frac \upsilon 2,  b_2=-\frac 12, b_3=0$, and the stochastic convolutions satisfy $|\E[(\mathcal{Z}_{p,n,1}+\upsilon \sigma \Delta B_n ) (\Delta W_n )^\top|\mathcal F_{t_n} ]|=\mathcal{O} (h^3)$ and  $|\E[ (\mathcal{Z}_{q,n,1}- \sigma \Delta B_n )(\Delta W_n )^\top | \mathcal{F}_{t_n} ]|=\mathcal{O}(h^{3})$ with $\Delta B_n:=\int^{t_{n+1}}_{t_n} \int^s_{t_n} dW(r)ds$.
		Then there exists a constant $C>0$ independent of $h$ such that
		\begin{align*}
        \sup_{n\geq0}\big|\E[g(P^h_n,Q^h_n)]-\E[g(P(t_n),Q(t_n))]\big|&\leq Ch^2,\\
			|\pi(g)-\pi_{h,\text{num}} (g)|  & \leq Ch^2.
		\end{align*}
	\end{theorem}
	The proof of Theorem~\ref{thm_weakconvergeconformalsym} is postponed to Section~\ref{sec_proofpart2}. Moreover, we would like to emphasize that the proposed approach could also be adapted to construct even higher-order explicit numerical methods that preserve both the conformal symplecticity and the ergodicity. %
    Note that in parts (\romannumeral 1) and (\romannumeral 2) of the above theorem the constants $C$ depend on the initial value.
	
    These bounds also control the distance to the ergodic limit at every finite time. Let $g$ be as in Theorem~\ref{thm_weakconvergeconformalsym}. Since $\pi$ is invariant for \eqref{eq:origin}, $\pi(g)=\int_{\RR^{2d}}\E[g(X_{0,y}(t))]\,\pi(\rd y)$ for every $t\geq0$, where $X_{0,y}$ denotes the solution of \eqref{eq:origin} started at $y$. Hence the derivative estimate \eqref{eq:Dku} with $l=1$, established in Section~\ref{sec:proffThem21}, gives $|\E[g(P(t),Q(t))]-\pi(g)|\leq C(1+|P(0)|^{\varrho}+|Q(0)|^{\varrho})e^{-\gamma t}$ for some $C,\varrho, \gamma>0$, because $\pi$ has moments of all orders. Together with Theorem~\ref{thm_weakconvergeconformalsym} (\romannumeral 1), this shows that there is a function $C(\cdot,\cdot)$ of polynomial growth, independent of $h$ and $n$, such that for sufficiently small $h>0$
	\begin{align} \label{equ_ensemble1}
		\big|\E[ g(P^h_n, Q^h_n)- \pi(g)]\big| \leq C( P(0),Q(0))(h+ e^{-\gamma nh}), \quad \text{ for any }\, n\in \N^+.
	\end{align}
	Similarly, 
	under the conditions of Theorem~\ref{thm_weakconvergeconformalsym} (\romannumeral 2)
	we have  
	\begin{align} \label{equ_ensemble2}
		\big|\E[ g(P^h_n, Q^h_n)- \pi(g)]\big| \leq C( P(0),Q(0))(h^2+ e^{-\gamma nh}), \quad \text{ for any }\, n\in \N^+.
	\end{align}
	
    The chain $\{(P^h_n, Q^h_n)\}_{n\in \N}$ also satisfies the strong law of large numbers. For sufficiently small $h>0$, the Lyapunov condition established in
the proof of Proposition~\ref{prop_momentboundedgeneral} (see \eqref{equ_numericalboundClaim}) and the positive transition density obtained in the proof of Theorem~\ref{thm_weakconvergeconformalsym} make $\{(P^h_n,Q^h_n)\}_{n\in\N}$ a positive Harris recurrent chain so that for any initial value $(P(0)^\top,Q(0)^\top)^\top\in\RR^{2d}$ and any measurable $g:\RR^{2d}\to\RR$ of at most polynomial growth, we have \cite[Chapter~17]{MeynTweedie2009}:
$$\lim_{N\to\infty}\frac1N\sum_{n=1}^{N}g(P^h_n,Q^h_n)=\pi_{h,\text{num}}(g), \quad \text{a.s.}$$  
For $g$ as in Theorem~\ref{thm_weakconvergeconformalsym}, the theorem then gives
	\begin{align} \label{equ_LLN}
		\lim_{h\to 0}	\lim_{N\to \infty}	\frac 1N \sum^N_{n=1} g(P^h_n, Q^h_n )=\lim_{h\to 0}	 \pi_{h,\text{num}}(g)= \pi(g), \quad \text{a.s.}
	\end{align}
	Theorem~\ref{thm_weakconvergeconformalsym} and equations \eqref{equ_ensemble1}--\eqref{equ_LLN} guarantee that 
	the ergodic limit can be approximated via both the ensemble-average estimator $\hat{g}^{\text{erg}}_T:= \frac 1M \sum^M_{j=1} g(P^{h,j}_N, Q^{h,j}_N)$ and the time-average estimator $\bar{g}^{\text{erg}}_{\tilde T}:=\frac 1K \sum^K_{n=1}g(P^h_n, Q^h_n)$. Here, $M$ denotes the number of independent sample paths for the Monte–Carlo simulation, 	$(P^{h,j}_N, Q^{h,j}_N)$ denotes
	the approximation of the $j$-th trajectory at the final time step $N$, $T=hN$, $\tilde T =hK$. Typically, $K \gg N$ or equivalently $\tilde T \gg T$.
	The corresponding numerical experiments are presented in Section~\ref{sec_numericalexperiment}.

\subsection{First order schemes}\label{sec:meth1}

Based on the general explicit numerical method and the weak convergence Theorem~\ref{thm_weakconvergeconformalsym}, we construct explicit methods of weak order 1 that simultaneously preserve the conformal symplecticity and the ergodicity in the following proposition, whose proof is provided in Section~\ref{subsec_proofconcreteexamples}.
\begin{prop}\label{prop_example1st}
    Let Assumptions \ref{assm_CandPoly} and \ref{assm_ErgodicU} hold. Consider the numerical methods as follows 

    (\romannumeral 1)
    Applying the conformal symplectic Euler method to  \eqref{eq:tamed} yields the following numerical method 
		\begin{equation} \label{eq:ConformalSym}
			\left\{\begin{aligned}
				P^h_{n+1} &= e^{-\upsilon h} (P^h_n-hf^h(Q^h_n)+\sigma \Delta W_n),
				\\
				Q^h_{n+1} &= Q^h_n+hP^h_{n+1},
			\end{aligned} \right.
		\end{equation}
		where $\Delta W_n :=W(t_{n+1})-W(t_n)$ %
		and $t_{n+1}-t_n=h$. 
    
   (\romannumeral 2) 
   	Consider the solution obtained by Lie--Trotter splitting $(P^{h, \text{LT}}_{n+1}, Q^{h, \text{LT}}_{n+1})  := \Phi^S_{t_n,t_{n+1}} \circ \Phi^D_{t_n, t_{n+1}} (P^{h, \text{LT}}_{n}, Q^{h, \text{LT}}_n)$, where $ \Phi^D_{t_1, t_2}(p,q)$ and $ \Phi^S_{t_1, t_2}(p,q)$ are  the flows of the  nonlinear Hamiltonian system and  the linear stochastic differential equation
\begin{align}
D:\quad & dP=-f^h(Q)\,dt, &&dQ=0,
\label{eq:splitting1}\\
S: \quad & dP=-\upsilon P\,dt+\sigma\,dW,&&dQ=P\,dt,
\label{eq:splitting2}
\end{align}
for $t\in [t_1,t_{2}]$ with $0 \leq t_1<t_2$.

    Then for sufficiently small $h>0$, numerical methods in (\romannumeral 1) and (\romannumeral 2) each preserve the conformal symplecticity and admit a unique invariant measure. 
    Furthermore, (\romannumeral 1) and (\romannumeral 2) converge weakly with order 1 uniformly in time in the sense of Theorem~\ref{thm_weakconvergeconformalsym} (\romannumeral 1).
\end{prop}

\subsection{Second order schemes based on splitting: stochastic Hamiltonian system + friction ODE}\label{sec:meth2a}

We first consider the first class of second-order schemes described above. By applying the general weak convergence result in Theorem~\ref{thm_weakconvergeconformalsym}, the following proposition shows that these explicit methods simultaneously preserve conformal symplecticity and ergodicity, and achieve second-order weak convergence uniformly in time. 
The conformal PRK method is taken as an example in the following proposition, whose proof is provided in Section~\ref{subsec_proofconcreteexamples}.
\begin{prop} \label{prop_exampleSplitting}
     Let Assumptions \ref{assm_CandPoly} and \ref{assm_ErgodicU} hold. Consider the following conformal PRK method 
		\begin{equation}\label{eq:splittingPRK}
			\left\{\begin{aligned}
				{Q}^{h,\text{PRK}}_{n+\frac 12} &= Q^{h,\text{PRK}}_n +\frac h2 e^{-\frac{\upsilon h}{2}} P^{h,\text{PRK}}_n, \\
				P^{h,\text{PRK}}_{n+1} &= e^{-\upsilon h} P^{h,\text{PRK}}_n - he^{-\frac{\upsilon h}{2}} f^h(Q^{h,\text{PRK}}_{n+\frac 12}) +\sigma \Delta W_n -\upsilon \sigma \Delta B_n, \\
				Q^{h,\text{PRK}}_{n+1} &= {Q}^{h,\text{PRK}}_{n+\frac 12}+\frac h2 e^{\frac{\upsilon h}{2}} P^{h,\text{PRK}}_{n+1} + \sigma \Delta B_n - \frac{\sigma h}{2} \Delta W_n,
			\end{aligned} \right.
		\end{equation}
		where $\Delta B_n := \int^{t_{n+1}}_{t_n} \int^s_{t_n} d W(r) ds$. 
        Then for sufficiently small $h>0$, this method preserves the conformal symplecticity and admits a unique invariant measure. 
        Furthermore, it converges weakly with order 2 uniformly in time in the sense of Theorem~\ref{thm_weakconvergeconformalsym} (\romannumeral 2).
\end{prop}

\subsection{Second order Ornstein–Uhlenbeck-based splitting schemes}\label{sec:meth2b}
Examples from the second class of second-order schemes are presented in the following proposition. The proof is provided in Section~\ref{subsec_proofconcreteexamples}.
\begin{prop} \label{prop_exampleOU}
Let Assumptions \ref{assm_CandPoly} and \ref{assm_ErgodicU} hold. Consider the numerical methods as follows 

(\romannumeral 1) Consider the solution based on the Strang splitting technique $(P^{h, \text{S}}_{n+1}, Q^{h, \text{S}}_{n+1}) := \Phi^D_{t_{n+1/2}, t_{n+1}}  \circ \Phi^S_{t_n,t_{n+1}}  \circ \Phi^D_{t_n, t_{n+1/2}} (P^{h, \text{S}}_{n}, Q^{h, \text{S}}_n)$, where $ \Phi^D_{t_1, t_2}(p,q)$  and
		$ \Phi^S_{t_1, t_2}(p,q)$ are the flows defined in \eqref{eq:splitting1}--\eqref{eq:splitting2}.

(\romannumeral 2) 
Consider the solution using the BAOAB splitting technique
$(P^{h, \text{S2}}_{n+1}, Q^{h, \text{S2}}_{n+1}) := \Phi^B_{t_{n+1/2}, t_{n+1}}  \circ \Phi^A_{t_{n+1/2}, t_{n+1}} \circ \Phi^O_{t_n,t_{n+1}}  \circ \Phi^A_{t_n, t_{n+1/2}}\circ \Phi^B_{t_n, t_{n+1/2}} (P^{h, \text{S2}}_{n}, Q^{h, \text{S2}}_n)$, where 
  $\Phi^A$, $\Phi^B$, and $\Phi^O$ are the flows of the free motion, force impulse, and Ornstein–Uhlenbeck subsystems, respectively:
\begin{align*}
A:\quad &dP=0, &&dQ=P\,dt,\\
B:\quad &dP=-f^h(Q)\,dt, &&dQ=0,\\
O:\quad &dP=-\upsilon P\,dt+\sigma\,dW, &&dQ=0.
\end{align*}

 Then for sufficiently small $h>0$, numerical methods in (\romannumeral 1) and (\romannumeral 2) each preserve the conformal symplecticity and admit a unique invariant measure. 
 Furthermore, (\romannumeral 1) and (\romannumeral 2) converge weakly with order 2 uniformly in time in the sense of Theorem~\ref{thm_weakconvergeconformalsym} (\romannumeral 2).
\end{prop}

\section{Proofs}\label{sec_proofs}
   
 \subsection{Proof of Theorem~\ref{thm_weakconvergencetamed}} \label{sec:proffThem21}  
  
	We first define the  Lyapunov function
	\begin{align} \label{equ_defVh}
		V^h(p,q) := \frac{1}{2} |p|^2 + U^h(q)+\upsilon \alpha_1 p^\top q
	\end{align}
	with a constant $ \alpha_1\in(0, 1)$.  With \eqref{equ_defVh}, we establish uniform-in-time moment bounds for $\{(P^h(t), Q^h(t))\}_{t\geq 0}$ in the following proposition, whose proof is provided in Section~\ref{subsec_proofMomentequ}.	
	\begin{prop} \label{prop_MomentboundedOfTamedequ}
		Let Assumptions \ref{assm_CandPoly} and \ref{assm_ErgodicU} hold.  Then, $P^h(t)$ and $Q^h(t)$ have uniform in time bounded moments: for every $m \geq 1$ and all sufficiently small $h>0$ there is a constant $C(m) > 0$ such that 
		$$ \sup_{t\geq0} 	\E \Big[ |P^h(t)|^m \Big] \leq C(m) \text{\quad and \quad}  \sup_{t\geq 0} 	\E \Big[  |Q^h(t)|^m \Big] \leq C(m). $$
	\end{prop}

	Define by $\mathcal{L}$ and $\mathcal{L}^h$ the infinitesimal generators of \eqref{eq:origin} and \eqref{eq:tamed},
    respectively. Then for any $\psi \in\mathcal C^2(\RR^{2d};\RR)$, we have
	\begin{align*}
		\mathcal{L} \psi(p,q) &=(\partial_p \psi(p,q))^\top (-\upsilon p - f(q))+ (\partial_q \psi(p,q))^\top p +\frac{\sigma^2}{2} \tr( \partial_{pp} \psi (p,q) ), \\
		\mathcal{L}^h \psi(p,q) &=( \partial_{p}\psi(p,q))^\top (-\upsilon p - f^h(q))+ (\partial_q \psi(p,q))^\top p +\frac{\sigma^2}{2}\tr(  \partial_{pp} \psi(p,q) ),
	\end{align*}
	where we denote by $\partial_p\psi$ the partial derivative with respect to the first variable, by $\partial_q\psi$ the partial derivative with respect to  the second one, and by $\partial_{pp}\psi$ the Hessian matrix. 
	For $g\in\mathcal C^{\infty}(\RR^{2d}; \RR)$ with derivatives having growth not faster than polynomial, we define $u(t,p,q):= \E\big[g(P(t), Q(t)) \big| (P(0), Q(0))=(p,q)\big] $, which is the solution of the backward Kolmogorov equation written in forward time: 
    \begin{align*}
    &\partial_t u(t,p,q) =\mathcal{L} u(t,p,q), \qquad
    u(0,p,q)=g(p,q).  
    \end{align*}
	Under Assumption~\ref{assm_ErgodicU}, 
	for the Lyapunov function $V(p,q):=\frac 12 |p|^2 +U(q)+\upsilon \alpha_1  p^\top q$, we obtain the standard dissipative estimate:
	\begin{align}\label{equ_Talay1}
		\mathcal L V(p,q) \leq -\delta V(p,q) +C
	\end{align}
	and Assumption~\ref{assm_CandPoly} implies
	\begin{align} \label{equ_Talay2}
		V(p,q)
		\geq \delta(|p|^2+|q|^{2})-C
	\end{align}
	with some $\delta>0$ and $C>0$.  
    Assumptions~\ref{assm_CandPoly} and~\ref{assm_ErgodicU}, which yield \eqref{equ_Talay1}--\eqref{equ_Talay2}, together with the assumptions made on smoothness and growth of $g$ ensure that the conditions of \cite[Theorem 3.1]{Talay2002} are satisfied. Hence, for any integer $l\geq 1$, 
	 there exist positive constants $C$, $\kappa_1$, $\kappa_2$, and $\gamma$ (they depend on $g$ and $l$) such that
	\begin{align}\label{eq:Dku}
		\|D^l u(t,p,q)\|_{\op} \leq C (1+|p|^{\kappa_1}+|q|^{\kappa_2})\exp(-\gamma t)
	\end{align}
	for any $t \ge 0$ and  $(p^\top,q^\top)^\top \in\RR^{2d}$. With these properties, the detailed proof of Theorem~\ref{thm_weakconvergencetamed} is put in Section~\ref{subsec_proofweakconvergenceequ}.

\subsection{Proof of Theorem~\ref{thm_weakconvergeconformalsym}} \label{sec_proofpart2}

    In this section, we present the proofs of the ergodicity and weak convergence results of the general numerical method \eqref{eq:generalP}--\eqref{eq:generalQ}. 
	
	Using the Lyapunov function $V^h$ defined by \eqref{equ_defVh},
	we show the uniform moment boundedness of $\{(P^h_n, Q^h_n)\}_{n\in\N}$ in the following proposition, whose proof is provided in Section~\ref{subsec_proofmomentnumerical}. 
	\begin{prop} \label{prop_momentboundedgeneral}
		Let  Assumptions  
		\ref{assm_CandPoly} and \ref{assm_ErgodicU} hold.  Then, for any $\alpha \geq 1$ and sufficiently small $h>0$, it holds that
		$$ \sup_{n\geq 0}	\E \big[  |P^{h}_n|^\alpha \big] \leq C(\alpha) \text{\quad and \quad}  \sup_{n\geq 0}  	\E \big[  |Q^{h}_n|^\alpha \big] \leq C(\alpha). $$
	\end{prop}
	Using Proposition~\ref{prop_momentboundedgeneral}, we next prove Theorem~\ref{thm_weakconvergeconformalsym}.
	\begin{proof}[Proof of Theorem~\ref{thm_weakconvergeconformalsym}] 
		The proof of Proposition~\ref{prop_momentboundedgeneral} (see \eqref{equ_numericalboundClaim}) implies the Lyapunov condition for the $m$-step chain $\{(P^h_{mn},Q^h_{mn})\}_{n\in\N}$ with $m$ from Assumption~\ref{assp:nondeg}.
        Thanks to Assumption~\ref{assp:nondeg}, the transition kernel of this chain has the density $\rho(x,y)=\phi_\xi\big(\Gamma_x^{-1}(y)\big)\,\big|\det D\Gamma_x^{-1}(y)\big|$, where $\phi_\xi>0$ is the Gaussian density of $\xi$. The density $\rho(x,y)$ is continuous in $(x,y)$ and positive everywhere, hence the minorization condition holds on every compact set 
        (cf. \cite[Lemma 2.3]{Mattingly2002SPA}). Thus, by \cite[Theorem 2.5]{Mattingly2002SPA}, the $m$-step chain has a unique invariant measure for sufficiently small $h>0$. Every invariant measure of $\{(P^h_n,Q^h_n)\}_{n\in\N}$ is invariant for the $m$-step chain, and 
        $$\frac1m\sum_{j=0}^{m-1}\int_{\mathbb R^{2d}}\mathcal P_j(x,A) \mu(\mathrm dx),\qquad A\in\mathcal B(\mathbb R^{2d}),$$ 
        where ${\mathcal P}_j$ is the $j$-step transition kernel of the chain, is invariant for the chain whenever $\mu$ is invariant for the $m$-step chain. Consequently, $\{(P^h_n,Q^h_n)\}_{n\in\N}$ has a unique        
        invariant measure $\pi_{h,\text{num}}$ for sufficiently small $h>0$. 
		
		Next, we prove the weak convergence results. With the uniform moment boundedness established in  Proposition~\ref{prop_momentboundedgeneral}, the proof follows the ideas in \cite{Milstein2004book}.
		We denote $X(t) :=(P(t)^\top, Q(t)^\top)^\top $ and $Y^h(t_n) := ((P^h_{n})^\top, (Q^h_{n})^\top)^\top$.  Furthermore, denote by $X_{s, x}(t)$ the solution $X(t)$ of \eqref{eq:origin}  with initial condition $X(s)=x$, and  by $Y^h_{t_k,y}(t_n)$ the numerical solution $Y^h(t_n)$ computed by \eqref{eq:generalP}--\eqref{eq:generalQ} with initial condition $Y^h(t_k) =y$. 
		For notational simplicity, we let $\text{Poly}(x_1, \ldots, x_r)$ denote a generic function of the form $C(1+|x_1|^{m_1}+\cdots+|x_r|^{m_r})$, 
		where $r \in \mathbb{N}^+$, and the constants $C > 0$ and $m_i \in \mathbb{N}^+$, $i=1, \ldots, r$ may vary from one place to another. 
		
		Notice that 
		\begin{align} \label{equ_weakerror}
			&\quad 	\big|\pi(g) -\pi_{h, \text{num}} (g)\big| \nn \\
			&=
			\Big| \lim_{N\to \infty}\frac{1}{N}\sum^N_{n=1} \E \big[ g(Y^h(t_n))- g(X(t_n )) \big] \Big|  \nn \\
			&= \Big|\lim_{N\to \infty}\frac{1}{N}\sum^N_{n=1} \sum^{n-1}_{i=0} \E \big[ u(\tau_{n,i+1}, Y^h_{t_i, Y^h(t_i)}(t_{i+1})) - u( \tau_{n,i+1}, X_{t_i,Y^h(t_i)}(t_{i+1})) \big]  \Big|,
		\end{align}
		where we denote $\tau_{n,i}:= t_n-t_i$ for $1\leq i\leq n$ and use the time-homogeneity in the last step.
		Using \eqref{equ_weakerror}, the weak convergence result is obtained once we derive	the estimate of $ \E\big[u\big(\tau_{n,i+1}, Y^h_{t_{i},Y^h(t_i)} (t_{i+1})\big)\big] -\E \big[	u\big(\tau_{n,i+1}, X_{t_{i},Y^h(t_i)} (t_{i+1})\big) \big]$.

		The proof of (\romannumeral 1) is analogous to that of (\romannumeral 2). As the latter is more involved, we provide its proof here as a representative case.	By the Taylor expansion, we have
		\begin{align} \label{equ_onestepu}
			&\quad \E\big[u\big(\tau_{n,i+1}, Y^h_{t_{i},x} (t_{i+1})\big)\big] -\E \big[	u\big(\tau_{n,i+1}, X_{t_{i},x} (t_{i+1})\big) \big] \nn \\
			&= \sum^{5}_{l=1} \sum^{2d}_{k_1,\ldots, k_l=1} \frac{1}{l!} \frac{\partial^lu (\tau_{n,i+1}, x)}{\partial x_{k_1}\ldots \partial x_{k_l}}\mathbb{E} \bigg[ \prod_{j=1}^{l} (\delta_{Y^h,x})_{k_j} - \prod_{j=1}^{l} (\delta_{X,x})_{k_j} \bigg] + \mathbb{E}[\mathcal{R}_{n,i+1}],
		\end{align}
		where $\delta_{X,x}:= X_{t_i,x}(t_{i+1})-x$,  $\delta_{Y^h,x}:=Y^h_{t_i,x}(t_{i+1})-x$, $(\delta_{X,x})_{k_j}$ and $(\delta_{Y^h,x})_{k_j}$ denote the corresponding increments of the $k_j$-th coordinate, and
		\begin{align*}
			\mathcal{R}_{n,i+1} := \sum_{k_1,\ldots, k_6=1}^{2d} \frac{1}{5!} \bigg( &\int_{0}^{1} (1-s)^{5}  \frac{\partial^6u(\tau_{n,i+1}, x + s\delta_{Y^h,x})}{\partial x_{k_1}\ldots \partial x_{k_6}}
			ds \prod_{j=1}^{6} (\delta_{Y^h,x})_{k_j} \\
			&- \int_{0}^{1} (1-s)^{5} \frac{\partial^6u(\tau_{n,i+1}, x + s\delta_{X,x})}{\partial x_{k_1}\ldots \partial x_{k_6}}
			ds \prod_{j=1}^{6} (\delta_{X,x})_{k_j} \bigg).
		\end{align*}
		By \eqref{eq:Dku} and the bounds, $\E[|\delta_{X,x}|^\alpha]\leq \text{Poly}(x)h^{\alpha/2}$  and $\E[|\delta_{Y^h,x}|^\alpha]\leq \text{Poly}(x)h^{\alpha/2}$ for $\alpha\geq 1$, we arrive at
		\begin{align} \label{equ_onestepremain}
			|\E[\mathcal R_{n,i+1}] | \leq e^{-\gamma \tau_{n,i+1}}\text{Poly}(x)h^3.
		\end{align}
		
		Applying the It\^o--Taylor expansion to \eqref{eq:origin} yields 
		\begin{align} \label{equ_ItoTaylor}
			\begin{pmatrix}
				P(t+h) \\
				Q(t+h)
			\end{pmatrix} 
			&= 	\begin{pmatrix}
				P(t) \\
				Q(t)
			\end{pmatrix} + \begin{pmatrix}
				-\upsilon P(t) -f(Q(t)) \\
				P(t)
			\end{pmatrix} h +\begin{pmatrix}
				\sigma \\ 0
			\end{pmatrix}  \int^{t+h}_t dW(s) \nn  \\
			&\quad  + \begin{pmatrix}
				-\upsilon \sigma \\ \sigma
			\end{pmatrix}\int^{t+h}_t \int^{s}_t dW(r) ds \nn\\
            &\quad+\int^{t+h}_{t} \int^s_{t} \begin{pmatrix}
				\upsilon^2 P(r) +\upsilon f(Q(r)) -  \nabla f(Q(r))P(r) \\
				-\upsilon P(r) -f(Q(r)) 
			\end{pmatrix} dr ds
		\end{align} 
		for $t\geq0$.
		With \eqref{equ_ItoTaylor}, for any given initial value $x=(p^\top,q^\top)^\top \in\RR^{2d}$, we define the auxiliary stochastic variable $Y_{\text{aux}}:=(P_{\text{aux}}^\top, Q_{\text{aux}}^\top)^\top$  by
		\begin{align} \label{equ_PhMT}
			P_{\text{aux}}
			&= \big(1-\upsilon h +\frac{\upsilon ^2 h^2}{2} \big) p  +\big(- h + \frac{\upsilon h^2}{2}\big)  f(q) -\frac{h^2}{2} \nabla f(q) p \nn\\
			&\quad+\sigma (W(t+h)-W(t)) -\upsilon \sigma\int^{t+h}_t \int^{s}_t dW(r) ds, \\
            \label{equ_QhMT}
			Q_{\text{aux}}
			&=q+ \big(h-\frac \upsilon 2  h^2\big) p   - \frac 12 h^2 f(q) +\sigma\int^{t+h}_t \int^{s}_t dW(r) ds,
		\end{align}
		and denote $\delta_{Y,x}:= Y_{\text{aux}}-x$. 
		Comparing \eqref{equ_PhMT}--\eqref{equ_QhMT} with \eqref{equ_ItoTaylor}, for $l=1,\ldots, 5$, the moment boundedness of $X$ gives
		$
			|	\E[ \prod_{j=1}^{l}(\delta_{Y,x} )_{k_j}-\prod_{j=1}^{l} (\delta_{X,x})_{k_j } ] | \leq  \text{Poly} (x) h^3. 
		$
		Since the triangle inequality gives
		\begin{align*} %
			&\quad \Big|	\E\Big[ \prod_{j=1}^{l}(\delta_{Y^h,x} )_{k_j}-\prod_{j=1}^{l} (\delta_{X,x})_{k_j } \Big] \Big| \nn \\
			&\leq 	\Big|	\E\Big[ \prod_{j=1}^{l}(\delta_{Y^h,x} )_{k_j}-\prod_{j=1}^{l} (\delta_{Y,x})_{k_j } \Big] \Big| + 	\Big|	\E\Big[ \prod_{j=1}^{l}(\delta_{Y,x} )_{k_j}-\prod_{j=1}^{l} (\delta_{X,x})_{k_j } \Big] \Big|,
		\end{align*}  
		we only need to prove that for $l=1,\ldots, 5$,
		\begin{align} \label{equ_weakerrorpart2}
			|	\E[ \prod_{j=1}^{l}(\delta_{Y^h,x} )_{k_j}-\prod_{j=1}^{l} (\delta_{Y,x})_{k_j } ] |	\leq\text{Poly} (x) h^3.
		\end{align}
		
		By Taylor’s formula with integral remainder, we get 
        \begin{align*} 
        \bar H_{j,n}-D^2U^h(Q^h_n)=\int_0^1(1-r)\,D^3U^h\big(Q^h_n+r(\mathcal{Q}_{j,n}-Q^h_n)\big)\big[\mathcal{Q}_{j,n}-Q^h_n\big]\,\rd r,  \;\; j=p,q,
        \end{align*} 
        where $D^3U^h(x)[v]$ denotes the matrix with entries $\sum_{k}\partial_{ijk}U^h(x)v_k$. Combining this  with the 
		given assumptions, we arrive at
		\begin{align} \label{equ_PhS22}
			P^h_{n+1} 
			&= \big(1-\upsilon h +\frac{\upsilon ^2 h^2}{2} \big) P^h_n  +\big(- h + \frac{\upsilon h^2}{2}\big)  f^h(Q^h_n) -\frac{h^2}{2} \nabla f^h(Q^h_n) P^h_n +\sigma \Delta W_n \nn \\
			&\quad -\upsilon \sigma \Delta B_n +(\mathcal{Z}_{p,n,1}+\upsilon \sigma \Delta B_n  )+h \bar H_{p,n} \mathcal{Z}_{p, n,2}+\tilde{\mathcal{R}}_{p,n},
		\end{align}
		and
		\begin{align} \label{equ_QhS22}
			Q^h_{n+1} 
			&=Q^h_n+ \big(h-\frac \upsilon 2  h^2\big) P^h_n   - \frac 12 h^2 f^h(Q^h_n) \nn \\
			&\quad +\sigma \Delta B_n+(\mathcal{Z}_{q,n,1}-\sigma \Delta B_n) + h^2 \bar H_{q,n} \mathcal{Z}_{q, n,2}+\tilde{\mathcal{R}}_{q,n}.
		\end{align} 
		Here, the higher-order remainders are defined by
		\begin{align*}
			\tilde{\mathcal{R}}_{p,n}:= \mathcal{R}_{p,n}+  a_{0,p,2} h^3 \bar H_{p,n} P^h_n 
            -\frac{1}{2}h^2 \big(\bar H_{p,n}-D^2U^h(Q^h_n)\big) P^h_n 
			+a_{0,q} h^3  \bar H_{p,n}  f^h(Q^h_n),
		\end{align*}
		and $\tilde{\mathcal{R}}_{q,n}:=  \mathcal{R}_{q,n}+ b_{0,p} h^3  \bar H_{q,n} P^h_n $. 
		
		Comparing \eqref{equ_PhS22}--\eqref{equ_QhS22} with \eqref{equ_PhMT}--\eqref{equ_QhMT}, the proof of \eqref{equ_weakerrorpart2} reduces to estimating the remainders  $\mathcal{Z}_{p,n,1}+\upsilon \sigma \Delta B_n, h  \bar H_{p,n} \mathcal{Z}_{p, n,2}, \tilde{\mathcal{R}}_{p,n}$, and $ \mathcal{Z}_{q,n,1}-\sigma \Delta B_n$, $h^2  \bar H_{q,n} \mathcal{Z}_{q, n,2}$, $\tilde{\mathcal{R}}_{q,n}$. 
          
		Under the assumptions, we naturally have $
			\E[\mathcal{Z}_{p,n,1}+\upsilon \sigma \Delta B_n   |\mathcal F_{t_n} ]=0$, $\E[|\mathcal{Z}_{p,n,1}+\upsilon \sigma \Delta B_n  |^2 |\mathcal F_{t_n} ]=\mathcal{O} (h^3)$, 
		$
			\E[ \mathcal{Z}_{q,n,1}-\sigma \Delta B_n | \mathcal{F}_{t_n}]=0$, $\E[ |\mathcal{Z}_{q,n,1}-\sigma \Delta B_n|^2 | \mathcal{F}_{t_n}]=\mathcal{O} (h^3)
		$,
        $
            \E[(\mathcal{Z}_{p,n,1}+\upsilon \sigma \Delta B_n ) (\Delta W_n )^\top|\mathcal F_{t_n} ]=\mathcal{O} (h^3)$, and $\E[(\mathcal{Z}_{q,n,1}-\sigma \Delta B_n )(\Delta W_n)^\top |\mathcal F_{t_n} ]=\mathcal{O} (h^3)$.
		Next we provide estimates for the other terms. 
         By Proposition~\ref{prop:fh-fS2}\,(\romannumeral 5), $\|D^3U^h\|_{\op}$ grows at most polynomially.
        Therefore, for any $\alpha\geq 1$, it follows that
		\begin{align*}
			&\quad \E\big[\big|h^2\big(\bar H_{p,n}-D^2U^h(Q^h_n)\big) P^h_n\big|^\alpha \big| \mathcal F_{t_n} \big]\\
            & \leq \frac 12 h^{2\alpha}\E\Big[\sup_{r\in[0,1]}\big\|D^3U^ h\big(Q^h_n+r(\mathcal{Q}_{p,n}-Q^h_n)\big)\big\|_{\op}^\alpha|\mathcal{Q}_{p,n}-Q^h_n|^\alpha\Big| \mathcal F_{t_n} \Big]|P^h_n|^\alpha\\
            &\leq \text{Poly} (P^h_n,Q^h_n)h^{3\alpha},
		\end{align*}
        where we used $|\mathcal{Q}_{p,n}-Q^h_n|\leq|Q^h_{n+1}-Q^h_n|+h|P^h_n|$ and $\E[ |Q^h_{n+1}-Q^h_n|^\alpha|\mathcal{F}_{t_n}]\leq \text{Poly}(P^h_n,Q^h_n)h^\alpha$.
        Similarly, we derive $|  h^3 \bar H_{j,n} P^h_n|^\alpha \leq \text{Poly} (P^h_n,Q^h_n,Q^h_{n+1})h^{3\alpha}, \,j=p,q$ and $| h^3 \bar H_{p,n}  f^h(Q^h_n)|^\alpha\leq \text{Poly} (P^h_n,Q^h_n,Q^h_{n+1}) h^{3\alpha}.$ 
        Combining these estimates with the assumptions on $\mathcal{R}_{p,n}$ and $\mathcal{R}_{q,n}$, for any $\alpha\geq 1$, we arrive at
        $\E[|\tilde{\mathcal{R}}_{p,n}|^\alpha ]\leq \E[\text{Poly} (P^h_n,\allowbreak Q^h_n,Q^h_{n+1})]h^{3\alpha}$ and $\E[|\tilde{\mathcal{R}}_{q,n}|^\alpha]\leq \E[ \text{Poly} (P^h_n,Q^h_n,Q^h_{n+1})]h^{3\alpha}.$ 
		For the term $h \bar H_{p,n} \mathcal{Z}_{p, n,2}$, 
		we get 
		\begin{align*}
			\big|	\E\big[	h \bar H_{p,n} \mathcal{Z}_{p, n,2} \big] \big|
			&\leq \big| \E\big[	hD^2U^h(Q^h_n) \mathcal{Z}_{p, n,2} \big]\big|+ 
            h\E\big[\|\bar H_{p,n}-D^2U^h(Q^h_n)\|_{\op} |\mathcal{Z}_{p, n,2}| \big]  \\
			&\leq\E[ \text{Poly}(P^h_n,Q^h_n, Q^h_{n+1})]h^{\frac 72}
		\end{align*}
		and
		$
		\E\big[	|h \bar H_{p,n} \mathcal{Z}_{p, n,2} |^{\alpha} \big]  \leq \E[ \text{Poly}(P^h_n,Q^h_n, Q^h_{n+1})] h^{5\alpha/2}.
		$
		For the term $h^2 \bar H_{q,n} \mathcal{Z}_{q, n,2}$, we similarly 
		derive 
		$\E[	|h^2 \bar H_{q,n} \mathcal{Z}_{q, n,2}|^{\alpha}]  \leq\E[ \text{Poly}(P^h_n,Q^h_n, Q^h_{n+1})] h^{7\alpha/2}$ for $\alpha\geq1$.
		
		Combining the estimates above and using Propositions \ref{prop:fh-fS2} and \ref{prop_momentboundedgeneral}, we obtain
		\begin{align*} %
			\big| 	\E\big[(\delta_{Y,x} )_{k_1} - (\delta_{Y^h,x})_{k_1}\big] \big| 
			&\leq \big(h+Ch^2  \big) | f(q) - f^h(q)| +Ch^2 | p | \| D^2 U(q) - D^2U^h(q)\|_{\op} \nn \\
			&\quad +\E[\text{Poly}(x,Y^h_{t_i,x}(t_{i+1}))] h^3 \nn\\
			&\leq \text{Poly}(x)  h^3
		\end{align*}
		and
		\begin{align}  \label{equ_weakpro2}
			&\quad \big| 	\E\big[ \prod_{j=1}^{2}(\delta_{Y,x} )_{k_j} - \prod_{j=1}^{2}(\delta_{Y^h,x})_{k_j}\big] \big|\nn\\
            &\leq \big| 	\E\big[ (\delta_{Y,x}  - \delta_{Y^h,x})_{k_1}(\delta_{Y,x}  )_{k_2}\big] \big|+\big| 	\E\big[ (\delta_{Y^h,x})_{k_1}(\delta_{Y,x}  - \delta_{Y^h,x})_{k_2}\big] \big| \leq   \text{Poly}(x)  h^3.
		\end{align}
		Following a similar argument as in \eqref{equ_weakpro2} and using that $\mathcal{Z}_{p,n,i}$, $\mathcal{Z}_{q,n,i}$,  $\Delta B_n$ and $\Delta W_n$ are jointly Gaussian with mean zero, we have 
        \begin{align}\label{equ_onestep}
        |	\E[ \prod_{j=1}^{l}(\delta_{Y^h,x} )_{k_j}-\prod_{j=1}^{l} (\delta_{Y,x})_{k_j } ] |	\leq  \text{Poly}(x) h^3
        \end{align}
        for $l=2,\ldots, 5$.
		
		Using Proposition~\ref{prop_momentboundedgeneral}, \eqref{eq:Dku}, \eqref{equ_onestepu}, \eqref{equ_onestepremain}, and \eqref{equ_onestep},  we  obtain
		\begin{align*}
			&\quad \big|	\E\big[u\big(\tau_{n,i+1}, Y^h_{t_{i},Y^h(t_i)} (t_{i+1})\big)\big] -\E \big[	u\big(\tau_{n,i+1}, X_{t_{i},Y^h(t_i)} (t_{i+1})\big) \big]  \big|\\
            &\leq e^{-\gamma\tau_{n,i+1}} \E[\text{Poly}(Y^h(t_i))] h^3 \leq e^{-\gamma\tau_{n,i+1}} C(P(0),Q(0) )h^3.
		\end{align*}
  Since $X_{t_n,y}(t_n)=y$, $Y^h(t_0)=X(0)$,  the inner sum in \eqref{equ_weakerror} collapses to  $\E[g(Y^h(t_n))]-\E[g(X(t_n))]$ for each $n$.
  Summing the last bound over $i$ with $\tau_{n,i+1}=(n-i-1)h$, we get
		\begin{align*}
			\big|\E[g(Y^h(t_n))]-\E[g(X(t_n))]\big|
			&\leq C(P(0),Q(0))\,h^3\sum_{j=0}^{n-1}e^{-\gamma jh}
            \leq C(P(0),Q(0))\,h^2
		\end{align*}
		for all $n\in\N$ and sufficiently small $h>0$, which is the uniform-in-time bound in (\romannumeral 2).
		Its average in \eqref{equ_weakerror} gives the required $\big|\pi(g) -\pi_{h, \text{num}} (g)\big|  \leq Ch^2$, which ends the proof of (\romannumeral 2). 
	\end{proof}
    
	\section{Numerical experiments} \label{sec_numericalexperiment}
	In this section, we present numerical experiments to illustrate our theoretical results, namely, the orders of uniform-in-time weak convergence of the five schemes of Section~\ref{sec_mainresultsofExplicitmethods}, the approximation of the Gibbs density by ensemble and time averages, a comparison of the mollified and tamed potential approaches, and the preservation of conformal symplecticity. As a representative non-convex example, we consider the Mexican-hat-type potential $U(q)=\frac 16 |q|^6 -\frac 14 |q|^4 +1 $, $q\in\RR^d$.  The function $\theta$ is chosen as in \eqref{eq:theta}.
    For this potential, we have $r_1=6$, $r_2=5$, and $L_2 =\frac 16$. 
	\subsection{Uniform-in-time weak convergence}
	 In this subsection, we numerically investigate the uniform-in-time weak convergence rates of the methods proposed in Section~\ref{sec_mainresultsofExplicitmethods}. Here, we set $T=10$, $\sigma=\sqrt{2}$, $d=2$, $\upsilon=1$, and $R_0=2$, with initial values $P(0)=(1,1)^\top$ and $Q(0)=(-1,-1)^\top$. We use the test function $g(p,q)=\sin(p)^\top\sin(q)$, and step sizes are chosen as $h_i=2^{-i}$ with $i=2,3,4,5$. In the following experiments, we also count how many times truncation occurs, which we call the truncation count. The truncation count for a method with step size $h$ is defined as the total number of times that $|q|>R(h)$ when computing the mollified force, accumulated over all sample paths and all time steps in $[0,T]$. Here, $q$ is the position used in each force calculation.
     \begin{figure}[ht]
    	\centerline{\includegraphics[scale=0.35]{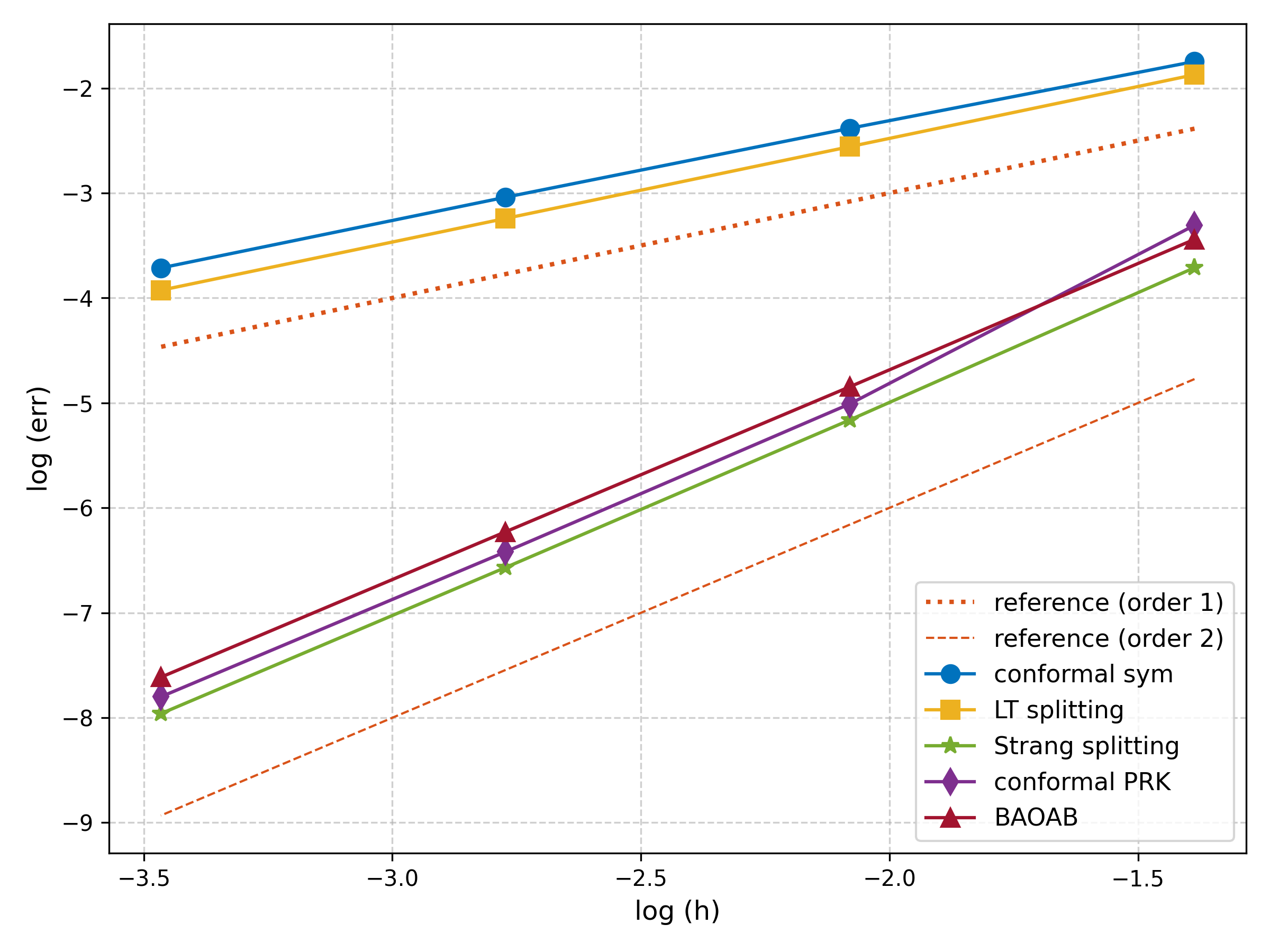}}
			\caption{Maximum weak error over $[0,10]$}
			\label{fig:weakorderEg1}
          
        \end{figure}

\textit{Maximum weak error over $[0,10]$.}
In Figure~\ref{fig:weakorderEg1}, using Monte Carlo simulation, we estimate the maximum weak error over $[0,10]$ computed as 
$$
    \max_{1\leq n\leq K_h}
    \Big|
        \frac 1M \sum^M_{j=1}\big[
            g(P_n^{h,j},Q_n^{h,j})
            -
            g(P_{m_n}^{h_\text{ref},j},Q_{m_n}^{h_\text{ref},j})
        \big]
    \Big|.
$$
Here, $K_h:=T/h$ denotes the total number of time steps, $M$ is the number of Monte Carlo sample paths, $(P_n^{h,j}, Q_n^{h,j})$ and
$(P^{h_{\text{ref}},j}_{m_n},Q^{h_{\text{ref}},j}_{m_n} )$
denote the numerical and reference solutions for the $j$-th trajectory, respectively, and $m_n:=\frac{nh}{h_{\text{ref}}}$  ensures that the numerical and
reference solutions are evaluated at the same time. We set $M=10^5$ and compute the reference solution $(P_{m}^{h_\mathrm{ref}},Q_{m}^{h_\mathrm{ref}})$  by using the BAOAB method proposed in Proposition~\ref{prop_exampleOU} with the fine step size $h_{\mathrm{ref}}=2^{-12}$. 
It can be observed in Figure~\ref{fig:weakorderEg1} that
the conformal symplectic Euler and Lie--Trotter splitting methods proposed in Proposition~\ref{prop_example1st} exhibit first-order weak convergence, whereas the conformal PRK, Strang splitting, and BAOAB methods proposed in Propositions \ref{prop_exampleSplitting}--\ref{prop_exampleOU} exhibit second-order weak convergence, which confirms our theoretical results.

For $h=2^{-2}$, the truncation counts are $288$, $1537$, $2506$, $4308$, and $1856$ for the conformal symplectic Euler, Lie--Trotter splitting, Strang splitting, conformal PRK, and BAOAB methods, respectively. No truncations occur for $h=2^{-3},2^{-4},2^{-5}$ or in the reference simulation. 
For $h=2^{-2}$, the mollification is active much more often than the Gibbs measure would suggest, since it gives the region $|q|>R(2^{-2})\approx2.38$ a mass of only about $7 \times 10^{-12}$. These excursions are caused by the large step. For $h \le 2^{-3}$, on every simulated path the mollified schemes coincide with the same schemes applied to the original potential and the results would be the same for any $R_0 \ge 2$. Here the mollification does not disturb the dynamics where it is not needed, which is the aim of its construction (see the discussion after \eqref{eq:molpot}). At the same time, our approach guarantees the uniform-in-time behaviour established in Theorem~\ref{thm_weakconvergeconformalsym}, which explicit schemes for the original potential lack in general \cite{Mattingly2002SPA,Hutzenthaler2011EulerDiverge}.
The times at which the maximum weak errors are attained are reported in Table~\ref{tab:maximumtimes}. For all methods and tested step sizes, these maxima occur strictly before the terminal time $T=10$.

\begin{table}[htbp]
    \centering
    \caption{Times of maximum weak errors}
    \label{tab:maximumtimes}
    \begin{tabular}{lcccc}
        \hline 
        Method & $h=2^{-2}$ & $h=2^{-3}$ & $h=2^{-4}$ & $h=2^{-5}$ \\
        \hline
        Conformal symplectic Euler & 1.75000 & 1.62500 & 1.62500  & 1.62500   \\
        Lie--Trotter splitting     & 1.25000 & 1.37500 & 1.43750 & 1.40625 \\
        Strang splitting          & 1.75000 & 1.75000  & 1.68750 & 0.87500   \\
        Conformal PRK             & 2.50000  & 2.50000   & 2.56250 & 2.53125 \\
        BAOAB                     & 1.00000  & 0.87500 & 0.75000   & 0.90625 \\
        \hline
    \end{tabular}
\end{table}
\begin{figure}[ht]
			\centerline{\includegraphics[scale=0.35]{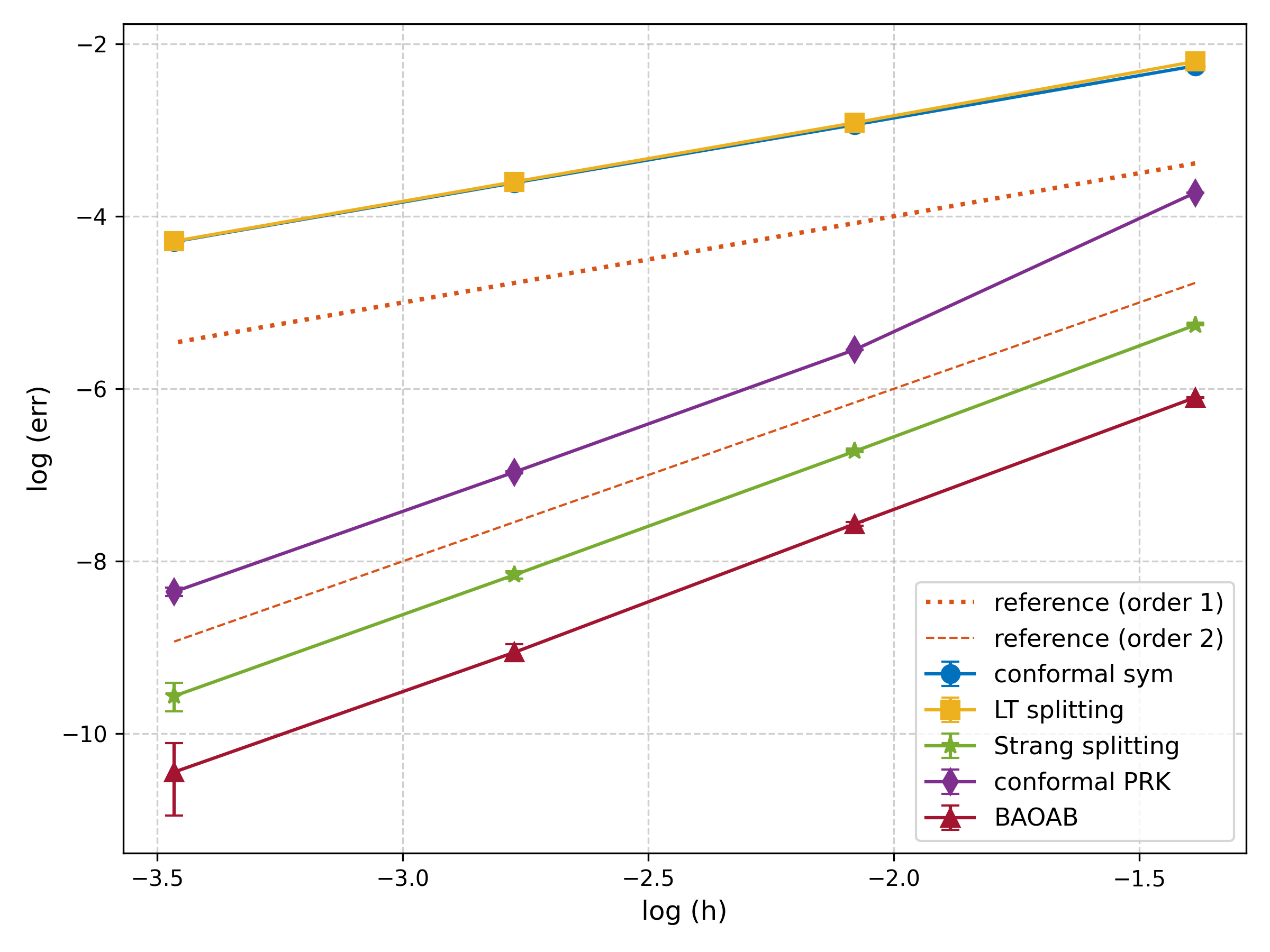}}
			\caption{Terminal-time weak error at $T=10$ relative to the ergodic limit}
			\label{fig:longtimeEg1}
	\end{figure}
\textit{Terminal-time weak error relative to the ergodic limit.}
In Figure~\ref{fig:longtimeEg1}, we estimate  the terminal-time weak error relative to the ergodic limit $\pi(g)$:
$
|\frac 1M \sum^M_{j=1} g(P_{K_h}^{h,j},Q_{K_h}^{h,j})-\pi(g)|
$, where $(P^{h,j}_{K_h},Q^{h,j}_{K_h} )$ denotes the numerical solution of the $j$-th trajectory at the final time step $K_h=T/h$. We take $M=10^{10}$. The error bars indicate the Monte Carlo sampling error.
Note that $\pi(g)=0$ for the test function $g(p,q)=\sin(p)^\top \sin(q)$. 
By \eqref{equ_ensemble1}--\eqref{equ_ensemble2}, this error is at most $C(h^{p}+e^{-\gamma T})$, where $p=1$ or $p=2$ is the order of the method. The observed rates show no saturation down to errors of about $3\times10^{-5}$, hence the relaxation term is negligible at $T=10$, and this experiment examines the error relative to the ergodic limit $\pi(g)$ in the near-equilibrium regime.
For $h=2^{-2}$, the truncation counts are $24445150$, $155624087$, $265505529$, $443494595$, and $156788505$ for the conformal symplectic Euler, Lie--Trotter splitting, Strang splitting, conformal PRK, and BAOAB methods, respectively. The observed convergence rates are consistent with those in Figure~\ref{fig:weakorderEg1}: first order for the conformal symplectic Euler and Lie--Trotter splitting methods, and second order for the conformal PRK, Strang splitting, and BAOAB methods. The errors of the two first-order methods are almost equal here resulting in their curves overlapping in Figure~\ref{fig:longtimeEg1}.

The two experiments provide complementary numerical evidence for uniform-in-time weak convergence. The first captures the maximum weak errors over $[0,10]$, which are attained during the initial transient phase, while the second examines the errors relative to the invariant average at a time when the system is expected to be close to equilibrium. The consistent convergence rates observed in these two regimes provide numerical support for the uniform-in-time weak error estimates established in Theorem~\ref{thm_weakconvergeconformalsym} and Propositions~\ref{prop_example1st}--\ref{prop_exampleOU}.

	\subsection{Density approximation}
	Recall that the Gibbs density associated with the exact solution is %
$
    \rho(p,q)=\frac{1}{Z} \exp\big(
-\frac{2\upsilon}{\sigma^2}H(p,q) 
\big),
$
    where
$
        Z=
\big(\frac{\pi\sigma^2}{\upsilon}\big)^{\frac d2}
\int_{\mathbb R^d}
\exp\big(
-\frac{2\upsilon}{\sigma^2}U(q)
\big)\,dq
$
    and $H(p,q) = \frac 12 |p|^2+ U(q) $.  In this subsection, we employ the ensemble-average estimator and the time-average estimator with Gaussian kernels to approximate the Gibbs density $\rho$. 

\begin{figure}[ht]
      \centering
      \subfloat[$T=1$]{%
        \label{fig:ensemblea}%
        \includegraphics[width=0.30\textwidth]
          {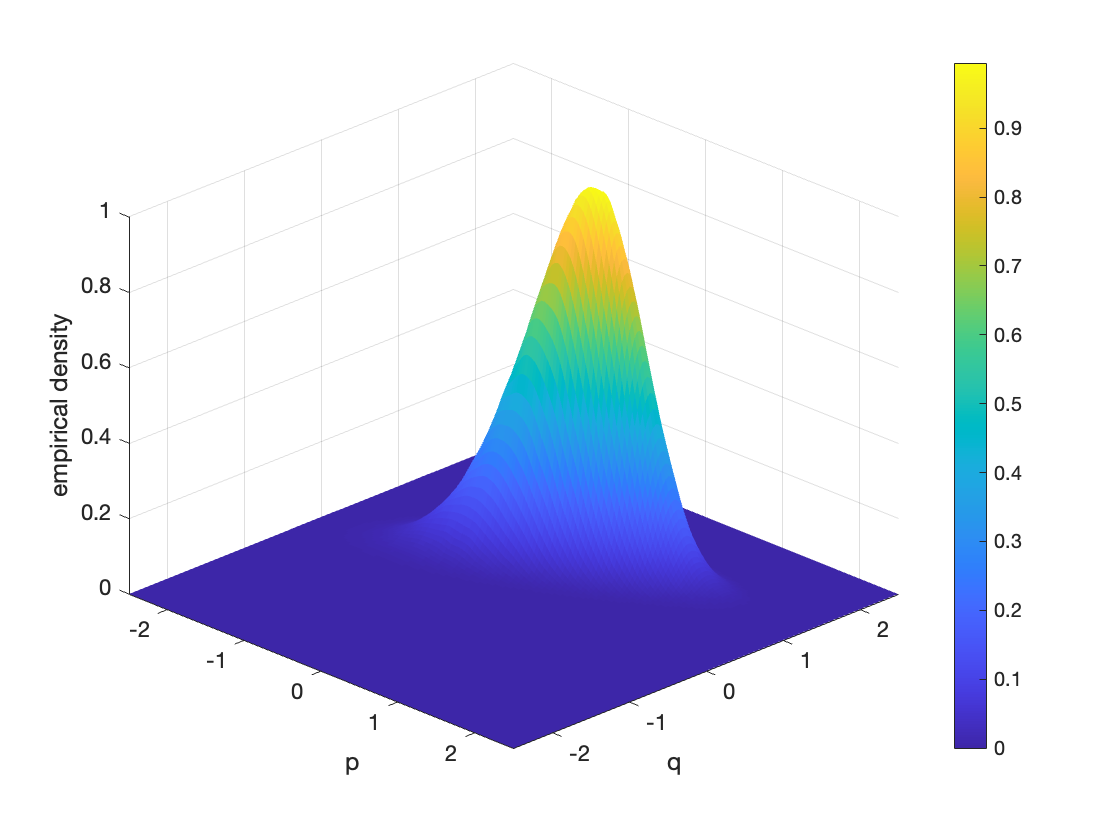}%
      }\hfill
      \subfloat[$T=100$]{%
        \label{fig:ensembleb}%
        \includegraphics[width=0.30\textwidth]
          {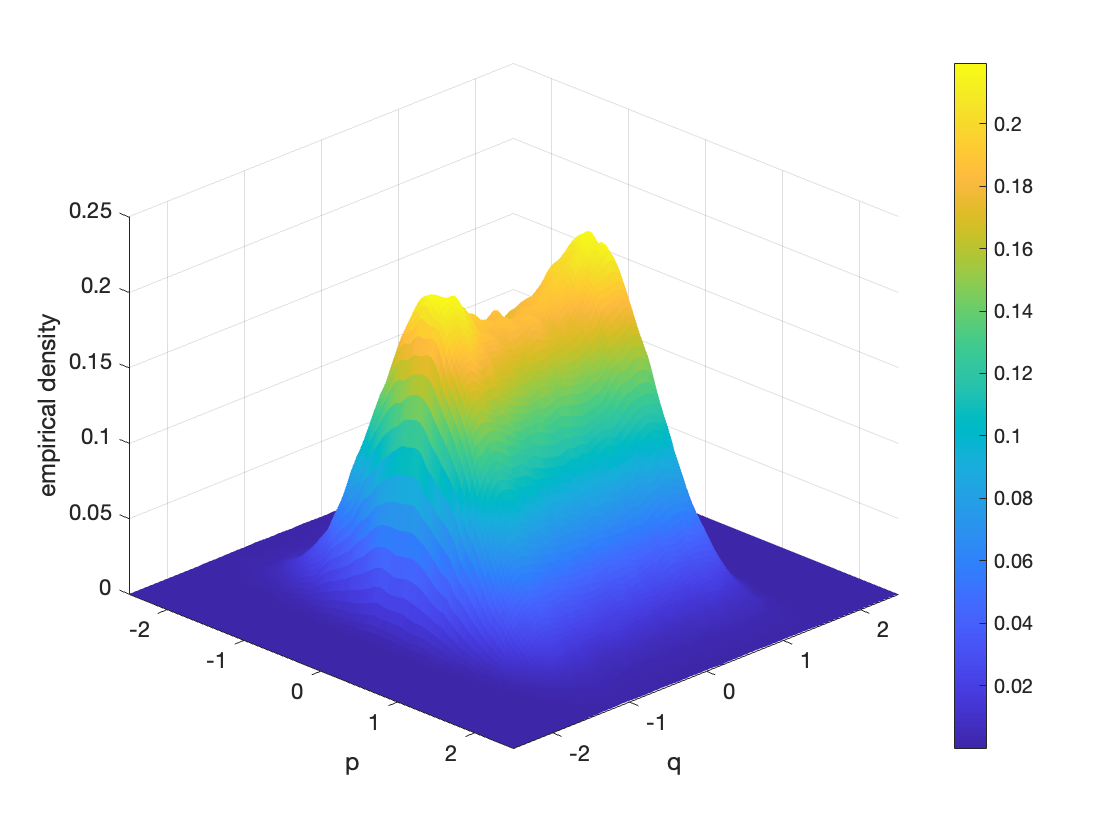}%
      }\hfill
      \subfloat[Gibbs density]{%
        \label{fig:ensembled}%
        \includegraphics[width=0.30\textwidth]
          {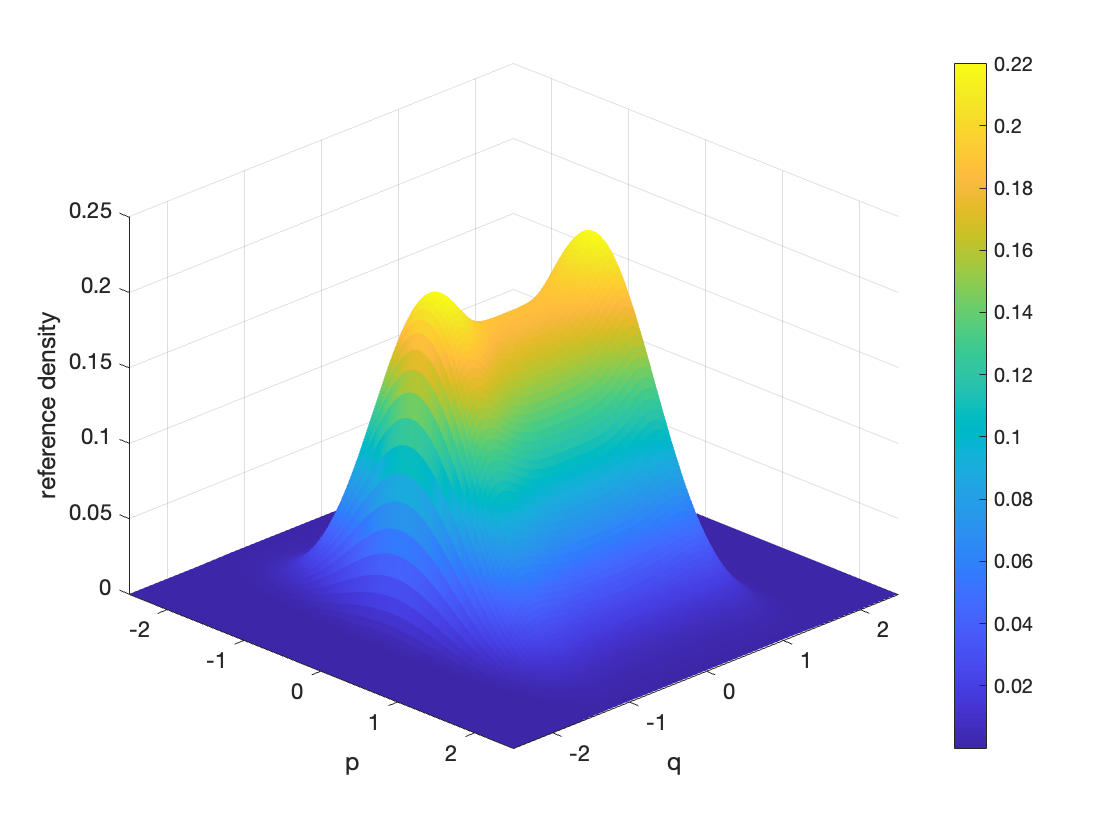}%
      }
      \caption{Simulation through the ensemble-averaging estimator}
      \label{ensemble}
    \end{figure}

	The ensemble-average estimator is $\hat{g}^{\text{erg}}_T:= \frac 1M \sum^M_{j=1} g(P^{h,j}_{K_h},Q^{h,j}_{K_h} )$, where $M$ denotes the number of sample paths for the Monte Carlo simulation, and $(P^{h,j}_{K_h},Q^{h,j}_{K_h} )$ denotes the numerical solution of the $j$-th trajectory at the final time step $K_h=T/h$. 
	We take the proposed conformal PRK method (see \eqref{eq:splittingPRK}) as a representative example, set $\sigma=1,\upsilon=1$, $d=1$,  $h=2^{-5}$, $M=10^6$, $R_0=1$, and employ a Gaussian kernel with the bandwidth $0.05$ as the test function $g$. 
	Setting the initial value to $(P(0),Q(0))^\top=(0,1)^\top$ and considering the time
	$T=1$ and  $100$, we illustrate the corresponding empirical measures in Figure~\ref{ensemble}(a)--(b) and the reference measure in Figure~\ref{ensemble}(c).
	We observe that the empirical densities asymptotically approach the target Gibbs density as $T$ increases, which is consistent with our theoretical result.

\begin{figure}[ht]
      \centering
      \subfloat[$T=1000$]{%
        \label{fig:time-averaginga}%
        \includegraphics[width=0.30\textwidth]
          {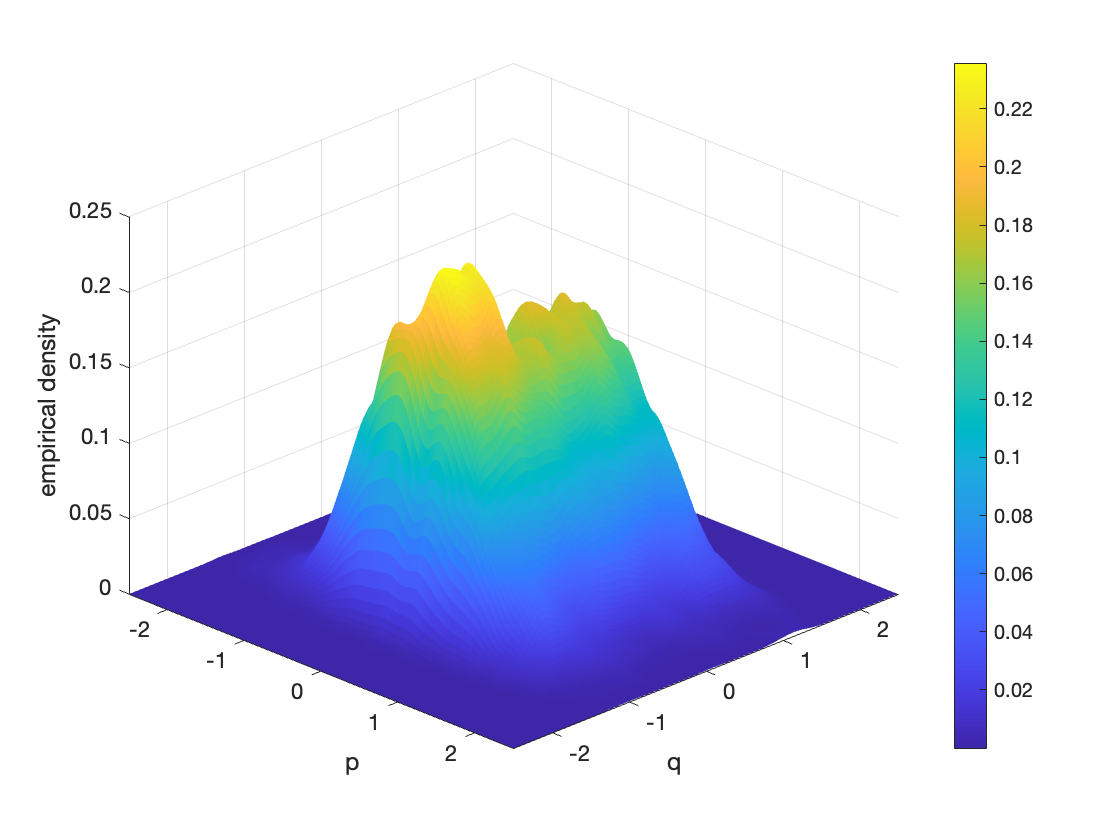}%
      }\hfill
      \subfloat[$T=10^5$]{%
        \label{fig:time-averagingb}%
        \includegraphics[width=0.30\textwidth]
          {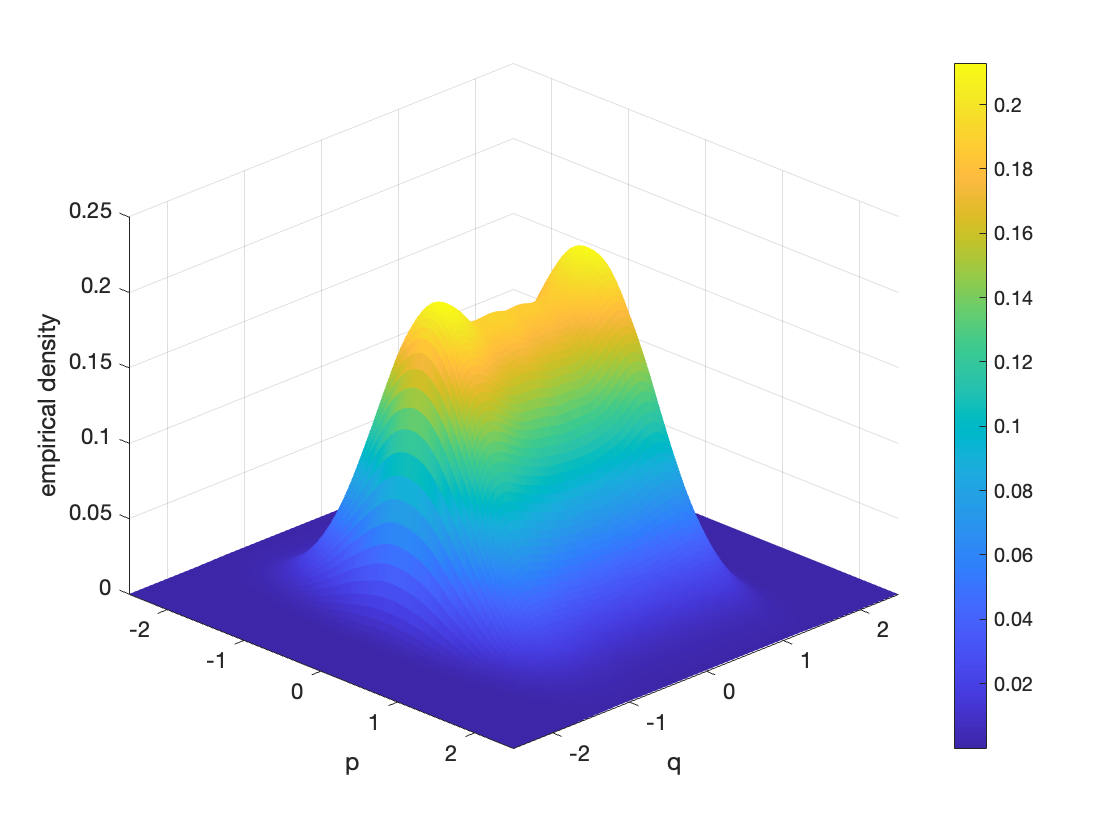}%
      }\hfill
      \subfloat[Gibbs density]{%
        \label{fig:time-averagingd}%
        \includegraphics[width=0.30\textwidth]
          {picture/density/Gibbs.png}%
      }
      \caption{Simulation through the time-averaging estimator}
      \label{time-averaging}
    \end{figure}

    The time-average estimator is defined by $\bar{g}_{T}^{\text{erg}}:= \frac {1}{K_h}\sum^{K_h}_{n=1} g(P^h_n, Q^h_n)$, where $K_h=T/h$.
	In this case, we still apply the proposed conformal PRK method,  set $\sigma=1,\upsilon=1$, $d=1$,  $h=2^{-5}$, $R_0=1$, and employ a Gaussian kernel with the bandwidth $0.05$ as the test function $g$. 
	Setting the initial value to $(P(0),Q(0))^\top=(0,1)^\top$ and considering the time
	$T=1000$ and $10^5$, we illustrate the corresponding empirical densities in Figure~\ref{time-averaging}(a)--(b) and the reference density in Figure~\ref{time-averaging}(c).
Both time-averaged densities are close to the Gibbs density. Here a single trajectory is used, which is never stopped or restarted. This is possible because the mollified schemes are ergodic and satisfy the strong law of large numbers \eqref {equ_LLN}. The rejection technique of \cite{MilsteinTretyakov2005SINUM, MilsteinTretyakov2007PhysD} does not allow this.

\subsection{Comparison of approaches to potential modification}
In this subsection, we compare two approaches to modifying the potential: the mollified potential approach introduced in the main text and the tamed potential approach presented in Appendix~\ref{suppsec_tamed} (see \eqref{Tame_equ_Uh}--\eqref{Tame_equ_fh} for their definitions). For both approaches, we employ the BAOAB method. In the tamed potential approach, we consider two choices of the parameter $\alpha$, namely $\alpha=1$ and $\alpha=2$.
For this comparison, we set $T=10$, $\sigma=1$, $d=2$,
$\upsilon=1$, and $R_0=2$, with initial values
$P(0)=(0,0)^\top$ and $Q(0)=(1,1)^\top$.  The step sizes are $h_i=2^{-i}$ for $i=3,4,5,6$. For the tamed potential approach with $\alpha=1$, the numerical estimate
$S_{U^h}\approx0.98$ gives $S_{U^h}<1=\upsilon$,
which supports Assumption~\ref{assm_v>CS} for the BAOAB method
with $a_{0,p,1}=-1/2$.
We choose the test function $g(p,q)=|q|^4$ to assess the effect
of the potential modification.
This choice is motivated by the fact that modifying the potential
affects only the position marginal of the invariant measure,
whose density is proportional to
$
    \exp\left(
        -\frac{2\upsilon}{\sigma^2}
        \left(\frac{1}{2}|p|^2+U^h(q)\right)
    \right).
$

We remark that the tamed potential approach may introduce a bias in the ergodic limit of the corresponding modified Langevin equation relative to that of the original Langevin equation. To illustrate this phenomenon, we consider the test function $g(p,q)=|q|^4$, for which numerical integration yields 
$\pi(g)\approx1.5026$. Let $\pi_{1,h}$ and $\pi_{2,h}$ denote
 the invariant measures of the modified equations
corresponding to the tamed potentials with $\alpha=1$ and $\alpha=2$, respectively.
The values of $\pi_{i,h}(g)$ and the corresponding errors
$|\pi_{i,h}(g)-\pi(g)|$, $i=1,2$, are reported
in Table~\ref{tab:ergodiclimit}.
These biases arise solely from the potential modification
and are independent of time discretization errors. By contrast, 
by Theorem~\ref{thm_weakconvergencetamed}~(\romannumeral 2) the bias of the mollified potential approach is at most $C\exp(-cR^{r_1}(h))$, because $U^h$ differs from $U$ only for $|q|>R(h )$. 
For the present example, numerical integration gives $|\pi_h(g)-\pi(g)|<10^{-50}$ for $h\le2^{-3} $, so this bias is negligible.

\begin{table}[htbp]
    \centering
    \caption{Ergodic limits and associated biases.}
    \label{tab:ergodiclimit}
    \begin{tabular}{ccccc}
        \hline
        Quantity & $h=2^{-3}$ & $h=2^{-4}$ & $h=2^{-5}$ & $h=2^{-6}$ \\
        \hline
        $\pi_{1,h}(g)$
        & 2.0606 & 1.7447 & 1.6156 & 1.5572 \\
        $|\pi_{1,h}(g)-\pi(g)|$
        & 0.5580 & 0.2421 & 0.1129 & 0.0546 \\
        $\pi_{2,h}(g)$
        & 1.5572 & 1.5159 & 1.5059 & 1.5034 \\
        $|\pi_{2,h}(g)-\pi(g)|$
        & 0.0546 & 0.0133 & $3.3\times10^{-3}$ & $8.3\times10^{-4}$ \\
        \hline
    \end{tabular}
\end{table}

    \begin{figure}[ht]
    \centerline{\includegraphics[scale=0.4]{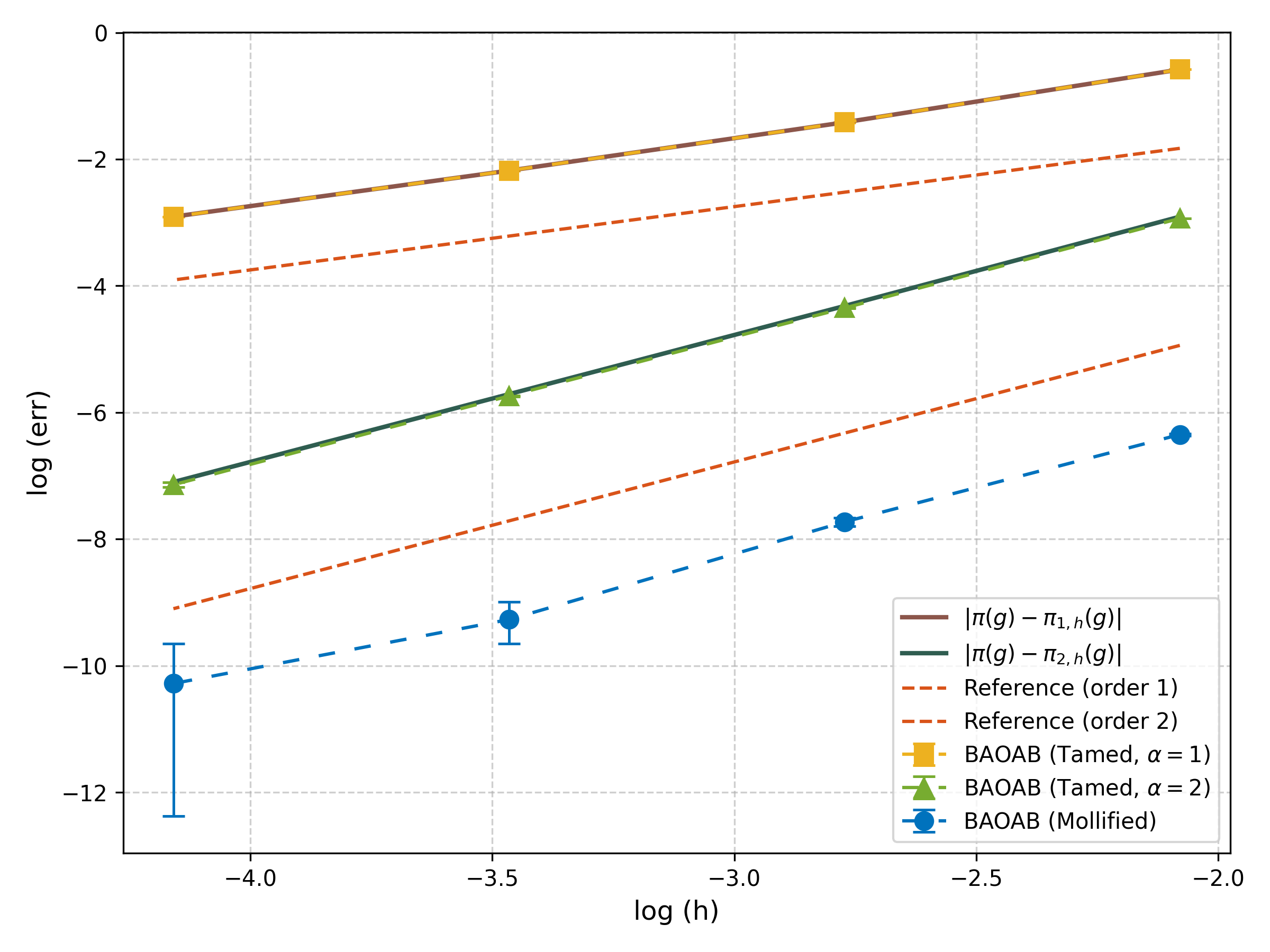}} 
        \caption{Comparison of modifying approaches}
		\label{fig:comparison}
	\end{figure}
In Figure~\ref{fig:comparison}, we employ $M=10^{10}$ sample paths for the Monte Carlo simulation and plot the errors of the BAOAB
methods using the mollified and tamed potentials relative
to $\pi(g)$, together with the taming biases
$|\pi_{i,h}(g)-\pi(g)|$, $i=1,2$, as reference curves.
The results indicate second-order convergence for the mollified
potential approach and the tamed potential approach with $\alpha=2$,
and first-order convergence for the tamed potential approach
with $\alpha=1$. No truncations are observed in the simulations
for $h=2^{-3},\ldots,2^{-6}$,
i.e., here on every simulated path the mollified BAOAB method coincides with the BAOAB method 
applied to the original potential and its error is the discretization error of BAOAB alone.
It is observed from Figure~\ref{fig:comparison} that the numerical method with the mollified potential approach gives the smallest errors over the tested range of step sizes, and the errors of the tamed potential approach are largely due to the inherent biases in the ergodic limits.
Note that with $\alpha=2$, the tamed potential approach also exhibits second-order convergence in this example, although our analysis does not cover it. However, the bias of the tamed potential is of order $h^{\alpha}$ (see Theorem~\ref{Tamed_thm_weakconvergencetamed}) and the analysis of the numerical methods there requires $\alpha\le1$. Therefore, taming cannot provide provably second-order methods, while the bias of the mollified potential, being exponentially small, does not restrict the order.

    \subsection{Conformal symplecticity}
	We also present numerical experiments to examine the preservation of the conformal symplecticity. We set $\sigma=1,\upsilon=0.5, d=1$,   $h=2^{-3}$, $R_0=1$, and choose initial values on the unit circle $p^2+q^2=1$. 
\begin{figure}[ht]
  \centering
  \subfloat[Conformal symplecticity]{%
    \label{fig:conform_sym2}%
    \includegraphics[width=0.42\textwidth]
      {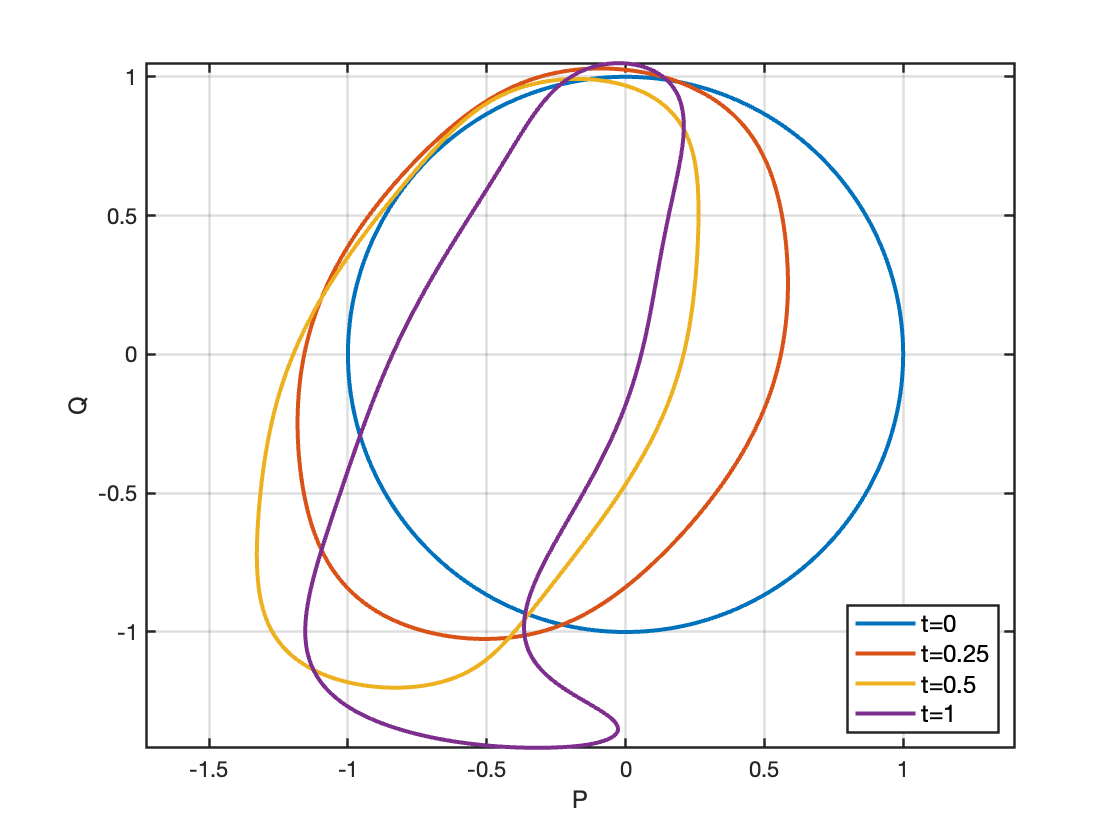}%
  }\hspace{15mm}
  \subfloat[Comparison]{%
    \label{fig:conform_compare}%
    \includegraphics[width=0.42\textwidth]
      {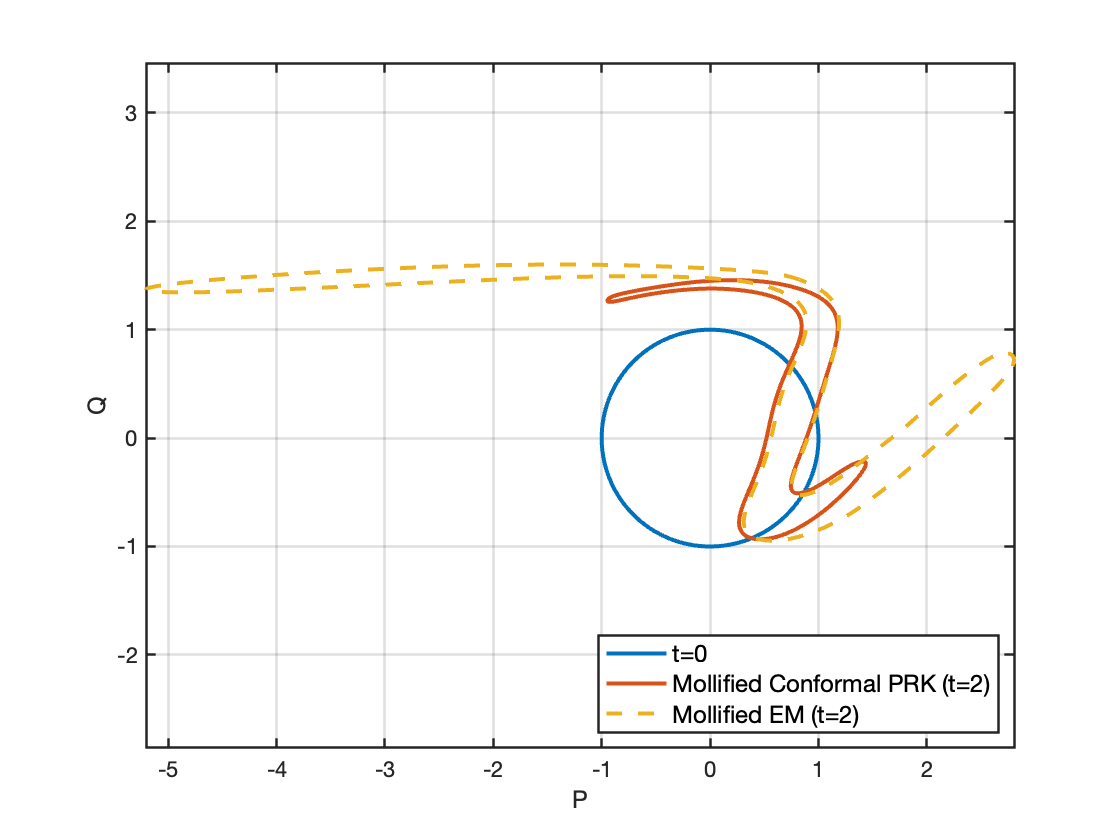}%
  }
  \caption{Evolution for area in the phase space, $\upsilon=0.5$.}
  \label{fig:conformEM}
\end{figure}
    In Figure~\ref{fig:conformEM}(a), we plot the curve of the proposed conformal PRK method (see \eqref{eq:splittingPRK}) in the phase plane at different time points $t=0, \frac{1}{4}, \frac 12$, and $1$.  %
	In Figure~\ref{fig:conformEM}(b), we compare the curves at $t=1$ for the proposed conformal PRK method and the Euler--Maruyama method applied to the proposed mollified approach. 
    
    To quantify the preservation of conformal symplecticity, we report the normalized phase-space areas $e^{\upsilon t}\frac{A_h(t)}{A_0}$ in Table~\ref{tab:conformal_areaEM}, which theoretically equal $1$ if the conformal symplecticity is preserved. Here, $A_h(t)$ is approximated by the area of the polygon formed by the evolved boundary points, and $A_0$ is the area of the initial polygon.
    \begin{table}[htbp]
    \centering
    \caption{Evolution of normalized phase-space areas}
    \label{tab:conformal_areaEM}
    \begin{tabular}{ccccc}
        \hline
        Method & t=0.5 & t=1 & t=2 & t=3 \\
        \hline
        Conformal PRK & 1.000 & 1.000 & 1.000 & 1.000 \\
        Euler--Maruyama   & 0.986 & 1.015 & 2.003 & 5.025\\
        \hline
    \end{tabular}
\end{table}
    It is observed in Table~\ref{tab:conformal_areaEM} that the conformal PRK method preserves the underlying conformal symplecticity, while the Euler--Maruyama method fails.  

    \section*{Acknowledgments}
    CC, XC and JH were supported by the National Key R\&D Program of China under Grant (No. 2024YFA1015900),  by the National Natural Science Foundation of China (No. 12031020, No. 12461160278, and No. 12471386), and by Youth Innovation Promotion Association CAS.
    MVT was partially supported by the London Mathematical Society under Scheme 5 (Grant No. 52423), and he is very grateful for the hospitality of ICMSEC, Academy of Mathematics and Systems Science, Chinese Academy of Sciences (Beijing), where this work was initiated during MVT's visit.

\bibliographystyle{plain}
\bibliography{references}

@book {HongWangbook2019,
	AUTHOR = {Hong, J. and Wang, X.},
	TITLE = {Invariant {M}easures for {S}tochastic {N}onlinear {S}chr\"odinger
	{E}quations: {N}umerical {A}pproximations and {S}ymplectic {S}tructures},
	SERIES = {Lecture Notes in Mathematics},
	VOLUME = {2251},
	PUBLISHER = {Springer},
	YEAR = {2019}
}

@article {CDHZ2025SINUM,
	AUTHOR = {Chen, C. and Dang, T. and Hong, J. and Zhang,
	F.},
	TITLE = {A new class of splitting methods that preserve ergodicity and
	exponential integrability for the stochastic {L}angevin
	equation},
	JOURNAL = {SIAM J. Numer. Anal.},
	FJOURNAL = {SIAM Journal on Numerical Analysis},
	VOLUME = {63},
	YEAR = {2025},
	NUMBER = {2},
	PAGES = {1000--1024},
	ISSN = {0036-1429,1095-7170},
	MRCLASS = {65C30 (37M25 60H35)},
	MRNUMBER = {4898691},
}

@article {MilsteinTretyakov2003Quasisym,
	AUTHOR = {Milstein, G. N. and Tretyakov, M. V.},
	TITLE = {Quasi-symplectic methods for {L}angevin-type equations},
	JOURNAL = {IMA J. Numer. Anal.},
	FJOURNAL = {IMA Journal of Numerical Analysis},
	VOLUME = {23},
	YEAR = {2003},
	NUMBER = {4},
	PAGES = {593--626}
	}

@book{Coffey2012book,
	title={The Langevin equation: with Applications to Stochastic Problems in Physics, Chemistry and Electrical Engineering},
	author={Coffey, W. and Kalmykov, Y. P.},
	volume={27},
	year={2012},
	publisher={World Scientific}
}

@misc{Dai2025arXiv,
	title={An explicit splitting {SAV} scheme for the kinetic {L}angevin dynamics}, 
	author={L. Dai and Y. Jiang and X. Wang},
	year={2025},
     howpublished={arXiv:2509.04251},
	eprint={2509.04251},
	archivePrefix={arXiv},
	primaryClass={math.NA}
}

@article {MilsteinTretyakov2005SINUM,
	AUTHOR = {Milstein, G. N. and Tretyakov, M. V.},
	TITLE = {Numerical integration of stochastic differential equations
	with nonglobally {L}ipschitz coefficients},
	JOURNAL = {SIAM J. Numer. Anal.},
	FJOURNAL = {SIAM Journal on Numerical Analysis},
	VOLUME = {43},
	YEAR = {2005},
	NUMBER = {3},
	PAGES = {1139--1154},
	ISSN = {0036-1429,1095-7170},
	MRCLASS = {60H35 (60H10 65C30)},
	MRNUMBER = {2177799}
}

@article {MilsteinTretyakov2007PhysD,
	AUTHOR = {Milstein, G. N. and Tretyakov, M. V.},
	TITLE = {Computing ergodic limits for {L}angevin equations},
	JOURNAL = {Phys. D},
	FJOURNAL = {Physica D. Nonlinear Phenomena},
	VOLUME = {229},
	YEAR = {2007},
	NUMBER = {1},
	PAGES = {81--95},
	ISSN = {0167-2789,1872-8022},
	MRCLASS = {82C31 (60H10)},
	MRNUMBER = {2340180},
	MRREVIEWER = {Peter\ E.\ Kloeden}
}

@article {Wang2025,
	AUTHOR = {Zhao, Y. and Wang, X. and Zhang, Z.},
	TITLE = {Weak error analysis for strong approximation schemes of {SDE}s
	with super-linear coefficients {II}: {F}initely many bounded
	moments and higher-order schemes},
	JOURNAL = {J. Sci. Comput.},
	FJOURNAL = {Journal of Scientific Computing},
	VOLUME = {105},
	YEAR = {2025},
	NUMBER = {2},
	EID = {58},
	ISSN = {0885-7474,1573-7691},
	MRCLASS = {65C30 (60H10 60H35)},
	MRNUMBER = {4973065}
}

@book {Milstein2004book,
	AUTHOR = {Milstein, G. N. and Tretyakov, M. V.},
	TITLE = {Stochastic {N}umerics for {M}athematical {P}hysics},
    edition={2nd},
	PUBLISHER = {Springer},
	YEAR = {2021}
}

@article {Hutzenthaler2011EulerDiverge,
	AUTHOR = {Hutzenthaler, M. and Jentzen, A. and Kloeden, P.
	E.},
	TITLE = {Strong and weak divergence in finite time of {E}uler's method
	for stochastic differential equations with non-globally
	{L}ipschitz continuous coefficients},
	JOURNAL = {Proc. R. Soc. Lond. Ser. A Math. Phys. Eng. Sci.},
	FJOURNAL = {Proceedings of The Royal Society of London. Series A.
	Mathematical, Physical and Engineering Sciences},
	VOLUME = {467},
	YEAR = {2011},
	NUMBER = {2130},
	PAGES = {1563--1576},
	ISSN = {1364-5021,1471-2946},
	MRCLASS = {65C30 (60H10 60H35)},
	MRNUMBER = {2795791},
	MRREVIEWER = {Henri\ Schurz}
}

@article {CHS2022,
    AUTHOR = {Cui, J. and Hong, J. and Sheng, D.},
     TITLE = {Density function of numerical solution of splitting {AVF}
              scheme for stochastic {L}angevin equation},
   JOURNAL = {Math. Comp.},
  FJOURNAL = {Mathematics of Computation},
    VOLUME = {91},
      YEAR = {2022},
    NUMBER = {337},
     PAGES = {2283--2333},
      ISSN = {0025-5718,1088-6842},
   MRCLASS = {60H10 (60H07 60H35)},
  MRNUMBER = {4451463},
MRREVIEWER = {Kai\ He}
}

@book {Hongsun2022,
	AUTHOR = {Hong, J. and Sun, L.},
	TITLE = {Symplectic {I}ntegration of {S}tochastic {H}amiltonian {S}ystems},
	SERIES = {Lecture Notes in Mathematics},
	VOLUME = {2314},
	PUBLISHER = {Springer, Singapore},
	YEAR = {2022},
	PAGES = {xii+298},
	ISBN = {978-981-19-7669-8; 978-981-19-7670-4},
	MRCLASS = {65P10},
	MRNUMBER = {4628103}
}

@article {Sun2017SINUM,
	AUTHOR = {Hong, J. and Sun, L. and Wang, X.},
	TITLE = {High order conformal symplectic and ergodic schemes for the
	stochastic {L}angevin equation via generating functions},
	JOURNAL = {SIAM J. Numer. Anal.},
	FJOURNAL = {SIAM Journal on Numerical Analysis},
	VOLUME = {55},
	YEAR = {2017},
	NUMBER = {6},
	PAGES = {3006--3029},
	ISSN = {0036-1429,1095-7170},
	MRCLASS = {60H35 (65C30 65P10)},
	MRNUMBER = {3730544}
}

@article {Owhadi2010SINUM,
	AUTHOR = {Bou-Rabee, N. and Owhadi, H.},
	TITLE = {Long-run accuracy of variational integrators in the stochastic
	context},
	JOURNAL = {SIAM J. Numer. Anal.},
	FJOURNAL = {SIAM Journal on Numerical Analysis},
	VOLUME = {48},
	YEAR = {2010},
	NUMBER = {1},
	PAGES = {278--297},
	ISSN = {0036-1429,1095-7170},
	MRCLASS = {65C30 (60H10 60J05 65C05 65P10)},
	MRNUMBER = {2608370},
	MRREVIEWER = {Andreas\ R\"o\ss ler}
}

@article {Sabanis2016AAP,
	AUTHOR = {Sabanis, S.},
	TITLE = {Euler approximations with varying coefficients: the case of
	superlinearly growing diffusion coefficients},
	JOURNAL = {Ann. Appl. Probab.},
	FJOURNAL = {The Annals of Applied Probability},
	VOLUME = {26},
	YEAR = {2016},
	NUMBER = {4},
	PAGES = {2083--2105},
	ISSN = {1050-5164,2168-8737},
	MRCLASS = {60H35 (65C30)},
	MRNUMBER = {3543890},
	MRREVIEWER = {Minoo\ Kamrani}
    }

@article {Jentzen2012AAP,
	AUTHOR = {Hutzenthaler, M. and Jentzen, A. and Kloeden, P.
	E.},
	TITLE = {Strong convergence of an explicit numerical method for {SDE}s
	with nonglobally {L}ipschitz continuous coefficients},
	JOURNAL = {Ann. Appl. Probab.},
	FJOURNAL = {The Annals of Applied Probability},
	VOLUME = {22},
	YEAR = {2012},
	NUMBER = {4},
	PAGES = {1611--1641},
	ISSN = {1050-5164,2168-8737},
	MRCLASS = {65C30 (60H35)},
	MRNUMBER = {2985171},
	MRREVIEWER = {Andreas\ R\"o\ss ler},
	}

@article {Talay2002,
	AUTHOR = {Talay, D.},
	TITLE = {Stochastic {H}amiltonian systems: exponential convergence to
	the invariant measure, and discretization by the implicit
	{E}uler scheme},
	JOURNAL = {Markov Process. Related Fields},
	FJOURNAL = {Markov Processes and Related Fields},
	VOLUME = {8},
	YEAR = {2002},
	NUMBER = {2},
	PAGES = {163--198},
	ISSN = {1024-2953},
	MRCLASS = {60H10 (37J99 60H35)},
	MRNUMBER = {1924934},
	MRREVIEWER = {Jan\ I.\ Seidler},
}

@article {Mattingly2002SPA,
	AUTHOR = {Mattingly, J. C. and Stuart, A. M. and Higham, D. J.},
	TITLE = {Ergodicity for {SDE}s and approximations: locally {L}ipschitz
	vector fields and degenerate noise},
	JOURNAL = {Stochastic Process. Appl.},
	FJOURNAL = {Stochastic Processes and their Applications},
	VOLUME = {101},
	YEAR = {2002},
	NUMBER = {2},
	PAGES = {185--232},
	ISSN = {0304-4149,1879-209X},
	MRCLASS = {60H10 (34F05 37H10 60J35)},
	MRNUMBER = {1931266},
	MRREVIEWER = {Jan\ I.\ Seidler}
}

@book{leimkuhler_mathews_15,
  publisher =   "Springer",
  year      =   "2015",
title ="Molecular dynamics with deterministic and stochastic numerical methods",
author = "Leimkuhler, B. and  Matthews, C."
}

@article{TZ13,
  author    =   "Tretyakov, M. V. and Zhang, Z." ,
  title     =   "A fundamental mean-square convergence theorem for {SDE}s with locally {L}ipschitz coefficients and its applications",
   journal  =   "SIAM J. Numer. Anal.",
      volume ="51",
     pages  =   "3135--3162",
  year      =   "2013"}

@article{BenM13,
	Author = {B. Leimkuhler and C. Matthews},
	Journal = {Appl. Math. Res. Express},
	Number = {1},
	Pages = {34--56},
	Title = {Rational construction of stochastic numerical methods for molecular sampling},
	Volume = {2013},
	Year = {2013}}

@article{MRT02,
    author  =   "G. N. Milstein and Yu. M. Repin and M. V. Tretyakov",
     title  =   "Symplectic integration of {H}amiltonian systems with additive noise",
   journal  =   "SIAM J. Numer. Anal.",
      year  =   "2002",
      volume ="39",
     pages  =   "2066--2088"}

@article{MRT03,
    author  =   "G. N. Milstein and Yu. M. Repin and M. V. Tretyakov",
     title  =   "Numerical methods for stochastic systems preserving symplectic structure",
   journal  =   "SIAM J. Numer. Anal.",
      year  =   "2003",
      volume ="40",
     pages  =   "1583--1604"}

@book{HLW02,
  author    =   "E. Hairer and C. Lubich and G. Wanner" ,
  title     =   "Geometric numerical integration: structure preserving algorithms for ordinary differential equations",
   publisher =   "Springer",
  address ="Berlin",
  year      =   "2002"}

@article{Handy24,
    author = {Wang, X. and Zhao, Y. and Zhang, Z.},
    title = {Weak error analysis for strong approximation schemes of {SDE}s with super-linear coefficients},
    journal = {IMA Journal of Numerical Analysis},
    volume = {44},
    number = {5},
    pages = {3153-3185},
    year = {2024}
    }

@book{HJ15,
  author    =   "M. Hutzenthaler and A. Jentzen" ,
  title     =   "Numerical approximations of stochastic differential equations with non-globally {L}ipschitz continuous coefficients",
   publisher =   "AMS",
   series ="Mem. Amer. Math. Soc.",
   volume = "236",
  address = "Providence",
  year      =   "2015"}

@article{KL18,
  author    =   "C. Kelly and G. J. Lord" ,
  title     =   "Adaptive time-stepping strategies for nonlinear stochastic systems",
   journal  =   "IMA J. Numer. Anal.",
      volume ="38",
     pages  =   "1523--1549",
  year      =   "2018"}

@article{FG19,
  author    =   "W. Fang and M. B. Giles" ,
  title     =   "Adaptive {E}uler-{M}aruyama method for {SDE}s with non-globally {L}ipschitz drift",
   journal  =   "Ann. Appl. Probab.",
   volume ="30",
     pages  =   "526--560",
  year      =   "2020"}

@book{HAS80,
  author    =   "R. Z. Khasminskii" ,
  title     =   "Stochastic stability of differential equations",
   publisher =   "Springer",
    edition={2nd},
  year      =   "2012"}

@incollection{Neal2011,
  author    = {Neal, R. M.},
  title     = {{MCMC} using {H}amiltonian dynamics},
  booktitle = {Handbook of Markov Chain Monte Carlo},
  editor    = {Brooks, S. and Gelman, A. and Jones, G. L. and Meng, X.-L.},
  publisher = {Chapman \& Hall/CRC},
  pages     = {113--162},
  year      = {2011}
}

@inproceedings{ChenFoxGuestrin2014,
  author    = {Chen, T. and Fox, E. and Guestrin, C.},
  title     = {Stochastic gradient {H}amiltonian {M}onte {C}arlo},
  booktitle = {Proc. 31st Intern. Confer. Machine Learn.},
  volume    = {32},
  pages     = {1683--1691},
  year      = {2014}
}

@book{MeynTweedie2009,
author = {Meyn, S. and Tweedie, R. L.},
title = {Markov Chains and Stochastic Stability},
edition = {2nd},
publisher = {Cambridge University Press},
address = {Cambridge},
year = {2009}
}

\appendix

\section{Proofs for the mollified Langevin equations}\label{suppsec_proofmollifiedequ}
\subsection{Proofs of  Propositions~\ref{prop:fh-fS2} and \ref{prop_ergodicfhS2}}
\label{sec_proofpart1}

\begin{proof}[Proof of Proposition~\ref{prop:fh-fS2}]
  		Throughout this proof, $C>0$ denotes a constant that may change from line to
line but is independent of $h$ and $q$ (and, in parts (\romannumeral1) and (%
\romannumeral4), of $k$); it depends only on $d$, the exponents $%
r_{1},\ldots,r_{3}$, the constants $L_{1},\ldots ,L_{4}$ and $C$ of Assumption~\ref{assm_CandPoly}, 
the polynomial growth of $D^3U$, $R_0$, and $\Vert \theta ^{\prime }\Vert _{\infty
},\Vert \theta ^{\prime \prime }\Vert _{\infty }, \Vert \theta ^{\prime \prime \prime }\Vert _{\infty }$. 
We abbreviate $R=R(h)$, $%
\lambda =\lambda (h)$, $\tilde{\theta}=\tilde{\theta}(q)=\theta (|q|/R)$ and 
$\hat{q}=q/|q|$, and we repeatedly use the scaling identities 
\begin{align}
R^{r_{2}-1}=R_0^{r_{2}-1} h^{-1/2},\ \ \ \lambda
=2^{r_{2}-1}L_{2}\,R^{r_{2}-1}=L_{2}(2R_0)^{r_2-1}h^{-1/2},\ \
  \label{equ_scaling}\\
h^{k}=R_0^{2k(r_{2}-1)} R^{-2k(r_{2}-1)},\quad k>0. \notag
\end{align}%
Also, note $R>1$ and hence $|q|\geq R$ implies $|q|>1$. Recall the three
regions determined by the definition of $\theta $: for $|q|\leq R$ we have $%
\tilde{\theta}=1$ and $U^{h}=U$; for $|q|\geq 2R$ we have $\tilde{\theta}=0$
and $U^{h}(q)=\lambda |q|^{2}$; and $\nabla \tilde{\theta},$ $D^{2}\tilde{\theta%
}$ are equal $0$ for $q$ outside $R<|q|<2R$.

To prove the statements of Proposition~\ref{prop:fh-fS2}, it is 
useful to first establish some auxiliary estimates. 
For $|q|\geq R,$ Assumption~\ref{assm_CandPoly} gives $|\nabla U(q)|\leq L_{3}|q|^{r_{2}}+C\leq
C|q|^{r_{2}}$ as $|q|>1$ and \eqref{equ_scaling} implies $2\lambda
|q|=2^{r_{2}}L_{2}R^{r_{2}-1}|q|\leq 2^{r_{2}}L_{2}|q|^{r_{2}}$, which together lead to
\begin{equation}
|\nabla U(q)-2\lambda q|\leq C|q|^{r_{2}}.  \label{equ_E1}
\end{equation}%
Next, let $R\leq |q|\leq 2R.$ Assumption~\ref{assm_CandPoly} gives $U\geq 0$
and $U(q)\leq L_{2}|q|^{r_{1}}+C=L_{2}|q|^{r_{2}+1}+C\leq C|q|^{r_{2}+1}$,
while $\lambda |q|^{2}=2^{r_{2}-1}L_{2}R^{r_{2}-1}|q|^{2}\leq C|q|^{r_{2}+1}$
by $R\leq |q|.$ Therefore, by the triangle inequality, we get%
\begin{equation}
|U(q)-\lambda |q|^{2}|\leq C|q|^{r_{2}+1}.  \label{equ_E2}
\end{equation}

Finally, from $\nabla \tilde \theta=\frac{\theta ^{\prime}(|q|/R)\hat{q}}{R}$, $D^{2}\tilde{\theta}=\frac{\theta ^{\prime \prime }(|q|/R)}{%
R^{2}}\,\hat{q}\hat{q}^{\top }+\frac{\theta ^{\prime }(|q|/R)}{R|q|}\big(I-%
\hat{q}\hat{q}^{\top }\big)$, and $R\geq |q|/2$ on the support of $D^{2}%
\tilde{\theta},$ we obtain 
\begin{equation}
|\nabla \tilde{\theta}|\leq \frac{\Vert \theta ^{\prime }\Vert _{\infty }}{R}%
,\qquad \Vert D^{2}\tilde{\theta}\Vert _{\text{op}}\leq \frac{\Vert \theta
^{\prime \prime }\Vert _{\infty }}{R^{2}}+\frac{\Vert \theta ^{\prime }\Vert
_{\infty }}{R|q|}\leq C|q|^{-2}.  \label{equ_D2phi}
\end{equation}

Introduce $F$ to be scalar-, vector-, or matrix-valued with $F(q)=0$ for $%
|q|\leq R$ and $\Vert F(q)\Vert \leq C_{0}|q|^{m}$ for $|q|>R$, where $m\geq
0$ and $C_{0}$ is independent of $h,$ $q$ and $k$. Since $2k(r_{2}-1)>0$
with  $k>0$ and for $|q|\geq R$ we have $R/|q|\leq 1,$ we get $1\leq
(|q|/R)^{2k(r_{2}-1)}.$ Therefore, by \eqref{equ_scaling}, we arrive at 
\begin{equation}
\Vert F(q)\Vert \leq C_{0}|q|^{m}\Big(\tfrac{|q|}{R}\Big)%
^{2k(r_{2}-1)}=C_{0}\,R^{-2k(r_{2}-1)}|q|^{\,m+2k(r_{2}-1)} \leq C_{0}%
\,h^{k}|q|^{\,m+2k(r_{2}-1)},  \label{eq:F}
\end{equation}%
where for $|q|<R$ the left-hand side vanishes.

Now we move to proving the statements of Proposition~\ref{prop:fh-fS2}.

\noindent \textbf{(\romannumeral2)} From the expression for the mollified
force \eqref{eq:fhS2}, we obtain%
\begin{equation*}
|f^{h}(q)|\leq \tilde{\theta}\,|\nabla U(q)|+(1-\tilde{\theta})\,2\lambda
|q|+|U(q)-\lambda |q|^{2}|\,|\nabla \tilde{\theta}|.
\end{equation*}%
The first term is bounded from above by $C(1+|q|^{r_{2}})$ by 
Assumption~\ref{assm_CandPoly}. The second term is supported in $|q|\geq R$, where $%
2\lambda |q|\leq C|q|^{r_{2}}$. The third is supported in $R\leq |q|\leq 2R$%
, where, by \eqref{equ_E2} and \eqref{equ_D2phi}, $|U(q)-\lambda
|q|^{2}||\nabla \tilde{\theta}|\leq C|q|^{r_{2}+1}\cdot \frac{\Vert \theta
^{\prime }\Vert _{\infty }}{R}\leq C|q|^{r_{2}}$ using $R\geq |q|/2$. Hence $%
|f^{h}(q)|\leq C(1+|q|^{r_{2}})$. 

For $|q|\leq R$, $|q|^{r_{2}}=|q|^{r_{2}-1}|q|\leq R^{r_{2}-1}|q|$ together
with $1\leq R^{r_{2}-1}$ yields $|f^{h}(q)|\leq C(1+|q|^{r_{2}})\leq
CR^{r_{2}-1}(1+|q|)=Ch^{-1/2}(1+|q|)$ using \eqref{equ_scaling}. For $R\leq
|q|\leq 2R$, $|q|^{r_{2}}\leq (2R)^{r_{2}-1}|q|=CR^{r_{2}-1}|q|$, so again $%
|f^{h}(q)|\leq CR^{r_{2}-1}(1+|q|)$. For $|q|\geq 2R$, using %
\eqref{equ_scaling}, we get $|f^{h}(q)|=2\lambda
|q|=2^{r_{2}}L_{2}R^{r_{2}-1}|q|\leq CR^{r_{2}-1}(1+|q|)$. 

Both bounds hold for every $q$, whence we arrive at (\romannumeral2). 

\noindent \textbf{(\romannumeral3)} Differentiating the expression for the mollified
force \eqref{eq:fhS2} and grouping the two cross-terms gives the symmetric
expression 
\begin{equation}
D^{2}U^{h}=\tilde{\theta}\,D^{2}U+(1-\tilde{\theta})\,2\lambda I+\big(\nabla
U-2\lambda q\big)\otimes \nabla \tilde{\theta}+\nabla \tilde{\theta}\otimes %
\big(\nabla U-2\lambda q\big)+\big(U-\lambda |q|^{2}\big)D^{2}\tilde{\theta}.
\label{equ_D2Uh_clean}
\end{equation}%
For $|q|\leq R$ we have $\Vert D^{2}U^{h}(q)\Vert _{\text{op}}=\Vert
D^{2}U(q)\Vert _{\text{op}}\leq L_{4}|q|^{r_{3}}+C\leq C(1+|q|^{r_{2}-1})$. For 
$|q|\geq $ $R$, the first term in \eqref{equ_D2Uh_clean} $\Vert \tilde{\theta%
}\,D^{2}U(q)\Vert _{\text{op}}\leq C|q|^{r_{3}}+C\leq C|q|^{r_{2}-1}$; the
second is bounded by $2\lambda =2^{r_{2}}L_{2}R^{r_{2}-1}\leq C|q|^{r_{2}-1}$
as $R\leq |q|$. Both cross-terms are non-zero only for $R\leq |q|\leq 2R$
and $|\nabla U(q)-2\lambda q|\,|\nabla \tilde{\theta}|\leq C|q|^{r_{2}}\cdot 
\frac{\Vert \theta ^{\prime }\Vert _{\infty }}{R}\leq C|q|^{r_{2}-1}$ thanks
to \eqref{equ_E1}, \eqref{equ_D2phi} and $R\geq |q|/2$. The last term $%
|U(q)-\lambda |q|^{2}|\,\Vert D^{2}\tilde{\theta}\Vert _{\text{op}}\leq
C|q|^{r_{2}+1}\cdot C|q|^{-2}=C|q|^{r_{2}-1}$ thanks to \eqref{equ_E2} and %
\eqref{equ_D2phi}. Thus $\Vert D^{2}U^{h}(q)\Vert _{\text{op}}\leq
C(1+|q|^{r_{2}-1})$ everywhere.

For $|q|\leq R$, $\Vert D^{2}U^{h}(q)\Vert _{\text{op}}=\Vert D^{2}U(q)\Vert _{%
\text{op}}\leq L_{4}|q|^{r_{3}}+C\leq L_{4}R^{r_{2}-1}+C\leq Ch^{-1/2}$
using $R^{r_{2}-1}=R_0^{r_{2}-1} h^{-1/2}\geq 1$. For $R\leq |q|\leq 2R$, the bound $\Vert
D^{2}U^{h}(q)\Vert _{\text{op}}\leq C|q|^{r_{2}-1}$ just proved and $|q|\leq 2R$
give $\Vert D^{2}U^{h}(q)\Vert _{\text{op}}\leq C(2R)^{r_{2}-1}=Ch^{-1/2}$. For 
$|q|\geq 2R$, $\Vert D^{2}U^{h}(q)\Vert _{\text{op}}=2\lambda \leq Ch^{-1/2}$
by \eqref{equ_scaling}. 

Both bounds hold at every $q$, whence we arrived at (\romannumeral3).

\noindent \textbf{(\romannumeral1)} Since $U^{h}=U$ for $|q|\leq R$, the
difference $f-f^{h}$ vanishes there. For $|q|\geq R$, \eqref{eq:fhS2} gives 
$f(q)-f^{h}(q)=(1-\tilde{\theta})\big(\nabla U(q)-2\lambda q\big)-\big(%
U(q)-\lambda |q|^{2}\big)\nabla \tilde{\theta}.$
The first term is bounded from above by $C|q|^{r_{2}}$ due to \eqref{equ_E1}%
; the second is supported in $R\leq |q|\leq 2R$, where thanks to  %
\eqref{equ_E2} and \eqref{equ_D2phi} it is bounded from above by $%
C|q|^{r_{2}+1}\cdot \frac{\Vert \theta ^{\prime }\Vert _{\infty }}{R}\leq
C|q|^{r_{2}}$. Hence $|f(q)-f^{h}(q)|\leq C|q|^{r_{2}}$ for $|q|\geq R$, and
it vanishes for $|q|<R$. Applying \eqref{eq:F} with $m=r_{2}$ yields (i).

\noindent \textbf{(\romannumeral4)} Since $U^{h}=U$ for $|q|\leq R$, the
difference $D^{2}U-D^{2}U^{h}$ vanishes there. For $|q|\geq R$, using $\Vert
D^{2}U(q)\Vert _{\text{op}}\leq C|q|^{r_{2}-1}$ thanks to $r_{3}\leq r_{2}-1$
and $|q|>1$, and using the bound $\Vert D^{2}U^{h}(q)\Vert _{\text{op}}\leq
C|q|^{r_{2}-1}$ established in (\romannumeral3), 
\begin{equation*}
\Vert D^{2}U(q)-D^{2}U^{h}(q)\Vert _{\text{op}}\leq \Vert D^{2}U(q)\Vert _{%
\text{op}}+\Vert D^{2}U^{h}(q)\Vert _{\text{op}}\leq C|q|^{r_{2}-1}.
\end{equation*}%
Applying \eqref{eq:F} with $m=r_{2}-1$ (so $m+2k(r_{2}-1)=(2k+1)(r_{2}-1)$)
yields (\romannumeral4).

\noindent \textbf{(\romannumeral5)} The bound on third derivatives is proved in the same way as
\noindent (\romannumeral3) and hence skipped here. 
\end{proof}  %
	
	\begin{proof} [Proof of Proposition~\ref{prop_ergodicfhS2}] We abbreviate $R=R(h)$ and $%
\lambda =\lambda (h)$.
    From (\ref{eq:fhS2}), we have
    \begin{align} \label{eq:prp221}
		  f^h(q)^\top q 
         & \ge \theta\Big(\frac{|q|}{R}\Big) (c_0 U(q)-c_1) + \Big[1-\theta\Big(\frac{|q|}{R}\Big) \Big]2\lambda  |q|^2 +  \frac{(U(q)- \lambda  |q|^2 ) |q|}{R} \theta'\Big(\frac{|q|}{R}\Big) \nn\\
         & \ge \min\{c_0 , 2\} U^h(q) - c_1 +  \frac{(U(q)- \lambda |q|^2 ) |q|}{R} \theta'\Big(\frac{|q|}{R}\Big), 
	\end{align}
where we used Assumption~\ref{assm_ErgodicU} and (\ref{eq:molpot}).

The third term in (\ref{eq:prp221}) is nonzero only for $|q|\in(R,2R)$. We have 
\begin{align*}
			U(q)- \lambda  |q|^2%
            &\leq L_2 \big( |q|^{r_1}  - (2R_0)^{r_2-1}h^{-1/2}|q|^2\big) +C \leq C,
		\end{align*}
   where we used Assumption~\ref{assm_CandPoly} and the relation $|q|^{r_1}= |q|^{r_2-1} |q|^2 \leq |2R|^{r_2-1} |q|^2
   = (2R_0)^{r_2-1} h^{-1/2} |q|^2$ (see (\ref{equ_scaling})).
   Combining it with $\theta' \leq 0$ and $|q|/R<2$ on the support of $\theta'$, we obtain 
    \begin{align*}
			\frac{(U(q)- (2R_0)^{r_2-1} L_2 h^{-1/2} |q|^2 ) |q|}{R}  \theta' \Big(\frac{|q|}{R}\Big) \geq -C \Vert \theta ^{\prime  }\Vert _{\infty }.
	\end{align*}

		Thus 
		$
			f^h(q)^\top q   \geq \min\{c_0 , 2\} U^h(q) -c_3
		$
        with $c_3=c_1+C \Vert \theta ^{\prime  }\Vert _{\infty }$.
	\end{proof}

 \subsection{Proof of Proposition~\ref{prop_MomentboundedOfTamedequ}}\label{subsec_proofMomentequ}
 
	We first show that ${V}^h(p,q)$  controls both  $|p|^2$ and $|q|^2$. Specifically, there exist  constants $C_1, \, C_2>0$ such that $V^h(p,q)\geq C_1(|p|^2+|q|^2)-C_2$. 

   	For $|q| \leq 2R(h)$, we obtain from (\ref{equ_scaling}): 
	$$\lambda(h)|q|^2= 2^{r_{2}-1}L_{2}R^{r_{2}-1}|q|^{2} \geq L_2|q|^{r_2-1}|q|^2=L_2|q|^{r_2+1}\geq L_1|q|^{r_1},$$ 
	which together with Assumption~\ref{assm_CandPoly} gives
	\begin{align} \label{equ_Uhlowerboundleq1}
		U^h(q) &\geq \theta\Big(\frac{|q|}{R(h)}\Big) \min\big\{ U(q), \lambda(h)|q|^2	\big\} +\Big[ 1- \theta\Big(\frac{|q|}{R(h)}\Big)\Big] \min\big\{ U(q), \lambda(h)|q|^2	\big\}\nn\\
		& \geq \min\big\{ U(q), \lambda(h)|q|^2	\big\} \geq L_1 |q|^{r_1} -C
	\end{align}
	with   constant $C>0$ independent of $h$. 
	Hence
	\begin{align} \label{equ_Vhlower1}
		{V}^h(p,q) %
		&\geq \big(\frac 12 -\frac 14 \big)|p|^2 + ( L_1  -\frac 12 L_1) |q|^{r_1} -C(\upsilon, \alpha_1, r_1, L_1)\nn\\
        &\geq \frac 14 |p|^2+\frac 12 L_1 |q|^2-C,
	\end{align}
	where we used  the $\varepsilon$-Young inequality to estimate (recall that $r_1>2$):
	\begin{align}\label{equ_estimateYoung1}
		\upsilon \alpha_1 |p| |q| \leq \frac 14 |p|^2 + \upsilon^2 \alpha_1^2 |q|^2 \leq \frac 14 |p|^2+ \frac 12 L_1 |q|^{r_1} +C(\upsilon, \alpha_1, r_1, L_1). 
	\end{align}
	
	For  $|q|\geq 2R(h)$, we get (see (\ref{eq:molpot}) and (\ref{equ_scaling})):
	\begin{align} \label{equ_Uhlowerboundgeq1}
		U^h(q) =\lambda(h)|q|^2 =
		\tilde{C}_{U^h} h^{-1/2} |q|^{2},
	\end{align}
	where $\tilde{C}_{U^h}:= (2R_0)^{r_2-1}L_2$.
	Apply the $\varepsilon$-Young inequality to obtain
	\begin{align} \label{equ_estimateYoung2}
		|p| |q| \leq \frac{\delta_1}{2} |p|^2 + \frac{1}{2\delta_1} |q|^2 
	\end{align}
	for $\delta_1>0$.
	By taking $\delta_1=h^{\gamma_1}$ with $\gamma_1 \in(0,1/2)$,  we have 
	\begin{align*} 
		{V}^h(p,q)
		&\geq \Big(\frac 12 -\frac{\upsilon\alpha_1}{2} h^{ {\gamma}_1}\Big)|p|^2 + \Big( \tilde{C}_{U^h}h^{-1/2}  - \frac{\upsilon\alpha_1}{2}h^{ -{\gamma}_1} \Big) |q|^{2} -C.
	\end{align*}
	Therefore, there exist constants $C_{2,1},C_{2,2},C_{2,3}>0$, independent of $h$, such that for all sufficiently small $h>0$,
	\begin{align} \label{equ_Vhlower2}
		{V}^h(p,q)
		\geq  C_{2,1}|p|^2+C_{2,2}|q|^2-C_{2,3}
	\end{align}
	holds.
    
	Using \eqref{equ_Vhlower1} and \eqref{equ_Vhlower2},
	we prove the long-time moment boundedness of $(P^h(t), Q^h(t))$ by showing the corresponding moment boundedness for the Lyapunov function $V^h(P^h(t), Q^h(t))$. %

 \begin{proof}[Proof of Proposition~\ref{prop_MomentboundedOfTamedequ}]
By \cite{Mattingly2002SPA} thanks to Propositions~\ref{prop:fh-fS2} and \ref{prop_ergodicfhS2} (i.e., that $U^h$ grows quadratically at infinity for every fixed $h>0$), 
all moments of $Q^h(t)$ and $P^h(t)$ are finite but bounded by a constant which may depend on $h$ since the potential $U^h$ depends on $h$.
The goal of the below proof is to show that the constant does not depend on $h$.
            
		By It\^o's formula, we get
		\begin{align*}
			&\quad d V^h(P^h(t), Q^h(t)) \\
			& = (\partial_p V^h(P^h(t) ,Q^h(t) ))^\top dP^h(t) +(\partial_q V^h(P^h(t) ,Q^h(t) ))^\top dQ^h(t) \\
			&\quad +\frac{\sigma^2}{2} \tr(\partial^2_{pp} V^h(P^h(t) ,Q^h(t) )) dt \\
			&= \Big(-\upsilon(1-\alpha_1)|P^h(t)|^2 -\upsilon^2 \alpha_1 P^h(t)^\top Q^h(t) - \upsilon \alpha_1 (Q^h(t))^\top f^h(Q^h(t))\Big) dt \\
			&\quad +d\mathcal M+\frac{\sigma^2 d}{2} dt \\
			&=\Big[-\upsilon\kappa V^h(P^h(t), Q^h(t)) -\upsilon \Big(\frac{2-2\alpha_1-\kappa}{2}\Big)|P^h(t)|^2 + \upsilon\kappa U^h(Q^h(t)) \\
			&\quad - \upsilon \alpha_1(Q^h(t))^\top f^h(Q^h(t)) +\upsilon^2 \alpha_1 (\kappa-1) P^h(t)^\top Q^h(t)\Big]dt +d\mathcal M +\frac{\sigma^2 d}{2} dt, 
		\end{align*}
		where $\kappa$ is an arbitrary constant and 
		$d\mathcal M =( P^h(t)+\upsilon\alpha_1Q^h(t))^\top\sigma dW(t)$ is the martingale term. Note that  $\mathcal M$ is a martingale since we have boundedness of 
        the moments as stated at the start of the proof. 
		
		Using Proposition~\ref{prop_ergodicfhS2}, we arrive at 
		\begin{align}\label{equ_dV}
			&\quad d V^h(P^h(t), Q^h(t))\nn \\
			&\leq \Big[-\upsilon\kappa V^h(P^h(t), Q^h(t)) -\upsilon \Big(\frac{2-2\alpha_1-\kappa}{2}\Big)|P^h(t)|^2 - ( \upsilon \alpha_1 c_2-\upsilon\kappa)U^h(Q^h(t)) \nn\\
			&\quad  +\upsilon^2 \alpha_1 (\kappa-1) P^h(t)^\top Q^h(t) \Big]dt+d\mathcal M +Cdt,
		\end{align}
        where the constant $C=c_3 \upsilon \alpha_1 + \frac{\sigma^2 d}{2}>0$ is independent of $h$.
		
        Choose $\kappa>0$ such that 
		$
		\kappa < \min (\alpha_1 c_2, (2-2\alpha_1)) , 
		$
		and define $C_{P^h}:= \upsilon(\frac{2-2\alpha_1-\kappa}{2}) \allowbreak>0$ and $C_{Q^h}:=  \upsilon (\alpha_1 c_2- \kappa ) >0$, both independent of $h$. 
         With this choice, \eqref{equ_dV} reads
\begin{equation}\label{equ_dVremainder}
dV^h(P^h(t),Q^h(t))
\le\Big[-\upsilon\kappa V^h(P^h(t),Q^h(t))+\mathcal{R}^h(P^h(t),Q^h(t))\Big]dt
+d\mathcal{M}+C\,dt,
\end{equation}
where 
$
\mathcal{R}^h(p,q):=-C_{P^h}|p|^2-C_{Q^h}U^h(q)+\upsilon^2\alpha_1(\kappa-1)\,p^\top q .
$
We show that there exist $h_0\in(0,1)$ and a constant $C>0$ independent of $h$
such that
\begin{equation}\label{equ_Rhbound}
\mathcal{R}^h(p,q)\le C
\qquad\text{for all }(p,q)\in\RR^{2d}\text{ and all }h\in(0,h_0),
\end{equation}
by considering separately the two regions determined by $|q|$.

{\textbf{Region \romannumeral 1: $|q|<2R(h)$.}}
Similar to \eqref{equ_estimateYoung1}, the $\varepsilon$-Young inequality gives, with a
constant $C>0$ depending on $C_{P^h},C_{Q^h},L_1,r_1,\upsilon^2\alpha_1(\kappa-1)$
but not on $h$,
$
\upsilon^2\alpha_1(\kappa-1)\,p^\top q
\le\tfrac12 C_{P^h}|p|^2+\tfrac12 C_{Q^h}L_1|q|^{r_1}+C .
$
Combining this with $U^h(q)\ge L_1|q|^{r_1}-C$ from \eqref{equ_Uhlowerboundleq1} yields
$
\mathcal{R}^h(p,q)
\le-\tfrac12 C_{P^h}|p|^2-\tfrac12 C_{Q^h}L_1|q|^{r_1}+C\le C .
$

{\textbf{Region \romannumeral 2: $|q|\ge 2R(h)$.}}
Here $U^h(q)=\tilde C_{U^h}h^{-1/2}|q|^2$ by \eqref{equ_Uhlowerboundgeq1}. Denote
$C_{(P^h,Q^h)}:=|\upsilon^2\alpha_1(\kappa-1)|$ and use \eqref{equ_estimateYoung2} with
$\delta_1=h^{\alpha_2}$, $\alpha_2\in(0,1/2)$, to obtain
$
\upsilon^2\alpha_1(\kappa-1)\,p^\top q
\le\tfrac12 C_{(P^h,Q^h)}h^{\alpha_2}|p|^2
+\tfrac12 C_{(P^h,Q^h)}h^{-\alpha_2}|q|^2,
$
so that for all $h\in(0,h_0)$,
$$
\mathcal{R}^h(p,q)\le
-\big(C_{P^h}-\tfrac12 C_{(P^h,Q^h)}h^{\alpha_2}\big)|p|^2
-\big(C_{Q^h}\tilde C_{U^h}h^{-1/2}-\tfrac12 C_{(P^h,Q^h)}h^{-\alpha_2}\big)|q|^2 \le 0.
$$

Regions {\romannumeral1} - {\romannumeral2} exhaust $(p,q) \in \RR^{2d}$, which proves \eqref{equ_Rhbound}.
Substituting \eqref{equ_Rhbound} into \eqref{equ_dVremainder}, the bound holds along the
trajectory for every $t\ge0$:
\begin{align} \label{eq:Vhextra}
dV^h(P^h(t),Q^h(t))
\le-\upsilon\kappa\,V^h(P^h(t),Q^h(t))\,dt+d\mathcal{M}+C\,dt .   
\end{align}

Taking expectations and using that $\mathcal{M}$ is a martingale, we arrive at
\begin{equation}\label{equ_dVequ}
d\,\E\big[V^h(P^h(t),Q^h(t))\big]
\le-\upsilon\kappa\,\E\big[V^h(P^h(t),Q^h(t))\big]\,dt+C\,dt,
\end{equation}
where $C$ is independent of $h$. Gr\"onwall's inequality then gives
$$\sup_{t\ge0}\E\big[V^h(P^h(t),Q^h(t))\big]\le C .$$
Here $C$ is independent of $h$ because the constant in \eqref{equ_dVequ} is $h$-independent as
shown above, and the initial value is bounded uniformly in $h$ since $R(h)\to\infty$ as
$h\to0$ implies $|q_0|\le R(h)$ for all sufficiently small $h>0$, whence $U^h(q_0)=U(q_0)$
and thus $V^h(p_0,q_0)$ is independent of $h$.

        We use induction to prove the result for higher moments.
		Following \eqref{equ_Vhlower1} and \eqref{equ_Vhlower2}, there exists $C_{V^h}>0$ independent of $h$ such that $\tilde{V}^h(p,q):= V^h(p,q) +C_{V^h} >0 .$  
		We first assume that for $m \ge 1$,
		$
			\sup_{t\geq 0} \E\big[\big|\tilde V^h(P^h(t), Q^h(t)) \big|^m\big] \leq C.
		$
Then one can obtain
\begin{align*}
&d  \E\Big[ \big|\tilde{V}^h(P^h(t), Q^h(t))  \big|^{m+1}\Big] 
=(m+1)\mathbb E \left[\big(\tilde V^h(P^h(t), Q^h(t))\big)^m d V^h \right] \\
&+\frac{m(m+1)}{2}\sigma^2 \mathbb E \left[ \big(\tilde V^h(P^h(t), Q^h(t))\big)^{m-1}
\big|P^h(t)+\upsilon\alpha_1Q^h(t)\big|^2 \right]dt \\
&\leq -\upsilon\kappa(m+1)  \E\Big[ \big|\tilde{V}^h(P^h(t), Q^h(t)) \big|^{m+1}\Big]  dt +C dt,
\end{align*}
where  we used \eqref{eq:Vhextra} and we also used  \eqref{equ_Vhlower1}  and \eqref{equ_Vhlower2} 
		to obtain %
		\begin{align*}
			\big| P^h(t)+\upsilon \alpha_1 Q^h(t)\big|^2 \leq2  |P^h(t)|^2 +2\upsilon^2 \alpha_1^2 |Q^h(t)|^2 \leq C V^h(P^h(t), Q^h(t)) +C
		\end{align*}
with $C>0$.
Hence, by the induction argument, we get for any $m\geq 1$:
		\begin{align*}
			\sup_{t\geq 0} \E\big[\big|\tilde V^h(P^h(t), Q^h(t)) \big|^m\big] \leq C .
		\end{align*}
        Since $|V^h(p,q)|\le \tilde V^h(p,q) + C_{V^h}$, the same bound holds for $V^h$.
Because $|p|^2+|q|^2\le C(\tilde V^h(p,q)+1)$, we get $\E[|Q^h(t)|^m]\le C \E[(\tilde V^h(P^h(t),Q^h(t))+1)^{m/2}]\le C$ and likewise for $P^h (t)$.
\end{proof}
\subsection{Proof of Theorem~\ref{thm_weakconvergencetamed}}\label{subsec_proofweakconvergenceequ}
		(\romannumeral 1) 
		It is obvious that the conformal symplectic structure  is preserved by \eqref{eq:tamed}, i.e.,
		\begin{align*}
			\rd P^h(t) \wedge \rd Q^h(t) = e^{-\upsilon t} \rd P^h(0) \wedge \rd Q^h(0).
		\end{align*}
        
        \noindent
		(\romannumeral 2) 
		Building upon  Proposition~\ref{prop_MomentboundedOfTamedequ} and the minorization condition,  %
		$\{(P^h(t), Q^h(t))\}_{t\geq 0}$ admits a unique invariant measure $\pi_h$ \cite{Mattingly2002SPA}. 
This measure is clearly Gibbsian with density $\rho_h(p,q) \propto \exp\big(-\Lambda\big(\frac12|p|^2+U^h(q)\big)\big)$, $\Lambda:=2\upsilon/\sigma^2$, as can be verified in the usual manner via the corresponding stationary Fokker--Planck equation.

		We now compare $\pi_h$ with $\pi$. Both measures have the $p$-marginal $\mathcal N(0,\Lambda^{-1}I_d)$ and they differ only in their $q$-parts. Let $\phi$ be the Gaussian density of $\mathcal N(0,\Lambda^{-1}I_d)$ and $\bar g(q):=\int_{\RR^d}g(p,q)\phi(p)\,\rd p$, and we have $|\bar g(q)|\leq  C(1+|q|^m)$ for some $C,m>0$. Introduce $w:=e^{-\Lambda U}$, $w_h:=e^{-\Lambda U^h}$, $z:=\int_{\RR^d}w\,\rd q$, and $z_h:=\int_{\RR^d}w_h\,\rd q$. Then we have $\pi(g)=z^{-1}\int_{\RR^d}\bar g\,w\,\rd q$ and $\pi_h(g)=z_h^{-1}\int_{\RR^d}\bar g\,w_h\,\rd q$.  Since $\int_{\RR^d}(\bar g (q) -\pi(g))\,w\,\rd q=0$, we get
		\begin{align*}
			\pi_h(g)-\pi(g)&=\frac{1}{z_h}\int_{\RR^d}\big(\bar g(q)-\pi(g)\big)w_h(q)\,\rd q \\
			&=\frac{1}{z_h}\int_{|q|>R(h)}\big(\bar g(q)-\pi(g)\big)\big(w_h(q)-w(q)\big)\,\rd q,
		\end{align*}
		where the domain of integration reduces to $\{|q|>R(h)\}$ because $U^h=U$ (and hence $w_h=w$) on $\{|q|\leq R(h)\}$. Since $R(h)\geq R_0$, we also have $z_h\geq\int_{|q|\leq R_0}e^{-\Lambda U(q)}\,\rd q=:z_0>0$ for every $h$.

		It remains to bound $w$ and $w_h$ on $\{|q|>R(h)\}$. By Assumption~\ref{assm_CandPoly}, $w(q)\leq e^{\Lambda C}e^{-\Lambda L_1|q|^{r_1}}$ for all $q$. For $R<|q|<2R$, \eqref{equ_Uhlowerboundleq1} gives $U^h(q)\geq L_1|q|^{r_1}-C$ and $w_h(q)\leq e^{\Lambda C}e^{-\Lambda L_1|q|^{r_1}}$. For $|q|\geq2R$, \eqref{equ_Uhlowerboundgeq1} gives $w_h(q)=e^{-\Lambda\lambda(h)|q|^2}\leq e^{-2\Lambda\lambda(h)R^2}e^{-\Lambda\lambda(h)|q|^2/2}$, where $2\lambda(h)R^2=2^{r_2}L_2R^{r_1}\geq L_1R^{r_1}$ and $\lambda(h)\geq2^{r_2-1}L_2$. Therefore, for any $c\in(0,\Lambda L_1)$:
		\begin{align*}
			&\quad |\pi_h(g)-\pi(g)|\\
			&\leq\frac{C}{z_0}\Big(\int_{|q|>R}(1+|q|^m)e^{-\Lambda L_1|q|^{r_1}}\,\rd q
			+e^{-\Lambda L_1R^{r_1}}\int_{\RR^d}(1+|q|^m)e^{-2^{r_2-2}\Lambda L_2|q|^2}\,\rd q\Big)\\
			&\leq Ce^{-cR^{r_1}},
		\end{align*}
		where the first integral accounts for $w$ on $\{|q|>R\}$ and for $w_h$ on $\{R<|q|<2R\}$ and the second for $w_h$ on $\{|q|\geq2R\}$. The first integral is small because of its domain: on $\{|q|>R\}$,
\[
  e^{-\Lambda L_1|q|^{r_1}}\leq e^{-cR^{r_1}}e^{-(\Lambda L_1-c)|q|^{r_1}},
\]
and $c<\Lambda L_1$ keeps the remaining integral finite. 
The second integral is for simplicity extended to all of $\RR^d$ since its smallness comes from the factor $e^{-\Lambda L_1R^{r_1}}$ in front of it.
The constant $C$ depends on $c$, $K$, $m$, $d$, $\Lambda$, $L_1$, $L_2$ and the constant $C$ of Assumption~\ref{assm_CandPoly}, but not on $h$. Since $R^{r_1}(h)=R_0^{r_1}h^{-r_1/(2(r_2-1))}$, we have proved the bound in (\romannumeral 2). This bound implies $|\pi(g)-\pi_h(g)|\leq C(k)h^k$ for every $k>0$ because $x^s e^{-bx}$ is bounded on $(0,\infty)$ for any $b,s>0$, applied here with $x=h^{-r_1/(2(r_2-1))}$, $b=cR_0^{r_1}$ and $s=2k(r_2-1)/r_1$. 
         
         Finally, the first tail estimate together with $z\geq z_0$ also gives $\pi(\{|q|>R(h)\})\leq Ce^{-cR^{r_1}}$, i.e., that the Gibbs measure mass outside the ball of radius $R(h)$ decays super-exponentially in $R_0$ as we discussed in Section~\ref{sec_tamingAndGeneralmethod}.
        
        \noindent(\romannumeral 3) Now we prove the weak convergence result.
		Note that $\E[u(0, P^h(T), Q^h(T))]= \E[ g(P^h(T), Q^h(T) )]$ and $\E[u(T, P^h(0), Q^h(0))]=\E[g( P(T), Q(T))]$.
		Integrating both sides over the interval $[0,T]$ and taking expectations, we arrive at
		\begin{align*}
			\big|\E[g(P^h(T), Q^h(T) )]-\E[ g( P(T), Q(T))]  \big|
			&=\Big| \E\Big[ \int^T_0 du(T-t,P^h(t),Q^h(t)) \Big] \Big|. 
		\end{align*}
		Using It\^o's formula, for any fixed $T>0$ and $t\in[0,T]$, we obtain
		\begin{align*}
			&\quad du(T-t,P^h(t),Q^h(t))  \\
			&= -\partial_t u(T-t,P^h(t),Q^h(t))  dt +(\partial_p u(T-t,P^h(t),Q^h(t)) )^{\top}dP^h(t) \\
			&\quad + (\partial_q u(T-t,P^h(t),Q^h(t)) )^{\top}dQ^h(t) +\frac{\sigma^2}{2} \tr( \partial_{pp}u (T-t,P^h(t),Q^h(t)) ) dt\\
			&=\big(\mathcal L^h- \mathcal L\big) u(T-t,P^h(t),Q^h(t))  dt + (\partial_p u(T-t,P^h(t),Q^h(t)) )^{\top}\sigma dW(t) \\
			&=\big( f(Q^h(t))-f^h(Q^h(t))\big)^\top \partial_p u (T-t,P^h(t),Q^h(t)) dt\\
            &\quad+ (\partial_p u(T-t,P^h(t),Q^h(t)) )^{\top}\sigma dW(t).
		\end{align*}
		It follows from Proposition~\ref{prop:fh-fS2}(i) and $\E[\int^T_0 (\partial_p u(T-t,P^h(t),Q^h(t)) )^{\top}\sigma dW(t)] =0$ that, for any $k>0$,
		\begin{align*}
			&\quad \big|\E[g(P^h(T), Q^h(T) )]-\E[ g( P(T), Q(T))]  \big|\\
			&=\Big|\int^T_0 \E\big[\big( f(Q^h(t))-f^h(Q^h(t))\big)^\top \partial_p u (T-t,P^h(t),Q^h(t)) \big]dt \Big| \\
			&\leq Ch^k \int^T_0 \E\Big[( 1+|Q^h(t)|^{(2k+1)r_2-2k} ) \big|\partial_p u (T-t,P^h(t),Q^h(t))\big| \Big]dt  \\
			&\leq Ch^k \int^T_0\Big( \big(\E\big[1+ |Q^h(t)|^{2\tilde{k}} \big] \big)^{\frac 12}\big( \E\big[ \big|\partial_p u (T-t,P^h(t),Q^h(t))\big|^2 \big]\big)^{\frac 12}\Big)dt \\
			&\leq C(k)h^k \int^T_0  \big( \E\big[ \big|\partial_p u (T-t,P^h(t),Q^h(t))\big|^2 \big]\big)^{\frac 12}dt
			\leq C(k)h^k,
		\end{align*}
		where $\tilde{k}:=(2k+1)r_2-2k$ and \eqref{eq:Dku} and Proposition~\ref{prop_MomentboundedOfTamedequ} are used.
        The constant $C(k)$ does not depend on $T$. Indeed, \eqref{eq:Dku} and Proposition~\ref{prop_MomentboundedOfTamedequ} imply that $\big(\E\big[|\partial_p u(T-t,P^h(t),Q^h(t))|^2\big]\big)^{1/2}\leq Ce^{-\gamma(T-t)}$ with $C$ independent of $t$, $T$ and $h$. Moreover,  $\int_0^Te^{-\gamma(T-t)}\,\rd t\leq1/\gamma$. Since $T>0$ is arbitrary, this proves (\romannumeral 3).		

\section{Proofs for the numerical methods}\label{suppsec_proofofnumerical}
\subsection{Proof of Proposition~\ref{prop_momentboundedgeneral}}\label{subsec_proofmomentnumerical}
		The definition of ${V}^h$ yields
		\begin{align} \label{equ_generalmoment1}
		{V}^h (P^{h}_{n+1},Q^{h}_{n+1}) 
			&= \frac 12 |P^{h}_{n+1}|^2 + U^h(Q^{h}_n)+( f^h(Q^{h}_n) )^\top ( Q^{h}_{n+1}-Q^{h}_n)  \nn \\
			&\quad + \frac 12 (Q^{h}_{n+1}-Q^{h}_n)^\top D^2 U^h(\xi_1) (Q^{h}_{n+1}-Q^{h}_n) + \upsilon\alpha_1 ( P^{h}_{n+1})^\top Q^{h}_{n+1},
		\end{align}
		where we use the Taylor expansion of $U^h(Q^{h}_{n+1})$ and $\xi_1:=\varepsilon Q^h_{n}+(1-\varepsilon)Q^h_{n+1}$ for some $\varepsilon\in(0,1)$. 
		We next expand each term in \eqref{equ_generalmoment1} by substituting the one-step update formulas \eqref{eq:generalP}--\eqref{eq:generalQ}.
		
		Using $h^{1/2}\|D^2 U^h (\cdot)\|_{\op} \leq C$, for the first term on the right-hand side of \eqref{equ_generalmoment1}, we derive
		\begin{align}\label{equ_momentgeneralP2}
		&\quad	\frac 12 \E[	|P^{h}_{n+1}|^2 |\mathcal{F}_{t_n}] \nn\\
			&= \frac 12 \E\Big[ \big|	\big(I_d +[-\upsilon I_d+a_{0,p,1} h \bar H_{p,n}] h +[a_1 I_d+ a_{0,p,2} h \bar H_{p,n}]h^2 \big)P^h_n \nn \\
			&\quad + (- h I_d+[a_2 I_d+a_{0,q} h \bar H_{p,n}]h^2 ) f^h(Q^h_n) + \sigma \Delta W_n + a_3 h^2 %
			+ \mathcal{Z}_{p, n,1} \nn \\
			&\quad +h \bar H_{p,n} \mathcal{Z}_{p, n,2}+\mathcal{R}_{p,n} \big|^2 \Big| \mathcal{F}_{t_n} \Big]   \nn \\
			&\leq \frac 12  \big(1-2\upsilon h +C_{p^2,h^{3/2},1}h^{3/2}
			\big) |P^h_n |^2+\frac 12 \big(-2h+2(\upsilon +a_2)h^2\big) ( P^h_n)^\top f^h(Q^h_n)  \nn \\
			&\quad  +\frac 12  h^2 |f^h(Q^h_n)|^2+ Ch^3|f^h(Q^h_n)|^2+ C h,
		\end{align}
		where $	C_{p^2,h^{3/2},1}>0$ is a constant. 
		Here, for the last inequality in \eqref{equ_momentgeneralP2}, we use the estimates %
		$$h^3|\bar H_{p,n} P^h_n| |f^h(Q^h_n)|\leq C h^{5/2} |P^h_n| |f^h(Q^h_n)|\leq  Ch^2|P^h_n|^2+Ch^3|f^h(Q^h_n)|^2,$$
		\begin{align*}
			\E\big[ |(P^h_n)^\top  \bar H_{p,n}  h^2 \Delta W_n | \big|\mathcal{F}_{t_n}\big] \leq Ch^{3/2}|P^h_n| \E[ |\Delta W_n|  |\mathcal{F}_{t_n}]  \leq Ch^3 |P^h_n|^2 +Ch,
		\end{align*}
		\begin{align*}
			\E\big[ | (P^h_n)^\top \bar H_{p,n}  h  \mathcal{Z}_{p,n,2}| \big|\mathcal{F}_{t_n}\big] \leq C h^3 |P^h_n|^2 +Ch,
		\end{align*}
        $\E[ |\sigma \Delta W_n|^2]\leq Ch$, \,
		$h^2|a_3||P^h_n| \leq Ch+Ch^3 |P^h_n|^2$, \,$|(P^h_n)^\top \mathcal{R}_{p,n}|\leq Ch^3 (|P^h_n|^2+ |f^h(Q^h_n)|^2)+Ch$, and $h|(f^h(Q^h_n))^\top \mathcal{R}_{p,n}|\leq Ch^4(|P^h_n|^2+ |f^h(Q^h_n)|^2)+Ch$.  %
		
		Similarly,  for the third term on the right-hand side of \eqref{equ_generalmoment1}, it holds that
		\begin{align}\label{equ_momentgeneralUq}
			&\quad	\E[  (f^h(Q^h_n))^\top (Q^h_{n+1}- Q^h_n) | \mathcal{F}_{t_n} ]  \nn\\
			&\leq  (h+b_1  h^2) (P^h_n)^\top  f^h(Q^h_n)  + b_2 h^2 |f^h(Q^h_n)|^2 + Ch^2|P^h_n|^2+ Ch^3|f^h(Q^h_n)|^2 +Ch.
		\end{align}
		
		The fourth term on the right-hand side of \eqref{equ_generalmoment1} has the following upper bound
		\begin{align}\label{equ_momentgeneralQ2}
			\frac 12\E[ (Q^{h}_{n+1}-Q^{h}_n)^\top D^2 U^h(\xi_1) (Q^{h}_{n+1}-Q^{h}_n)  | \mathcal{F}_{t_n}  ] 
			&\leq C h^{3/2} | P^h_n|^2
			+ Ch^3|f^h(Q^h_n)|^2 +Ch.
		\end{align}
		
		The last term on the right-hand side of \eqref{equ_generalmoment1} can be estimated as 
		\begin{align}\label{equ_momentgeneralPQ}
			&\quad \upsilon \alpha_1\E[	( P^{h}_{n+1})^\top  Q^{h}_{n+1}  | \mathcal{F}_{t_n}  ]\nn \\
			&\leq \upsilon \alpha_1h|P^h_n|^2+\upsilon \alpha_1(P^h_n)^\top  Q^{h}_{n}+C_{pq,h} h|P^h_n| |Q^{h}_{n}| +\upsilon \alpha_1 \big(- h +a_2 h^2 \big) (f^h(Q^h_n))^\top Q^h_n  \nn\\
			&\quad +\upsilon\alpha_1(b_2-1) h^2 (P^h_n)^\top  f^h(Q^{h}_{n})+C_{q^2,h^2} h^2|Q^h_n|^2 +Ch^2|P^h_n|^2 + Ch^3|f^h(Q^h_n)|^2 +Ch 
		\end{align}
		with positive constants $C_{pq,h}$ and $  C_{q^2, h^2}$,
		where the estimate  $C h^{5/2}|f^h(Q^h_n)||Q^h_n| \leq Ch^2|Q_n^h|^2 +Ch^3|f^h(Q^h_n)|^2$ is used.
		
		Combining  \eqref{equ_momentgeneralP2}--\eqref{equ_momentgeneralPQ}, 
		we obtain
		\begin{align*}
			&\quad 	\E\big[	{V}^h (P^h_{n+1},Q^h_{n+1}) \big| \mathcal{F}_{t_n} \big] \\
			&\leq \Big( \frac 12 +(\upsilon\alpha_1-\upsilon )h+C_{p^2,h^{3/2}}h^{3/2}\Big) |P^h_n|^2 +C_{pq,h^2}h^2 
			|P^h_n| |f^h(Q^h_n)|\\
			&\quad  +\upsilon\alpha_1(P^h_n)^\top Q^h_n  + C_{pq,h}
			h 	|P^h_n| |Q^h_n|+\big( \frac 12 +b_2 \big)h^2 |f^h(Q^h_n)|^2 +C_{q^2,h^2} h^2|Q^h_n|^2\\
			&\quad + \upsilon\alpha_1 \big(-h+a_2h^2\big)  (f^h(Q^h_n))^\top Q^h_n  + U^h(Q^h_n)+ C_{f^2,h^3}h^3|f^h(Q^h_n)|^2+ Ch \\
			&= (1-\kappa h) {V}^h (P^h_{n},Q^h_{n}) \\
			&\quad + \Big( (\frac 12 \kappa +\upsilon\alpha_1-\upsilon  )h +C_{p^2,h^{3/2}}h^{3/2}\Big) |P^h_n|^2 +C_{pq,h^2}h^2
			|P^h_n| |f^h(Q^h_n)|  \\
			&\quad    +\upsilon\alpha_1\kappa h (P^h_n)^\top Q^h_n
			+ C_{pq,h}	h|P^h_n| |Q^h_n|+\big( \frac 12 +b_2\big)h^2 |f^h(Q^h_n)|^2 +C_{q^2,h^2} h^2|Q^h_n|^2\\
			&\quad + \upsilon\alpha_1 \big(-h+a_2h^2\big)  (f^h(Q^h_n))^\top Q^h_n  + \kappa h U^h(Q^h_n)+ C_{f^2,h^3}h^3|f^h(Q^h_n)|^2+ Ch,
		\end{align*}
		where $C_{p^2,h^{3/2}}, C_{pq, h^2}, C_{pq, h}, C_{q^2,h^2},C_{f^2,h^3}, \kappa$ are positive constants.
		By using Proposition~\ref{prop_ergodicfhS2} and denoting
		$A_h:=-\frac 12\kappa -\upsilon\alpha_1+\upsilon -C_{p^2,h^{3/2}}h^{1/2} $, $B_h:=C_{pq,h^2}$, $C_h:=  \frac 12 +b_2 $, $D_h:=\upsilon\alpha_1\kappa + C_{pq,h}$, and $E_h:=-\kappa+\upsilon\alpha_1 c_2 -\upsilon\alpha_1c_2a_2h$, we arrive at
		\begin{align*}
			\E\big[	{V}^h (P^h_{n+1},Q^h_{n+1}) \big| \mathcal{F}_{t_n} \big] \leq  (1-\kappa h) {V}^h (P^h_{n},Q^h_{n})+\tilde{\mathcal{R}}^h(P^h_n,Q^h_n),
		\end{align*}
    where
    \begin{align*}
        \tilde{\mathcal{R}}^h(p,q)&:= - A_h h |p|^2
			+ B_h h^2 |p| |f^h(q)| +C_{h} h^2 |f^h(q)|^2  +D_h h |p| |q| -E_h h U^h(q)\\
            &\quad+C_{q^2,h^2} h^2 |q|^2 +  C_{f^2,h^3}h^3|f^h(q)|^2+Ch.
    \end{align*}
We claim that for sufficiently small $h>0$, there exists a constant $C>0$ such that
\begin{align}\label{equ_numericalboundClaim}
    \frac{\tilde{\mathcal {R}}^h(p,q)}{h}\leq C \quad \text{for all } (p,q)\in\mathbb{R}^{2d}.
\end{align}
        
		By taking $$\kappa\in\Big(0,\frac{2c_2\upsilon}{2+c_2} \Big),$$
		there exists $\alpha_1\in(0,1)$ such that 
		\begin{align*} %
			\frac{\kappa}{c_2 \upsilon}< \alpha_1 <  1 -\frac {1}{2\upsilon}\kappa, 
		\end{align*}
		which guarantees $A_h>0$ and $E_h>0$ for sufficiently small $h>0$.  
		
		It remains to absorb the cross terms $B_h h^2|p||f^h(q)|$, and $D_h h |p| |q|$,  and the square terms $C_hh^2|f^h(q)|^2$, $C_{q^2,h^2} h^2|q|^2$, $C_{f^2,h^3} h^3|f^h(q)|^2$. We decompose the estimates into two regions determined by $|q|$ below.
		
		{\textbf{Region \romannumeral 1: $|q|<2R(h)$.}} In this case, it follows from Proposition~\ref{prop:fh-fS2} (\romannumeral 2) that $|f^h(q)|\leq C(1+|q|^{r_2})$.  Since $(R(h))^{r_2-1} = R_0^{r_2-1} h^{-1/2}$, we have
		$$ |q|^{2r_2}= |q|^{ r_2-1} |q|^{r_2+1} \leq(2R(h))^{r_2-1} |q|^{r_2+1}\leq Ch^{-1/2} |q|^{r_1}$$ 
		for  $|q|<2R(h)$, which leads to $|f^h(q)|^{2} \leq C|q|^{2r_2} +C\leq Ch^{-1/2} |q|^{r_1} +C$ in this case.
		The $\varepsilon$-Young inequality yields
		\begin{align*} %
			B_h h^2|p||f^h(q)| \leq \frac{B_h}{2} h^2 |p|^2 +\frac{B_h}{2} h^2 |f^h(q)|^{2} \leq  \frac{B_h}{2} h^2 |p|^2 + C\frac{B_h}{2} h^{\frac 32} |q|^{r_1} +Ch^2,
		\end{align*}
		\begin{align*}
            C_h h^2 |f^h(q)|^2\leq Ch^{3/2} |q|^{r_1} +Ch^2,
        \end{align*}
		\begin{align*} %
			D_h h |p| |q| \leq \frac 12 A_h h  |p|^2 + C({A_h,D_h}) h |q|^2 \leq  \frac 12 A_h h  |p|^2 + \frac 12 E_h L_1 h|q|^{r_1}+Ch,
		\end{align*}
		$
			C_{f^2,h^3}h^3 |f^h(q)|^2 \leq  Ch^{5/2} |q|^{r_1} +Ch^3,
		$
		and $C_{q^2,h^2} h^2|q|^2 \leq h^2 |q|^{r_1} +Ch^2.$
		Using \eqref{equ_Uhlowerboundleq1},
		we arrive at
		\begin{align*}
		&\quad	\frac{\tilde{\mathcal R}^h(p,q)}{h}\\
        &\leq \Big( -A_h   +\frac 12 A_h +\frac{B_h}{2} h \Big) |p|^2+C \\
			&\quad +\Big(- E_h L_1+\frac 12 E_h L_1 +  C\frac{B_h}{2} h^{1/2}+Ch^{1/2} +h +Ch^{3/2}\Big)  |q|^{r_1} \\
			&\leq \Big(-\frac 12A_h  +\frac{B_h}{2} h \Big) |p|^2  + \Big( -\frac 12 E_h L_1+ C\frac{B_h}{2}h^{1/2}+Ch^{1/2} +h+C h^{3/2} \Big)  |q|^{r_1}+C.
		\end{align*}
		For sufficiently small  $h>0$, it holds that 
		$-\frac 12 A_h+\frac{B_h }{2}h<0$ and $-\frac 12 E_h L_1+ C\frac{B_h}{2}h^{1/2} +Ch^{1/2}+h + Ch^{3/2}<0$. Therefore,
		we obtain
		$
		\frac{\tilde{\mathcal R}^h(p,q)}{h} \leq  C.
		$
		
		{\textbf{Region \romannumeral 2: $|q|\geq 2R(h)$.}}
		In this case, Proposition~\ref{prop:fh-fS2} (\romannumeral 2) implies $|f^h(q)|\leq Ch^{-1/2} |q|$. We apply the $\varepsilon$-Young inequality to derive
		$$	B_h h^2|p||f^h(q)| \leq \frac{B_h}{2} h^2 |p|^2 +\frac{B_h}{2} h^2 |f^h(q)|^{2}\leq \frac{B_h}{2} h^2 |p|^2 +C\frac{B_h}{2}  h|q|^2,
		$$ 
        $
            C_hh^2|f^h(q)|^2 \leq C |C_h| h|q|^2,
        $
		$
			D_h h |p| |q| \leq \frac{1}{2} A_h h |p|^2 +C(A_h,D_h) h |q|^{2},
		$
        and
        $
            C_{f^2,h^3}h^3 \allowbreak |f^h(q)|^2 \leq CC_{f^2,h^3}h^2 |q|^2.
       $
		Using  the lower bound estimate \eqref{equ_Uhlowerboundgeq1}, we obtain
		\begin{align*}
			\frac{\tilde{\mathcal R}^h (p,q)}{h}
			&\leq  \Big(-A_h+\frac 12 A_h +\frac{B_h}{2}h \Big) |p|^2  + \Big( -E_h \tilde{C}_{U^h}+\frac{CB_h}{2} h^{1/2}+ C(A_h,D_h)h^{1/2}
            \\
			&\quad +C|C_h| h^{1/2}+C_{q^2,h^2} h^{3/2 } +CC_{f^2,h^3}h^{3/2} \Big)  h^{-1/2} |q|^{2} + C,
		\end{align*}
		which implies $\frac{\tilde{\mathcal{R}}^h(p,q)}{h}\leq C$.
		
    {Regions \romannumeral 1}--{ \romannumeral 2} exhaust $(p,q)\in \mathbb{R}^{2d}$, which proves \eqref{equ_numericalboundClaim} and leads to the upper bound
    $\sup_{n\in\N} \E[ {V}^h(P^h_n,Q^h_n)] \leq C$.

Higher moments follow from the same expansion. Let $\tilde V_n:=V^h(P^h_n,Q^h_n)+C_0$ and choose $C_0>0$ by \eqref{equ_Vhlower1}--\eqref{equ_Vhlower2} so that $\tilde V_n\geq1$ and $|P^h_n|^2+|Q^h_n|^2+U^h(Q^h_n)\leq C\tilde V_n$. Introduce $\Delta_n:=\tilde V_{n+1}-\tilde V_n$. Carrying out the expansion \eqref{equ_generalmoment1} without conditional expectations, the terms $-h(P^h_n)^\top f^h(Q^h_n)$ and $h(f^h(Q^h_n))^\top P^h_n$ cancel exactly. The remaining terms are estimated with $|f^h(q)|^2\leq Ch^{-1/2}(U^h(q)+1)$, $|(f^h(q))^\top q|\leq C(U^h(q)+1)$ and $h^{1/2}\|D^2U^h(q)\|_{\op}\leq C$, which follow from Proposition~\ref{prop:fh-fS2} (\romannumeral 2)--(\romannumeral 3), \eqref{equ_Uhlowerboundleq1} and \eqref{equ_Uhlowerboundgeq1}. This gives
		\begin{align*}
			|\Delta_n|&\leq C\big(h\tilde V_n+h^{1/2}\tilde V_n^{1/2}\zeta_n+h\zeta_n^2\big),\\
			\zeta_n&:=1+h^{-1/2}|\Delta W_n|+h^{-3/2}\sum_{i=1}^2\big(|\mathcal Z_{p,n,i}|+|\mathcal Z_{q,n,i}|\big),
		\end{align*}
		where $\E[(\zeta_n)^j|\mathcal F_{t_n}]\leq C(j)$ for every $j\geq1$. Fix $k\geq2$. By Taylor's formula, $(\tilde V_n+\Delta_n)^k\leq\tilde V_n^k+k\tilde V_n^{k-1}\Delta_n+C(k)\big(\tilde V_n^{k-2}\Delta_n^2+|\Delta_n|^k\big)$ and $\E[\Delta_n|\mathcal F_{t_n}]\leq-\kappa h\tilde V_n+Ch$ by \eqref{equ_numericalboundClaim}. Together with the bound on $|\Delta_n|$ and Young's inequality, this yields for sufficiently small $h>0$: 
		\begin{align*}
			\E\big[\tilde V_{n+1}^k\big|\mathcal F_{t_n}\big]&\leq(1-k\kappa h)\tilde V_n^k+C\big(h\tilde V_n^{k-1}+h^2\tilde V_n^k+h^{k/2}\tilde V_n^{k/2}\big)+Ch\\
			&\leq\Big(1-\frac{k\kappa h}{2}\Big)\tilde V_n^k+Ch.
		\end{align*}
		Therefore, $\sup_{n\in\N}\E[\tilde V_n^k]\leq\tilde V_0^k+2C/(k\kappa)$ for every $k\geq2$.
        		Using \eqref{equ_Vhlower1} and \eqref{equ_Vhlower2}, we end the proof. 

    \subsection{Proof of Propositions \ref{prop_example1st}--\ref{prop_exampleOU}}\label{subsec_proofconcreteexamples}
     \begin{proof}[Proof of Proposition~\ref{prop_example1st}]
    (\romannumeral 1) The conformal symplecticity follows from the identity $	\rd P^h_{n+1} \wedge \rd Q^h_{n+1} 
		=\rd P^h_{n+1} \wedge \rd Q^h_n = e^{-\upsilon h} \rd  P^h_{n} \wedge  \rd Q^h_n.$
        
    Then 
    we  verify  the conditions in Theorem~\ref{thm_weakconvergeconformalsym} for \eqref{eq:ConformalSym}.
		By rewriting \eqref{eq:ConformalSym} as
		\begin{align} \label{eq:PEuler+R}
			P^h_{n+1}
			&= (1-\upsilon h+\frac 12 \upsilon^2 h^2)P^h_n+(-h+ \upsilon h^2)f^h(Q^h_n)  +\sigma \Delta W_n+\mathcal{Z}_{p,n,1}+ \mathcal R_{p,n,1} 
		\end{align}
		and
		\begin{align}\label{eq:QEuler+R}
			Q^h_{n+1} %
			&=Q^h_n+hP^h_n-h^2 \big(\upsilon P^h_n+ f^h(Q^h_n) \big) +h\sigma e^{-\upsilon h} \Delta W_n +  \mathcal R_{q,n,1} \nn \\
			&=(h-\upsilon h^2) P^h_n +Q^h_n -h^2 f^h(Q^h_n) +\mathcal{Z}_{q,n,1}+ \mathcal R_{q,n,1},
		\end{align}
		and comparing \eqref{eq:PEuler+R}--\eqref{eq:QEuler+R} with the general framework \eqref{eq:generalP}--\eqref{eq:generalQ}, we derive $a_1=\frac 12 \upsilon^2, a_2=\upsilon ,b_1=-\upsilon, b_2=-1$, and $a_3=0, a_{0,p,1}=0, a_{0,p,2}=0, a_{0,q}=0, b_3=0, b_{0,p}=0$. Here, $\mathcal{Z}_{p,n,1} =(e^{-\upsilon h} -1) \sigma \Delta W_n$, $\mathcal{Z}_{q,n,1}=h\sigma e^{-\upsilon h} \Delta W_n $, $\mathcal R_{p,n,1}=\big[e^{-\upsilon h} -(1-\upsilon h+\frac 12 \upsilon^2 h^2) \big]P^h_n-\big[e^{-\upsilon h}- (1-\upsilon h) \big] h f^h(Q^h_n)$, $\mathcal R_{q,n,1} =\big[e^{-\upsilon h} -(1-\upsilon h) \big] hP^h_n-( e^{-\upsilon h}- 1 ) h^2 f^h(Q^h_n)$. 
        Assumption~\ref{assp:nondeg}  holds with $m=2$ and $\xi=(\Delta W_0,\Delta W_1)$: given $(P^h_2,Q^h_2)$ and $x$, \eqref{eq:ConformalSym} determines successively $Q^h_1=Q^h_2-hP^h_2$, $P^h_1=(Q^h_1-Q^h_0)/h$, $\Delta W_0=\sigma^{-1}\big(e^{\upsilon h}P^h_1-P^h_0+hf^h(Q^h_0)\big)$ and $\Delta W_1=\sigma^{-1}\big(e^{\upsilon h}P^h_2-P^h_1+hf^h(Q^h_1)\big)$. Hence $\Gamma_x$ is a bijection with continuously differentiable inverse.
		Then the ergodicity and the weak convergence result follow from Theorem~\ref{thm_weakconvergeconformalsym} (\romannumeral 1).

        (\romannumeral 2) For the numerical scheme based on the Lie--Trotter splitting technique,  we 
        obtain
		$\rd P^{h, \text{LT}}_{n+1} \wedge \rd Q^{h, \text{LT}}_{n+1}
        =  e^{-\upsilon h} \rd  P^{h, \text{LT}}_n  \wedge \rd Q^{h, \text{LT}}_{n},$
		which shows that the conformal symplectic structure is inherited. 
        
        Then we prove the weak convergence result. We first derive
		\begin{align}\label{eq:solutionofsplitting1P}
			P^{h, \text{LT}}_{n+1} %
			&=e^{-\upsilon h} P^{h, \text{LT}}_n - he^{-\upsilon h} f^h({Q}^{h, \text{LT}}_{n})+\sigma \int^{t_{n+1}}_{t_n} e^{-\upsilon(t_{n+1}-s)} dW(s), 
		\end{align}
		and 
		\begin{align}\label{eq:solutionofsplitting1Q}
			 Q^{h, \text{LT}}_{n+1}
			&=Q^{h, \text{LT}}_n+ \frac{1-e^{-\upsilon h}}{\upsilon} P^{h, \text{LT}}_n -h \frac{1-e^{-\upsilon h}}{\upsilon}  f^h({Q}^{h, \text{LT}}_{n})  \nn\\
            &\quad +  \sigma\int^{t_{n+1}}_{t_n} \frac{1-e^{-\upsilon(t_{n+1}-s)}}{\upsilon} dW(s).
		\end{align}
		Furthermore, \eqref{eq:solutionofsplitting1P} can be rewritten as
		\begin{align*}
			P^{h, \text{LT}}_{n+1} 
			&=(1-\upsilon h +\frac 12\upsilon^2 h^2)P^{h, \text{LT}}_n +(- h +\upsilon h^2) f^h({Q}^{h, \text{LT}}_{n})+\sigma \Delta W_n + \mathcal{Z}_{p,n,1} +\mathcal{R}_{p,n,1}, 
		\end{align*}
		where 
		$	\mathcal{Z}_{p,n,1} := \sigma \int^{t_{n+1}}_{t_n} \big(e^{-\upsilon(t_{n+1}-s)} -1 \big)dW(s)$
		satisfying
		$\E[	\mathcal{Z}_{p,n,1}  |\mathcal{F}_{t_n}]=0$ and 
		$\E\big[ |	\mathcal{Z}_{p,n,1}|^{\alpha} |\mathcal{F}_{t_n} \big]  = \mathcal{O} (h^{3\alpha/2})$ for $\alpha\geq 2$, and $\mathcal{R}_{p,n,1} =\big[e^{-\upsilon h}- (1-\upsilon h +\frac 12\upsilon^2 h^2)\big] P^{h, \text{LT}}_n - h\big[e^{-\upsilon h}- (1-\upsilon h )\big] f^h(Q^{h, \text{LT}}_n)$.
		\eqref{eq:solutionofsplitting1Q} can be rewritten as
		\begin{align*}
			Q^{h, \text{LT}}_{n+1}
			&=Q^{h, \text{LT}}_n+ (h -\frac 12 \upsilon h^2)P^{h, \text{LT}}_n -h^2 f^h({Q}^{h, \text{LT}}_{n})  + \mathcal{Z}_{q,n,1} + \mathcal{R}_{q,n,1},
		\end{align*}
		where
		$		\mathcal{Z}_{q,n,1} :=\frac{\sigma}{\upsilon} \int^{t_{n+1}}_{t_n} \big(1-e^{-\upsilon(t_{n+1}-s)}\big) dW(s) = -\frac{1}{\upsilon} 	\mathcal{Z}_{p,n,1}$
		satisfying
		$\E[	\mathcal{Z}_{q,n,1}  |\mathcal{F}_{t_n}]=0$ and $\E [ |	\mathcal{Z}_{q,n,1}|^{\alpha}   |\mathcal F_{t_n}]= \mathcal{O} (h^{3\alpha/2})$ for $\alpha\geq 2$, and  $\mathcal{R}_{q,n,1}= \big[ \frac{1-e^{-\upsilon h}}{\upsilon} - (h -\frac 12 \upsilon h^2)\big]P^{h, \text{LT}}_n \allowbreak -h\big(\frac{1-e^{-\upsilon h}}{\upsilon} - h \big) f^h (Q^{h, \text{LT}}_n)$.
		Comparing with 
		\eqref{eq:generalP}--\eqref{eq:generalQ}, we derive the concrete coefficients 
		$a_1=\frac 12 \upsilon^2, a_2=\upsilon, b_1=-\frac{\upsilon}{2}, b_2=-1$, and $a_3=0, a_{0,p,1}=0, a_{0,p,2}=0, a_{0,q}=0, b_3=0, b_{0,p}=0$. 
        Assumption~\ref{assp:nondeg} holds with $m=1$: by \eqref{eq:solutionofsplitting1P}-\eqref{eq:solutionofsplitting1Q}, $(P^{h,\text{LT}}_1,Q^{h,\text{LT}}_1)$ is a function of $x$ and $\sigma\xi$ with $\xi=(\xi_1,\xi_2):=\big(\int_0^he^{-\upsilon(h-s)}\,d W(s), \allowbreak \int_0^h\frac{1-e^{-\upsilon(h-s)}}{\upsilon}\,d W(s)\big)$, which is a non-degenerate Gaussian vector in $\RR^{2d}$ because the two integrands are linearly independent. Therefore $\Gamma_x$ is an affine bijection.
        The ergodicity and the weak convergence follow from Theorem~\ref{thm_weakconvergeconformalsym} (\romannumeral 1).
    \end{proof}
    
    \begin{proof}[Proof of Proposition~\ref{prop_exampleSplitting}]
        The conformal symplecticity of \eqref{eq:splittingPRK} follows from the fact that 
		$
			\rd P^{h,\text{PRK}}_{n+1} \wedge \rd Q^{h,\text{PRK}}_{n+1} %
			= e^{-\upsilon h} \rd P^{h,\text{PRK}}_n \wedge \rd Q^{h,\text{PRK}}_n.
		$

		The Taylor expansion of \eqref{eq:splittingPRK} leads to
		\begin{align} \label{equ_PhS2}
		P^{h,\text{PRK}}_{n+1} 
			&= \Big(1+\big[-\upsilon -\frac 12 h \bar H_{p,n} \big]h +\big[\frac{\upsilon ^2 }{2} +\frac{\upsilon }{2} h \bar H_{p,n} \big] h^2 \Big) P^{h,\text{PRK}}_n   \nn \\
			&\quad +\big(- h + \frac{\upsilon h^2}{2}\big)  f^h(Q^{h,\text{PRK}}_n) +\sigma \Delta W_n +\mathcal{Z}_{p,n,1} +\mathcal{R}_{p,n,1}
            + h \bar H_{p,n} \mathcal{R}_{p,n,2},
		\end{align}
		where $\bar H_{p,n}$ is defined with $\mathcal{Q}_{p,n}=Q^{h,\text{PRK}}_{n+\frac 12}$, and
		\begin{align} \label{equ_QhS2}
			Q^{h,\text{PRK}}_{n+1} %
			&=\big(h+[-\frac \upsilon 2 -\frac 14 h\bar H_{q,n}] h^2\big) P^{h,\text{PRK}}_n  + Q^{h,\text{PRK}}_n - \frac 12 h^2 f^h(Q^{h,\text{PRK}}_n) + \mathcal{Z}_{q, n,1} \nn\\
			&\quad +\mathcal{R}_{q,n,1} + h^2 \bar H_{q,n} \mathcal{R}_{q,n,2}, 
		\end{align} 
		where $\bar H_{q,n}:=\bar H_{p,n}$, the stochastic convolutions $\mathcal{Z}_{p,n,1} :=-\upsilon \sigma\Delta B_n= -\upsilon \sigma\int^{t_{n+1}}_{t_n} \int^s_{t_n} \allowbreak dW(r)ds$ satisfying $\E[ \mathcal{Z}_{p,n,1} |\mathcal F_{t_n}]=0$ and $\E[ |\mathcal{Z}_{p,n,1}|^{\alpha} |\mathcal F_{t_n}] = \mathcal{O}(h^{3\alpha/2})$ for $\alpha\geq 2$, and $\mathcal{Z}_{q,n,1}:= -\frac{1}{\upsilon} \mathcal{Z}_{p,n,1}+ \frac{\sigma h}{2} (e^{\frac{\upsilon h}{2}}-1)\Delta W_n - \frac{\upsilon \sigma h}{2} e^{\frac{\upsilon h}{2}} \Delta B_n $ satisfying $\E[ \mathcal{Z}_{q,n,1} |\mathcal F_{t_n}]=0$, $\E[ |\mathcal{Z}_{q,n,1}|^{\alpha} |\mathcal F_{t_n}] = \mathcal{O}(h^{3\alpha/2})$ for $\alpha\geq 2$, and $\E[ |\mathcal{Z}_{q,n,1}- \sigma \Delta B_n |^2  |\mathcal{F}_{t_n}]=\mathcal{O}(h^{5})$, which implies $|\E[ (\mathcal{Z}_{q,n,1}- \sigma \Delta B_n )(\Delta W_n)^\top |\mathcal{F}_{t_n}]|=\mathcal{O}(h^{3})$.
        The remainders are $\mathcal{R}_{p,n,1}=\big[ e^{-\upsilon h} - (1-\upsilon h +\frac 12 \upsilon^2 h^2)\big]P^{h,\text{PRK}}_n - h\big[e^{-\frac{\upsilon h}{2}} -(1-\frac{\upsilon h}{2}) \big] f^h(Q^{h,\text{PRK}}_n)$, $\mathcal{R}_{p,n,2}=-\frac h2 \big[ e^{-\upsilon h} -(1-\upsilon h) \big] P^{h,\text{PRK}}_n$, $\mathcal{R}_{q,n,1}=h\big[e^{-\frac{\upsilon h}{2}} -(1-\frac{\upsilon h}{2}) \big]P^{h,\text{PRK}}_n$, and $ \mathcal{R}_{q,n,2}= -\frac{h}{4} ( e^{-\frac{\upsilon h}{2}} -1)P^{h,\text{PRK}}_n$. 
		Comparing \eqref{equ_PhS2}--\eqref{equ_QhS2} with \eqref{eq:generalP}--\eqref{eq:generalQ}, we obtain 
		$ a_1=\frac{\upsilon^2}{2}, a_2=\frac{\upsilon}{2}, a_3=0, a_{0,p,1}=-\frac 12, a_{0,p,2}= \frac{\upsilon}{2}, a_{0,q}=0$, and $b_1= -\frac \upsilon 2, b_{0,p}= -\frac 14, b_2=-\frac 12, b_3=0$. 
        Assumption~\ref{assp:nondeg} holds with $m=1$ and $\xi=(\Delta W_0,\Delta B_0)$, which is non-degenerate, and by \eqref{eq:splittingPRK}, $(P^{h,\text{PRK}}_1,Q^{h,\text{PRK}}_1)$ is an affine function of $\xi$ whose linear part has determinant $\sigma^{2d}(1-\frac{\upsilon h}{2})^d\neq0$ for $h<2/\upsilon$.
        The ergodicity and the weak convergence follow from Theorem~\ref{thm_weakconvergeconformalsym} (\romannumeral 2).
    \end{proof}

    \begin{proof}[Proof of Proposition~\ref{prop_exampleOU}]
        (\romannumeral 1) 
        We have 
		$
			\rd P^{h, \text{S}}_{n+1} \wedge \rd Q^{h, \text{S}}_{n+1}
            =  e^{-\upsilon h}  \rd  P^{h, \text{S}}_n  \wedge \rd Q^{h, \text{S}}_{n},
		$
		which leads to the conformal symplecticity.
   
		It follows that   \begin{align}\label{eq:solutionofsplitting3P}
			P^{h, \text{S}}_{n+1}%
			&=e^{-\upsilon h} P^{h, \text{S}}_n -\frac 12 h e^{-\upsilon h} f^h({Q}^{h, \text{S}}_{n})+\sigma \int^{t_{n+1}}_{t_n} e^{-\upsilon(t_{n+1}-s)} dW(s)-\frac 12 h f^h( {Q}^{h, \text{S}}_{n+1}), 
		\end{align}
		and
		\begin{align}\label{eq:solutionofsplitting3Q}
			Q^{h, \text{S}}_{n+1}
			&=Q^{h, \text{S}}_n+ \frac{1-e^{-\upsilon h}}{\upsilon }P^{h, \text{S}}_n -\frac 12 h \frac{1-e^{-\upsilon h}}{\upsilon } f^h({Q}^{h, \text{S}}_{n}) \nn\\
            &\quad +  \sigma\int^{t_{n+1}}_{t_n} \frac{1-e^{-\upsilon(t_{n+1}-s)}}{\upsilon} dW(s).
		\end{align}
		Applying the Taylor expansion to \eqref{eq:solutionofsplitting3P}--\eqref{eq:solutionofsplitting3Q}  yields
		\begin{align} \label{equ_PhS3}
			P^{h, \text{S}}_{n+1} &= (1-\upsilon h +\frac 12 \upsilon^2 h^2)	P^{h, \text{S}}_{n}-h f^h({Q}^{h, \text{S}}_{n} ) + \frac 12 \upsilon h^2 f^h({Q}^{h, \text{S}}_{n} )+\sigma \Delta W_n +\mathcal{Z}_{p,n,1}  \nn \\
			&\quad +\mathcal{R}_{p,n,1} -\frac 12 h \bar H_{p,n} \big((h-\frac 12 \upsilon h^2) P^{h, \text{S}}_n -\frac 12 h^2 f^h({Q}^{h, \text{S}}_{n} )+\mathcal{Z}_{p,n,2} +\mathcal{R}_{p,n,2} \big), 
		\end{align}
		where $\bar H_{p,n}$ is defined with $\mathcal{Q}_{p,n}=Q^{h,\text{S}}_{n+1}$, and
		\begin{align} \label{equ_QhS3}
			Q^{h, \text{S}}_{n+1} = Q^{h, \text{S}}_n +(h-\frac 12 \upsilon h^2) P^{h, \text{S}}_n -\frac 12 h^2 f^h({Q}^{h, \text{S}}_{n} )+\mathcal{Z}_{q,n,1} +\mathcal{R}_{q,n,1},
		\end{align}
		where  $\mathcal{Z}_{p,n,1}:=\sigma \int^{t_{n+1}}_{t_n} (e^{-\upsilon(t_{n+1}-s)}-1) dW(s)$,  $\mathcal{Z}_{q,n,1} := \sigma\int^{t_{n+1}}_{t_n} \frac{1-e^{-\upsilon(t_{n+1}-s)}}{\upsilon} dW(s)$, and $\mathcal{Z}_{p,n,2}:=\mathcal{Z}_{q,n,1}$.
		Here, $\mathcal{Z}_{p,n,1}$ 
		satisfies $\E[ \mathcal{Z}_{p,n,1} |\mathcal F_{t_n}]=0$, $\E[ |\mathcal{Z}_{p,n,1}|^\alpha |\mathcal F_{t_n}] = \mathcal{O}(h^{3\alpha/2})$ for $\alpha\geq 2$,
		and 
		\begin{align*}
			\E \big[ |\mathcal{Z}_{p,n,1}+\upsilon \sigma \Delta B_n |^2  \big|\mathcal{F}_{t_n} \big]=  d \sigma^2 \int^{t_{n+1}}_{t_n} \big| e^{-\upsilon(t_{n+1}-s)}-1+ \upsilon(t_{n+1}-s) \big|^2 ds  =\mathcal{O}(h^{5}),
		\end{align*}
        which implies $|\E [ (\mathcal{Z}_{p,n,1}+\upsilon \sigma \Delta B_n )(\Delta W_n )^\top |\mathcal{F}_{t_n} ] |= \mathcal O(h^3).$
		The term $\mathcal{Z}_{q,n,1}$ satisfies
		$\E[ \mathcal{Z}_{q,n,1} |\mathcal F_{t_n}]  =0$, $\E[ |\mathcal{Z}_{q,n,1}|^\alpha |\mathcal F_{t_n}] = \mathcal{O}(h^{3\alpha/2})$ for $\alpha\geq 2$, and 
		$	|\E[ (\mathcal{Z}_{q,n,1}- \sigma \Delta B_n )(\Delta W_n)^\top|\mathcal{F}_{t_n}]|=\mathcal{O}(h^{3})$. 
        The remainders are $\mathcal{R}_{p,n,1}=\big[e^{-\upsilon h} -(1-\upsilon h +\frac 12 \upsilon^2 h^2) \big] P^{h, \text{S}}_{n}-\frac{h}{2} \big[e^{-\upsilon h}- (1-\upsilon h)\big] f^h({Q}^{h, \text{S}}_{n} )$, $\mathcal{R}_{p,n,2} =  \big[\frac{1-e^{-\upsilon h}}{\upsilon }-(h-\frac 12 \upsilon h^2) \big]P^{h, \text{S}}_{n} -\frac 12 h(\frac{1-e^{-\upsilon h}}{\upsilon }-h )f^h(Q^{h, \text{S}}_{n})$, and $\mathcal{R}_{q,n,1} =   \big[\frac{1-e^{-\upsilon h}}{\upsilon }-(h-\frac 12 \upsilon h^2) \big]P^{h, \text{S}}_{n} -\frac 12h(\frac{1-e^{-\upsilon h}}{\upsilon }-h )f^h(Q^{h, \text{S}}_{n})$.
		Comparing \eqref{equ_PhS3}--\eqref{equ_QhS3} with \eqref{eq:generalP}--\eqref{eq:generalQ}, we obtain 
		$ a_1=\frac{\upsilon^2}{2}, a_2=\frac{\upsilon}{2}, a_3=0, a_{0,p,1}=-\frac 12, a_{0,p,2}= \frac{\upsilon}{4}, a_{0,q}=\frac 14$, and $b_1= -\frac \upsilon 2, b_{0,p}= 0, b_2=-\frac 12, b_3=0$, where $-\frac12\mathcal{Z}_{p,n,2}$ and $-\frac12\mathcal{R}_{p,n,2}$ play the roles of $\mathcal{Z}_{p,n,2}$ and $\mathcal{R}_{p,n,2}$ in \eqref{eq:generalP}.  
        Assumption~\ref{assp:nondeg} holds with $m=1$ and $\xi=(\xi_1,\xi_2)$ as for the Lie--Trotter scheme: by \eqref{eq:solutionofsplitting3P}--\eqref{eq:solutionofsplitting3Q}, $Q^{h,\text{S}}_1=c(x)+\sigma\xi_2$ and $P^{h,\text{S}}_1=a(x)+\sigma\xi_1-\frac h2f^h(Q^{h,\text{S}}_1)$ with $a$, $c$ depending only on $x$, hence $\Gamma_x$ is a bijection with continuously differentiable inverse and Jacobian determinant $\sigma^{2d}$.
        The ergodicity and the weak convergence follow from Theorem~\ref{thm_weakconvergeconformalsym} (\romannumeral 2).

       (\romannumeral 2)  
       The conformal symplecticity follows from 
		$
			\rd P^{h, \text{S2}}_{n+1} \wedge \rd Q^{h, \text{S2}}_{n+1}
            =  e^{-\upsilon h}  \rd  P^{h, \text{S2}}_n  \wedge \rd Q^{h, \text{S2}}_{n}.
		$

    It follows from the splitting technique that
    \begin{align}\label{equ_PBAOonestep}
      P^{h, \text{S2}}_{n+1} 
      &= e^{-\upsilon h} \big( P_n^{h,\text{S2}} - \frac{h}{2} f^h(Q_n^{h,\text{S2}})\big)+ \sigma \int_{t_n}^{t_{n+1}} e^{-\upsilon (t_{n+1}-s)} dW(s) -\frac{h}{2} f^h(Q_{n+1}^{h,\text{S2}})
    \end{align}
    and
\begin{align}\label{equ_QBAOonestep}
Q_{n+1}^{h,\text{S2}} 
&= Q_n^{h,\text{S2}} + \frac{h}{2} (1+ e^{-\upsilon h}) \big(P_n^{h,\text{S2}} - \frac{h}{2} f^h(Q_n^{h,\text{S2}}) \big)+ \frac{\sigma h}{2} \int_{t_n}^{t_{n+1}} e^{-\upsilon (t_{n+1}-s)} dW(s).  
\end{align}
 Applying the Taylor expansion to \eqref{equ_PBAOonestep}--\eqref{equ_QBAOonestep} gives
 \begin{align} \label{equ_TaylorPOAB}
     P^{h, \text{S2}}_{n+1} &= \big(1-\upsilon h+\frac{\upsilon^2h^2}{2}\big) P_n^{h, \text{S2}} -h \big(1-\frac{\upsilon h}{2} \big) f^h(Q_n^{h, \text{S2}}) +\sigma \Delta W_n +\mathcal Z_{p,n,1} +\mathcal R_{p,n,1} \nn \\
     &\quad -\frac{h}{2} \bar H_{p,n} \big( (h-\frac{1}{2}\upsilon h^2) P_n^{h, \text{S2}} - \frac{1}{2}h^2f^h( Q_n^{h, \text{S2}}) +\mathcal Z_{p,n,2} +\mathcal R_{p,n,2} \big)  
 \end{align}
where $\bar H_{p,n}$ is defined with $\mathcal{Q}_{p,n}=Q^{h,\text{S2}}_{n+1}$, and
\begin{align} \label{equ_TaylorQOAB}
    Q_{n+1}^{h,S2} = Q_n^{h,S2} + (h-\frac 12 \upsilon h^2) P_n^{h,S2} - \frac{1}{2}h^2 f^h(Q_n^{h,S2}) + \mathcal{Z}_{q,n,1} + \mathcal{R}_{q,n,1} 
\end{align}
where $\mathcal{Z}_{p,n,1} := \sigma \int_{t_n}^{t_{n+1}} (e^{-\upsilon (t_{n+1}-s)} - 1) dW(s)$, $\mathcal{Z}_{q,n,1} := \frac{\sigma h}{2} \int_{t_n}^{t_{n+1}} e^{-\upsilon (t_{n+1}-s)} dW(s)$, and $\mathcal{Z}_{p,n,2}:= \mathcal{Z}_{q,n,1} $. Here, $\mathcal{Z}_{p,n,1}$ 
		has the same form as in (\romannumeral 1),
		and $\mathcal{Z}_{q,n,1}$ satisfies
		$\E[ \mathcal{Z}_{q,n,1} |\mathcal F_{t_n}]=0$, $\E[ |\mathcal{Z}_{q,n,1}|^\alpha |\mathcal F_{t_n}] = \mathcal{O}(h^{3\alpha/2})$ for $\alpha\geq 2$, and 
        \begin{align*}
        \big| \E[ (\mathcal{Z}_{q,n,1}- \sigma \Delta B_n )(\Delta W_n)^\top  |\mathcal{F}_{t_n}] \big| = \sigma \int^{t_{n+1}}_{t_n} \big(\frac h2 -(t_{n+1}-s) \big) ds + \mathcal{O}(h^3) = \mathcal{O}(h^3).
        \end{align*}
		 The remainders are $\mathcal{R}_{p,n,1}=\big[e^{-\upsilon h} -(1-\upsilon h +\frac 12 \upsilon^2 h^2) \big] P^{h, \text{S2}}_{n}-\frac{h}{2} \big[e^{-\upsilon h}- (1-\upsilon h)\big] f^h({Q}^{h, \text{S2}}_{n} )$,  $\mathcal{R}_{p,n,2} = \frac{h}{2}\left(e^{-\upsilon h} - 1 + \upsilon h\right) P_n^{h, \text{S2}} - \frac{h^2}{4}\left(e^{-\upsilon h} - 1\right) f^h(Q_n^{h, \text{S2}})$, and \allowbreak $\mathcal{R}_{q,n,1} = \frac{h}{2}\left(e^{-\upsilon h} - 1 + \upsilon h\right) P_n^{h, \text{S2}} - \frac{h^2}{4}\left(e^{-\upsilon h} - 1\right) f^h(Q_n^{h, \text{S2}})$. Combining \eqref{equ_TaylorPOAB}--\eqref{equ_TaylorQOAB} with \eqref{eq:generalP}--\eqref{eq:generalQ}, we obtain 
		$ a_1=\frac{\upsilon^2}{2}, a_2=\frac{\upsilon}{2}, a_3=0, a_{0,p,1}=-\frac 12, a_{0,p,2}= \frac{\upsilon}{4}, a_{0,q}=\frac 14$, and $b_1= -\frac \upsilon 2, b_{0,p}= 0, b_2=-\frac 12, b_3=0$,
        where $-\frac12\mathcal{Z}_{p,n,2}$ and $-\frac12\mathcal{R}_{p,n,2}$ play the roles of $\mathcal{Z}_{p,n,2}$ and $\mathcal{R}_{p,n,2}$ in \eqref{eq:generalP}.  
        It is not difficult to check that Assumption~\ref{assp:nondeg} holds with $m=2$ and $\xi=(\chi_0,\chi_1)$, $\chi_n:=\sigma\int_{t_n}^{t_{n+1}}e^{-\upsilon(t_{n+1}-s)}\,\rd W(s)$.
        The ergodicity and the weak convergence follow from Theorem~\ref{thm_weakconvergeconformalsym} (\romannumeral 2).
	\end{proof}
	
\section{Tamed potential approach}\label{suppsec_tamed}

In the main text, we construct convergent explicit numerical methods for the (kinetic) Langevin equation via the mollified potential approach, preserving both its intrinsic geometric structure and long-time dynamical behavior. The key structural principle is to modify the potential $U$ itself rather than directly taming its gradient $f$ \cite{Jentzen2012AAP, Sabanis2016AAP}. A natural alternative to the mollification is to apply a taming procedure to the potential $U$. In this supplement, we develop this tamed-potential approach and establish its first-order weak convergence uniformly in time. Although the tamed-potential method also preserves the conformal symplectic structure, we identify a substantial obstacle to extending the convergence analysis to higher orders in comparison with the mollified potential approach considered in the main text. 

In this section, we present the construction of the tamed potential approach, and discuss its smoothness and related properties.
\subsection{Construction of the tamed potential approach}
We define the tamed potential function $U^h$ as
\begin{align}\label{Tame_equ_Uh} 
	U^h(q) := \frac{\tilde{U}(q)}{1+h^\alpha|q|^{\beta r_2} }, \quad q\in\RR^d,
\end{align}
where we decompose $U(q) = \tilde{U}(q) + C$ such that $\tilde{U}(0) = 0$. Here, $h \in (0, 1)$ is the step size,  $\beta$ is the taming coefficient, and $\alpha > 0$. The tamed exponent $\beta$ is chosen so that $r_1 - \beta r_2 = 2$, which ensures that the superquadratic growth of $U$ is suppressed to quadratic growth.  The tamed force is therefore defined by 
$f^h:=\nabla U^h$. Specifically, 
given $U\in \mathcal C^\infty(\RR^d; \RR_+)$, for any $q\in \RR^d$ satisfying $|q|\neq0$,  the tamed force $f^h$ takes the form
\begin{align} \label{Tame_equ_fh}
	f^h(q) =\frac{ \nabla U(q) (1+h^\alpha|q|^{\beta r_2})- h^\alpha \beta r_2 |q|^{\beta r_2-2}\tilde{U}(q)  q}{(1+h^\alpha|q|^{\beta r_2})^2},
\end{align}
where we use $ \nabla(|q|^{\beta r_2}) =\beta r_2 |q|^{\beta r_2-2} q$ for $|q|\neq0$.  
At $q=0$, we have $f^h(0)=f(0)$, since $U^h(q)-\tilde U(q)=-h^\alpha|q|^{\beta r_2}\tilde U(q)/(1+h^\alpha|q|^{\beta r_2})=o(|q|)$ as $q\to 0$.

\subsection{Smoothness of the tamed potential approach}\label{subsec_smoothTamed}
Since \eqref{Tame_equ_fh} gives $f^h(q)\to f(0)$ as $q\to 0$, the tamed potential $U^h$ is $\mathcal C^1$. 
We emphasize, however, that an additional assumption on $U$ is required to ensure  higher-order smoothness of $U^h$.
\begin{assp}\label{assm_zeropoint}
	Assume that for $q\to 0$, it holds that $\tilde{U}(q) = \mathcal O (|q|^k)$. %
\end{assp}
Assumption~\ref{assm_zeropoint} requires the lowest-order nonzero term in the expansion of $\tilde U$ at the origin to have degree at least $k$.
Under Assumption~\ref{assm_zeropoint}, 
we have $U^h \in\mathcal C^{k}$ for $k\geq 2$. For example, in view of the expression
\begin{align}\label{Tame_equ_expressD2Uh}
	&\quad D^2 U^h(q) \nn \\
    &= \frac{1}{(1+h^\alpha|q|^{\beta r_2})^3} \Big[ (1+h^\alpha| q|^{\beta r_2})^2 D^2 U(q)- (1+h^\alpha| q|^{\beta r_2}) h^\alpha \beta r_2 | q|^{\beta r_2-2}  \Big(\nabla U(q) q^\top \nn \\
	&\quad + \tilde{U}(q) I+ q \nabla U (q)^\top \Big) + h^\alpha \beta r_2 \tilde{U}(q) | q|^{\beta r_2-4} \Big( (2-\beta r_2) + (\beta r_2+2) h^\alpha| q|^{\beta r_2} \Big) q q^\top  \Big],
\end{align}
verifying the continuity of $D^2U^h$ with respect to $q$ reduces to analyzing the term  $|q|^{\beta r_2 -2} \tilde U(q)$ near the origin.   
As $q \to {0}$,
Assumption~\ref{assm_zeropoint} implies that $\tilde{U}(q)= \mathcal{O}(|q|^2)$. Consequently, $|q|^{\beta r_2 -2} \tilde U(q) =\mathcal{O} (|q|^{\beta r_2})\to 0$ which yields $U^h\in\mathcal C^2$. Moreover, if $r_1>3$, then $\beta r_2=r_1-2>1$ and $U^h\in \mathcal C^2$ without Assumption~\ref{assm_zeropoint}, since the terms of \eqref{Tame_equ_expressD2Uh} that are singular at the origin are then $\mathcal O(|q|^{\beta r_2-1})$ for every $\tilde U$ with $\tilde U (0)=0$. Higher-order regularity can be established by analogous arguments.
\subsection{Properties of the tamed force}
The estimates in Proposition \ref{prop:fh-fS2} for the mollified force can be similarly derived for the tamed force, which are stated in the following proposition.
\begin{prop} \label{prop:fh-f} Under Assumption \ref{assm_CandPoly}, it holds that\\
	(\romannumeral 1)
	$
	| f(q)-f^h(q)| \leq  Ch^\alpha(1+|q|^{(1+\beta)r_2}),
	$ \\
	(\romannumeral 2)
	$
	|f^h(q)| \leq C\min\{ 1+|q|^{r_2},  h^{-\alpha}(1+|q|)\}.
	$\\
	(\romannumeral 3) If $r_1\leq 3$, let Assumption~\ref{assm_zeropoint} hold with $k\geq 2.$ Then $U^h\in \mathcal C^2$ satisfies $\| D^2 U^h (q)\|_{\op}  \leq C\min\{ 1+|q|^{r_2-1},  h^{-\alpha}\}.$
\end{prop}

\begin{proof} 
	(\romannumeral 1)
    At $q=0$, we have $f^h(q)=f(q)$ (see after \eqref{Tame_equ_fh}). For $q\neq 0$,
    the expression  \eqref{Tame_equ_fh} yields
	\begin{align*}
		f(q)-f^h(q) &=  \nabla \tilde U(q) -\frac{ \nabla U(q) (1+h^\alpha |q|^{\beta r_2})-  h^\alpha\beta r_2 |q|^{\beta r_2-2} \tilde U(q)  q}{(1+h^\alpha |q|^{\beta r_2})^2} \\
		&=  \frac{ h^\alpha |q|^{\beta r_2} \nabla U(q)+h^{2\alpha}|q|^{2\beta r_2}\nabla U(q)+  h^\alpha \beta r_2 |q|^{\beta r_2-2} \tilde U(q) q}{(1+h^\alpha |q|^{\beta r_2})^2}, 
	\end{align*}
	which %
	leads to
	\begin{align*}
		| f(q)-f^h(q)| &\leq C\frac{h^\alpha (1+|q|^{(1+\beta)r_2})+ h^{2\alpha}(1+|q|^{(2\beta+1)r_2} )+ h^\alpha (1+|q|^{r_1+\beta r_2-1 })}{(1+h^\alpha|q|^{\beta r_2})^2}.
	\end{align*}
	For $|q|<1$, we have $| f(q)-f^h(q)|\leq Ch^\alpha$ directly.  
	For $|q|\geq 1$, it follows from $r_1=r_2+1$ that
	\begin{align*}
		| f(q)-f^h(q)| &\leq Ch^\alpha \frac{|q|^{(1+\beta) r_2}(1+ h^\alpha |q|^{\beta r_2})}{(1+h^\alpha|q|^{\beta r_2})^2}
		\leq Ch^\alpha (1+|q|^{(1+\beta)r_2}),
	\end{align*}
	which proves (\romannumeral 1).
	
	\noindent(\romannumeral 2) %
	For  $|q|<1$, it immediately follows that
	$	|	f^h(q)| \leq C$. 
	Hence, it suffices to prove the result for $|q|\geq 1$. In this case,  we  obtain %
	\begin{align*}
		|	f^h(q)| %
		&\leq C\frac{(1+|q|^{r_2}) (1+h^\alpha |q|^{\beta r_2})+ h^\alpha (1+|q|^{r_1+\beta r_2 -1})}{(1+h^\alpha |q|^{\beta r_2})^2} \\
		&\leq C \min\Big\{ (1+|q|^{r_2})+  	\frac{h^\alpha(1+|q|^{(1+\beta)r_2})}{1+h^\alpha|q|^{\beta r_2}}, \frac{1+|q|^{r_2}}{1+h^\alpha|q|^{\beta r_2}} + 	\frac{ h^\alpha(1+|q|^{(1+\beta)r_2})}{(1+h^\alpha|q|^{\beta r_2})^2} \Big\}\\
        &\leq C \min\{ 1+|q|^{r_2}, h^{-\alpha}(1+|q|)\} .
	\end{align*}
	
  \noindent(\romannumeral 3) 
  Similar to (\romannumeral 2), it suffices to prove the result for $|q|>1$. 
  It follows from \eqref{Tame_equ_expressD2Uh} that 
    \begin{align*}
  	\|D^2U^h(q)\|_{\op}
  	&\leq
  	\frac{\|D^2U(q)\|_{\op}}
  	{1+h^\alpha|q|^{\beta r_2}} +
  	\frac{Ch^\alpha|q|^{\beta r_2-2}}
  	{\left(1+h^\alpha|q|^{\beta r_2}\right)^2}
  	\big(|q||\nabla U(q)|+|\tilde U(q)|\big) \\
  	&\leq
  	C\frac{1+|q|^{r_2-1}}
  	{1+h^\alpha|q|^{\beta r_2}}
  	+C\frac{h^\alpha
  		\left(1+|q|^{\beta r_2+r_2-1}\right)}
  	{\left(1+h^\alpha|q|^{\beta r_2}\right)^2} 
  	\leq C\min\left\{
  	1+|q|^{r_2-1},\,h^{-\alpha}
  	\right\}.
  \end{align*}
  
Here, the first bound follows from
 $
  	\frac{1+|q|^{r_2-1}}
  	{1+h^\alpha|q|^{\beta r_2}}
  	\leq 1+|q|^{r_2-1}, 
 $
  and 
  $
  	\frac{h^\alpha|q|^{\beta r_2+r_2-1}}
  	{\left(1+h^\alpha|q|^{\beta r_2}\right)^2}
  	=
  	|q|^{r_2-1}
  	\frac{h^\alpha|q|^{\beta r_2}}
  	{\left(1+h^\alpha|q|^{\beta r_2}\right)^2}
  	\leq C|q|^{r_2-1}
  $
while the second bound is derived using $\beta r_2= r_2-1$, 
  $
 	\frac{1+|q|^{r_2-1}}
 	{1+h^\alpha|q|^{\beta r_2}}
 	=\frac{1+|q|^{\beta r_2}}
 	{1+h^\alpha|q|^{\beta r_2}}
 	\leq Ch^{-\alpha},
 $
 and
 $
 	\frac{h^\alpha
 		\left(1+|q|^{\beta r_2+r_2-1}\right)}
 	{\left(1+h^\alpha|q|^{\beta r_2}\right)^2}
 	\leq
 	h^\alpha
 	+
 	h^{-\alpha}
 	\Big(
 	\frac{h^\alpha|q|^{\beta r_2}}
 	{1+h^\alpha|q|^{\beta r_2}}
 	\Big)^2
 	\leq Ch^{-\alpha}.
 $
	The proof is complete.
\end{proof}

\begin{rem}
	For the mollified potential approach, we can achieve $|f(q)-f^h(q)|\leq Ch^k\operatorname{Poly}(q)$ for any $k>0$, while the factor $h^{-1/2}$ in the upper bound for the mollified force $f^h(q)$ remains unchanged. In contrast, for the tamed potential approach, the approximation error is only of order $h^\alpha$, whereas the relevant bounds for the tamed force and its derivatives involve the factor $h^{-\alpha}$. In particular, the uniform-in-time moment estimates for the numerical method require $h\|D^2U^h(\cdot)\|_{\op}\leq C$, restricting us to $\alpha\leq1$. Since the approximation error $f(q)-f^h(q)$ contributes a term of order $h^\alpha$ to the weak error, the resulting weak convergence order is at most one, which prevents a second-order result within this framework.
\end{rem}

The tamed potential also inherits the coercivity condition. However, unlike in the mollified-potential approach, this inheritance requires the following additional assumption on $U$.
\begin{assp} \label{assm_ErgodicUh}
	Under Assumptions \ref{assm_CandPoly} and \ref{assm_ErgodicU}, 
	we further assume that  $c_0>\beta r_2=r_1-2$. 
\end{assp}
\begin{rem}
	Assumption~\ref{assm_ErgodicUh} is satisfied, e.g., 
    for $U(q)=\frac 1k|q|^{k}$ with an even integer $k\geq 4$, for which $L_1=L_2=\frac 1k$, $r_2=k-1$, $\beta r_2=k-2$, and $c_0=k$.
\end{rem}
\begin{prop} \label{prop_ergodicfh} 
	Under %
	Assumption~\ref{assm_ErgodicUh}, the inequality holds
	\begin{align*}
		f^h(q)^\top q \geq c_2 U^h(q)-c_3,
	\end{align*}
	where the constant $c_2$ satisfies $0<c_2\leq c_0-\beta r_2$ %
    and $c_3>0$ is independent of $h$.
\end{prop}
\begin{proof}
Write $x:=h^\alpha|q|^{\beta r_2}$. By \eqref{Tame_equ_Uh}, \eqref{Tame_equ_fh} and Assumption~\ref{assm_ErgodicU}, which gives $\nabla U(q)^\top q\geq c_0\tilde U(q)-c_1$,
\begin{align*}
f^h(q)^\top q-c_2U^h(q)&=\frac{\big(\nabla U(q)^\top q-c_2\tilde U(q)\big)(1+x)-\beta r_2x\tilde U(q)}{(1+x)^2}\\
&\geq\frac{\tilde U(q)\big((c_0-c_2)+(c_0-c_2-\beta r_2)x\big)-c_1(1+x)}{(1+x)^2}.
\end{align*}
Since $c_0-c_2-\beta r_2\geq0$ and $\tilde U\geq-U(0)$, the numerator is at least $-c_3(1+x)$ with $c_3:=U(0)(c_0-c_2)+c_1$, and hence $f^h(q)^\top q-c_2U^h(q)\geq-c_3$.
\end{proof}
\section{Ergodicity and uniform-in-time weak convergence of the tamed approach}\label{suppsec_ergodicTame}
Having introduced the tamed potential $U^h$, the corresponding \textit{tamed Langevin equations} take the same form as \eqref{eq:tamed}, namely,
	\begin{equation}\label{Tame_eq:tamed}
		\left\{\begin{aligned}
			dP^h(t) &= -\upsilon P^h(t)dt-\nabla U^h(Q^h(t)) dt+\sigma dW(t), \,\, P^h(0)=p_0,
			\\
			dQ^h(t) &= P^h(t)dt, \,\, Q^h(0)=q_0.
		\end{aligned} \right.
	\end{equation}
Here, with no risk of confusion, we denote by $(P^h(t), Q^h(t))$ the solution of the tamed Langevin equation \eqref{Tame_eq:tamed} with the tamed potential $U^h$. We also denote by $(P^h_n, Q^h_n)$ the corresponding numerical solution generated by \eqref{eq:generalP}--\eqref{eq:generalQ}.
Since conformal symplecticity of both the tamed Langevin equations and the corresponding numerical methods follows by the same arguments as in the mollified-potential approach, we focus in this section on establishing their ergodicity and the uniform-in-time weak convergence of the numerical methods. 
\subsection{Ergodicity and  uniform-in-time weak convergence of the tamed solution} 
In this section, we show that
the tamed Langevin dynamics also inherit the ergodicity of the original equation. %
The following theorem provides the counterpart of Theorem \ref{thm_weakconvergencetamed} for the  tamed-potential approach. 

\begin{theorem}[Tamed-potential version of Theorem \ref{thm_weakconvergencetamed}] \label{Tamed_thm_weakconvergencetamed}
If $r_1\leq 3$, let Assumption~\ref{assm_zeropoint} hold with $k\geq 2$.\\
(\romannumeral 1) 
		The phase flow of \eqref{Tame_eq:tamed} preserves the conformal symplecticity, i.e., 
		$\rd P^h(t)  \wedge \rd Q^h(t) = e^{-\upsilon t} \rd P(0) \wedge \rd Q(0).$
		\\
        (\romannumeral 2) 
    Let %
	Assumption~\ref{assm_ErgodicUh} hold and let $\alpha >0$. Then, for all sufficiently small $h>0$, the process $\{(P^h(t), Q^h(t))\}_{t\geq 0}$ is ergodic with the unique invariant measure $\pi_h$ being Gibbsian whose density is proportional to $\exp\big(-\frac{2\upsilon}{\sigma^2}\big(\frac12|p|^2+U^h(q)\big)\big)$.
		Furthermore, for any $g\in \mathcal C^\infty(\RR^{2d}; \RR)$ with derivatives growing no faster than polynomials at infinity,  there exists a constant $C>0$ independent of $h$ (but dependent on the initial value $(p_0,q_0)$) such that for sufficiently small $h>0$:
	\begin{align*}
		|\pi(g)-\pi_h(g)| &\leq Ch^\alpha , \\
        \sup_{t\geq0}\big|\E[g(P^h(t),Q^h(t))]-\E[g(P(t),Q(t))]\big| &\leq Ch^\alpha.
	\end{align*}
\end{theorem}
\begin{proof}
(\romannumeral 1) The preservation of the conformal symplectic structure is obvious.

\noindent(\romannumeral 2) We first introduce lower-bound estimates for the tamed potential $U^h$, which replace \eqref{equ_Uhlowerboundleq1} and \eqref{equ_Uhlowerboundgeq1} for the mollified-potential approach.	We divide the derivation into two cases according to the value of $|q|$.
For $q\in \RR^d$ satisfying $h^\alpha |q|^{\beta r_2}\leq 1$, we obtain
\begin{align} \label{Tame_equ_Uhlowerboundleq1}
	U^h(q) \geq \frac{L_1|q|^{r_1}-C}{1+h^\alpha |q|^{\beta r_2}} \geq  \frac{L_1}{2}|q|^{r_1}-C  =:C_{U^h} |q|^{r_1} -C.
\end{align}
For  $q\in \RR^d$ satisfying $h^\alpha |q|^{\beta r_2}> 1$, we have
\begin{align} \label{Tame_equ_Uhlowerboundgeq1}
	U^h(q) \geq  \frac{L_1|q|^{r_1}-C}{1+  h^\alpha |q|^{\beta r_2}} \geq \frac{L_1}{2}h^{-\alpha} |q|^{r_1-\beta r_2}-C =:\tilde{C}_{U^h} h^{-\alpha} |q|^{2}-C,
\end{align}
where we use $r_1 - \beta r_2 = 2$.
With \eqref{equ_Uhlowerboundleq1} and \eqref{equ_Uhlowerboundgeq1} replaced by \eqref{Tame_equ_Uhlowerboundleq1}--\eqref{Tame_equ_Uhlowerboundgeq1}, the uniform-in-time moment bounds established in Proposition~\ref{prop_MomentboundedOfTamedequ} remain valid for the tamed solution $P^h$ and $Q^h$ under Assumption  \ref{assm_ErgodicUh}. We then obtain the second bound by following the same argument as in Theorem~\ref{thm_weakconvergencetamed} (\romannumeral 3). The first bound follows from the second by letting $t\to \infty$, since $\E[g(P^h(t),Q^h(t))]\to\pi_h(g)$ and $\E[g(P(t),Q(t))]\to\pi(g)$ by the geometric ergodicity used in the proof of Theorem~\ref{thm_weakconvergencetamed}.
\end{proof}

\subsection{Ergodicity and  uniform-in-time weak convergence  of  a generic numerical method $(P^h_n, Q^h_n)$ under the tamed approach}
As noted in Proposition~\ref{prop:fh-f}, the estimates for $U^h$ and $f^h$ depend critically on $h^\alpha$. Since the subsequent analysis involves $\|D^2U^h\|_{\op}$, we impose Assumption~\ref{assm_zeropoint} with $k\geq2$ when $r_1\leq 3$ so that $U^h\in\mathcal C ^2$ (see Section~\ref{subsec_smoothTamed}). Then  Proposition~\ref{prop:fh-f}~(\romannumeral 3) shows that there exists a constant $S_{U^h}>0$ independent of $h$ such that
$
\|D^2U^h(\cdot)\|_{\op}\leq S_{U^h}h^{-\alpha}.
$
For the case $\alpha=1$, this dependence requires an additional condition on the damping coefficient to establish numerical ergodicity, as specified below. 

\begin{assp} \label{assm_v>CS} Assume that $\upsilon> (\frac 12 + |a_{0,p,1}|) S_{U^h}$.
\end{assp}
As an example, 
for $U(q)=\frac 16|q|^6-\frac 14|q|^4+1$, we have $S_{U^h}\approx 0.98$. %
Then Assumption~\ref{assm_v>CS}  requires $\upsilon >  0.98 (\frac 12 + |a_{0,p,1}|)$. 
	\begin{theorem}[Tamed-potential version of Theorem \ref{thm_weakconvergeconformalsym}] \label{Tamed_thm_weakconvergeconformalsym}
	Let Assumptions \ref{assp:nondeg} %
    and \ref{assm_ErgodicUh} hold, and if $r_1\leq 3$, let Assumption \ref{assm_zeropoint} hold with $k\geq2$. If $\alpha=1$, additionally let  Assumption~\ref{assm_v>CS} hold and let $b_2 \leq -\frac12$. Then for sufficiently small $h>0$, $\{(P^h_n, Q^h_n)\}_{n\in \N}$ under the tamed potential approach with the index $\alpha\in(0,1]$ is ergodic with a unique invariant measure $\pi_{h, \text{num}}$. 
	Furthermore, let the scheme start from $(P^h_0,Q^h_0)=(P(0),Q(0))$ and let $g\in \mathcal{C}^\infty(\RR^{2d};\RR)$ have all its derivatives growing not faster than polynomials at infinity. Then, there exists a constant $C>0$ independent of $h$ such that
	\begin{align*}
    \sup_{n\geq0}\big|\E[g(P^h_n,Q^h_n)]-\E[g(P(t_n),Q(t_n))]\big|&\leq Ch^\alpha,\\
		|\pi(g)-\pi_{h,\text{num}} (g)| &\leq Ch^\alpha.
	\end{align*}
\end{theorem}

\begin{proof} 
We reduce the proof of Theorem~\ref{Tamed_thm_weakconvergeconformalsym} to establishing the corresponding uniform-in-time moment bounds (Tamed-potential version  of Proposition \ref{prop_momentboundedgeneral}) for the tamed potential approach. Once these bounds are available, the only difference in the weak convergence analysis arises from the approximation estimate in Proposition~\ref{prop:fh-f}~(\romannumeral 1). Indeed, following the proof of Theorem~\ref{thm_weakconvergeconformalsym}, the estimate $|f-f^h|\leq Ch^\alpha(1+|q|^{(1+\beta)r_2})$ yields a weak-error bound of order $h^\alpha$.

The proof of the moment bounds for the numerical solution $(P_n^h,Q_n^h)$ follows the same argument as that of Proposition \ref{prop_momentboundedgeneral}. We therefore only highlight the modifications required when $\alpha=1$.
The proof still relies on the Lyapunov function $V^h$ defined in \eqref{equ_defVh}, and the identity \eqref{equ_generalmoment1} remains valid. 
	
	For the first term in \eqref{equ_generalmoment1}, we derive
	\begin{align}\label{Tame_equ_momentgeneralP2}
		\frac 12\E[	|P^{h}_{n+1}|^2 |\mathcal{F}_{t_n}]  
		&\leq  \Big(\frac 12+(-\upsilon+|a_{0,p,1}|S_{U^h} ) h +
		C_{p^2,h^2,1} h^2 \Big) |P^h_n |^2-h ( P^h_n)^\top f^h(Q^h_n)\nn \\
		&\quad +C_{pq,h^2,1} h^2 |P^h_n||f^h(Q^h_n)|
        + \frac 12 h^2 |f^h(Q^h_n)|^2+ Ch^3 |f^h(Q^h_n)|^2+ C h
	\end{align}
	with positive constants $	C_{p^2,h^2,1}$ and $C_{pq,h^2,1}$. 
	Similarly,  for the third term on the right-hand side of \eqref{equ_generalmoment1}, it holds that 
	\begin{align}\label{Tame_equ_momentgeneralUq}
		&\quad  \E[  (f^h(Q^h_n))^\top (Q^h_{n+1}- Q^h_n) | \mathcal{F}_{t_n} ]  \nn \\
		&\leq  h (P^h_n)^\top  f^h(Q^h_n)  + C_{pq,h^2,2} h^2 |P^h_n||f^h(Q^h_n)| +b_2 h^2 |f^h(Q^h_n)|^2 \nn \\
        &\quad +Ch^2|P^h_n|^2+Ch^3|f^h(Q^h_n)|^2 +Ch,
	\end{align}
    where the constant $C_{pq,h^2,2}>0$.
			The fourth term on the right-hand side of \eqref{equ_generalmoment1} has the following upper bound
	\begin{align}\label{Tame_equ_momentgeneralQ2}
		&\quad \frac 12 \E[ (Q^{h}_{n+1}-Q^{h}_n)^\top D^2 U^h(\xi_1) (Q^{h}_{n+1}-Q^{h}_n)  | \mathcal{F}_{t_n}  ] \nn \\
		&\leq \frac 12 S_{U^h} \big[ h| P^h_n|^2+ 2|b_2|h^2 |P^h_n| |  f^h(Q^h_n)| + C_{p^2,h^2,2}h^2|P^h_n|^2+Ch^3|f^h(Q^h_n)|^2+Ch\big]
	\end{align}
    with a positive constant $C_{p^2,h^2,2}$.
	The last term can be estimated as 
	\begin{align}\label{Tame_equ_momentgeneralPQ}
		&\quad \upsilon \alpha_1 \E[	( P^{h}_{n+1})^\top  Q^{h}_{n+1}  | \mathcal{F}_{t_n}  ]\nn \\
		&\leq \big(\upsilon \alpha_1h+C_{p^2,h^2,2}h^2 \big) |P^h_n|^2  +\upsilon \alpha_1 (P^h_n)^\top  Q^{h}_{n} +  C_{pq,h} h|P^h_n| |Q^{h}_{n} | \nn \\
		&\quad + \upsilon \alpha_1\big(- h +a_2h^2\big) (f^h(Q^h_n))^\top Q^h_n  + \upsilon \alpha_1a_{0,q}h^3 (f^h(Q^h_n))^\top\E[\bar H_{p,n}  |\mathcal{F}_{t_n}]   Q^h_n    \nn\\
        &\quad+\upsilon \alpha_1(b_2-1) h^2 (P^h_n)^\top  f^h(Q^{h}_{n}) +C_{q^2,h^2} h^2|Q^h_n|^2+Ch^2 |P^h_n|^2+Ch^3 |f^h(Q^h_n)|^2+Ch
	\end{align}
	with  positive constants $C_{pq, h}$ and $C_{q^2, h^2}$.
	Combining  \eqref{Tame_equ_momentgeneralP2}--\eqref{Tame_equ_momentgeneralPQ} with  the inequality $h\|D^2 U^h (\cdot)\|_{\op} \leq S_{U^h}$, for sufficiently small $h>0$, we obtain
\begin{align*}
	&\quad 	\E\big[	{V}^h (P^h_{n+1},Q^h_{n+1}) \big| \mathcal{F}_{t_n} \big] \\
	&\leq (1-\kappa h) {V}^h (P^h_{n},Q^h_{n}) + \Big( \Big[\frac 12 \kappa +\big(\frac 12+|a_{0,p,1}| \big)S_{U^h}+\upsilon\alpha_1-\upsilon  \Big]h +
	C_{p^2, h^2} h^2\Big) |P^h_n|^2\\
	&\quad  +C_{pq,h^2}h^2
	|P^h_n| |f^h(Q^h_n)| + \upsilon\alpha_1\kappa h (P^h_n)^\top Q^h_n
	+ C_{pq, h}h|P^h_n| |Q^h_n| \\
	&\quad    +\big( \frac 12 +b_2\big)h^2 |f^h(Q^h_n)|^2 +C_{q^2,h^2} h^2|Q^h_n|^2+ \upsilon \alpha_1\big(- h +a_2h^2\big) (f^h(Q^h_n))^\top Q^h_n  \\
	&\quad + \upsilon \alpha_1a_{0,q}h^3 (f^h(Q^h_n))^\top\E[\bar H_{p,n}  |\mathcal{F}_{t_n}]   Q^h_n +\kappa h U^h(Q^h_n)+ C_{f^2,h^3} h^3 |f^h(Q^h_n)|^2+ Ch
\end{align*}
where $C_{p^2, h^2}, C_{pq, h^2}, C_{pq, h}, C_{q^2,h^2}, C_{f^2,h^3}, \kappa$ are positive constants. 
Since $h\|\bar H_{p,n}\|_{\op}\leq S_{U^h}$, we have \begin{align*} \upsilon \alpha_1|a_{0,q}|h^3 \big|(f^h(Q^h_n))^\top\E[\bar H_{p,n} |\mathcal{F}_{t_n}] Q^h_n\big| &\leq \upsilon \alpha_1|a_{0,q}|S_{U^h}h^2 |f^h(Q^h_n)||Q^h_n| \\
&=:C_{fq,h^2}h^2 |f^h(Q^h_n)||Q^h_n|. \end{align*} Using Proposition~\ref{prop_ergodicfh} and denoting $A_h:=-\frac 12\kappa-(\frac 12+|a_{0,p,1}|)S_{U^h}-\upsilon\alpha_1+\upsilon-C_{p^2,h^{2}}h $, $B_h:=C_ {pq,h^2}$, $C_h:=\frac 12 +b_2 $, $D_h:=\upsilon\alpha_1\kappa + C_{pq,h}$, and $E_h:=-\kappa+\upsilon\alpha_1 c_2-\upsilon\alpha_1c_2a_2h$,
we arrive at
		\begin{align*}
			\E\big[	{V}^h (P^h_{n+1},Q^h_{n+1}) \big| \mathcal{F}_{t_n} \big] \leq  (1-\kappa h) {V}^h (P^h_{n},Q^h_{n})+\tilde{\mathcal{R}}^h(P^h_n,Q^h_n),
		\end{align*}
    where
    \begin{align*}
        \tilde{\mathcal{R}}^h(p,q)&:= - A_h h |p|^2
			+ B_h h^2 |p| |f^h(q)| +C_{h} h^2 |f^h(q)|^2  +D_h h |p| |q| -E_h h U^h(q)\\
            &\quad+C_{q^2,h^2} h^2 |q|^2 +C_{fq,h^2}h^2|f^h(q)||q|+   C_{f^2,h^3}h^3|f^h(q)|^2+Ch.
    \end{align*}
We claim that under Assumption~\ref{assm_v>CS} and the condition $b_2\leq -\frac 12$, for sufficiently small $h>0$, there exists a constant $C>0$ such that
\begin{align}\label{Tame_equ_numericalboundClaim}
    \frac{\tilde{\mathcal {R}}^h(p,q)}{h}\leq C \quad \text{for all } (p,q)\in\mathbb{R}^{2d}.
\end{align} 
We emphasize here that these two additional requirements, compared with the mollified case, arise from the estimates for $f^h$ in Proposition~\ref{prop:fh-f} and the lower bound estimates for $U^h$ in \eqref{Tame_equ_Uhlowerboundleq1}--\eqref{Tame_equ_Uhlowerboundgeq1}, which affect how we absorb the positive terms. 

Under Assumption~\ref{assm_v>CS}, by taking $$\kappa\in\Big(0,\frac{c_2[2\upsilon -(1+2|a_{0,p,1}| ) S_{U^h}]}{2+c_2} \Big),$$
		there exists $\alpha_1\in(0,1)$ such that 
		\begin{align*} %
			\frac{\kappa}{c_2 \upsilon}< \alpha_1 <  1 -\frac {1}{2\upsilon}[\kappa+(1+2|a_{0,p,1}|) S_{U^h}], 
		\end{align*}
which guarantees $A_h>0$ and $E_h> 0$ for sufficiently small $h>0$. Furthermore,  $C_h \leq 0$ follows immediately from the assumption that $b_2\leq- \frac 12$. 

It remains to absorb the cross terms $B_h h^2|p||f^h(q)|$,  $D_h h |p| |q|$,  $C_{fq,h^2}h^2|f^h(q)||q|$, and the square terms $C_{q^2,h^2} h^2|q|^2$, $C_{f^2,h^3} h^3|f^h(q)|^2$. We decompose the estimates into two regions determined by $|q|$ below.
		
\textbf{Region i: $h|q|^{\beta r_2} \leq 1$.} In this case, it follows from Proposition~\ref{prop:fh-f} (\romannumeral 2) that $|f^h(q)|\leq C(1+|q|^{r_2})$.  Since $ |q|^{2r_2}= |q|^{\beta r_2} |q|^{(2-\beta)r_2} \leq h^{-1} |q|^{r_1}$ for  $h|q|^{\beta r_2} \leq 1$, we obtain $|f^h(q)|^{2} \leq Ch^{-1} |q|^{r_1} +C$. The $\varepsilon$-Young inequality yields 
\begin{align*} %
    B_h h^2|p||f^h(q)| %
    &\leq  \frac{B_h}{2} h^{3/2} |p|^2 + C\frac{B_h}{2} h^{3/2} |q|^{r_1} +Ch^{5/2},
\end{align*}
\begin{align*} %
    D_h h |p| |q| \leq \frac 12 A_h h  |p|^2 + C({A_h,D_h}) h |q|^2 \leq  \frac 12 A_h h  |p|^2 + \frac 12 E_h C_{U^h} h|q|^{r_1}+Ch,
\end{align*}
$
    C_{f^2,h^3}h^3 |f^h(q)|^2 \leq  Ch^{2} |q|^{r_1} +Ch^3,
$
$C_{fq,h^2}h^2|f^h(q)||q|\leq Ch^2|q|^{r_1}+Ch^2$ and $C_{q^2,h^2} h^2 \allowbreak|q|^2 \leq h^2 |q|^{r_1} +Ch^2.$
Using \eqref{Tame_equ_Uhlowerboundleq1},
we arrive at
\begin{align*}
    \frac{\tilde{\mathcal R}^h(p,q)}{h}&\leq \Big( -A_h   +\frac 12 A_h +\frac{B_h}{2} h^{1/2} \Big) |p|^2+C_h h|f^h(q)|^2+C \\
    &\quad +\Big(- E_h C_{U^h}+\frac 12 E_h C_{U^h} +  C\frac{B_h}{2} h^{1/2}+Ch \Big)  |q|^{r_1} \\
    &\leq \Big(-\frac 12A_h  +\frac{B_h}{2} h^{1/2} \Big) |p|^2  + \Big( -\frac 12 E_h C_{U^h}+ C\frac{B_h}{2}h^{1/2} +Ch\Big)  |q|^{r_1}+C.
\end{align*}
For sufficiently small  $h>0$, it holds that 
$-\frac 12 A_h+\frac{B_h }{2}h^{1/2}<0$ and $-\frac 12 E_h C_{U^h}+ C\frac{B_h}{2}h^{1/2} +Ch<0$. Therefore,
we obtain
$
\frac{\tilde{\mathcal R}^h(p,q)}{h} \leq  C.
$

{\textbf{Region \romannumeral 2: $h|q|^{\beta r_2} > 1$.}}
In this case, Proposition~\ref{prop:fh-f} (\romannumeral 2) implies $|f^h(q)|\leq Ch^{-1} |q|$. We apply the $\varepsilon$-Young inequality to derive
$$	B_h h^2|p||f^h(q)| \leq \frac{B_h}{2h^{1/2}} h^2 |p|^2 +\frac{B_h h^{1/2}}{2} h^2 |f^h(q)|^{2}\leq \frac{B_h}{2} h^{3/2} |p|^2 +C\frac{B_h}{2} h^{1/2} |q|^2,
$$ 
$
    D_h h |p| |q| \leq \frac{1}{2} A_h h |p|^2 +C(A_h,D_h) h |q|^{2},
$
and
$
    C_{f^2,h^3}h^3 |f^h(q)|^2 \leq CC_{f^2,h^3}h |q|^2, \,\,\, C_{fq,h^2} \allowbreak h^2|f^h(q)||q|\leq Ch| q|^2.
$
Using  the lower bound estimate \eqref{Tame_equ_Uhlowerboundgeq1}, we obtain
\begin{align*}
    \frac{\tilde{\mathcal R}^h (p,q)}{h}
    &\leq  \Big(-A_h+\frac 12 A_h +\frac{B_h}{2}h^{1/2} \Big) |p|^2  + \Big( -E_h \tilde{C}_{U^h}+\frac{CB_h}{2} h^{1/2}+ C(A_h,D_h)h
    \\
    &\quad +C_{q^2,h^2} h^2 +CC_{f^2,h^3}h + Ch\Big)  h^{-1} |q|^{2} + C,
\end{align*}
which implies $\frac{\tilde{\mathcal{R}}^h(p,q)}{h}\leq C$ for sufficiently small $h>0$.

{Regions \romannumeral 1}--{ \romannumeral 2} exhaust $(p,q)\in \mathbb{R}^{2d}$, which proves \eqref{Tame_equ_numericalboundClaim} and leads to the upper bound
$\sup_{n\in\N} \E[ {V}^h(P^h_n,Q^h_n)] \leq C$.
The remainder of the proof proceeds analogously to the proof of Proposition~\ref{prop_momentboundedgeneral} and is therefore omitted.
\end{proof}

\end{document}